\documentclass[10pt]{amsart} 
\newtheorem{thm}{Theorem}[section]
\newtheorem{lem}[thm]{Lemma}
\newtheorem{prop}[thm]{Proposition}
\newtheorem{corr}[thm]{Corollary}
\newtheorem{Def}[thm]{Definition}

\newcommand{\thinqed}  {{\hfill \rule{1mm}{3mm}\vspace{2mm}}}

\newcommand{\frameqed} {{\hfill \framebox[2mm]{} \vspace{6mm}}}
\renewenvironment{proof}{\vspace{3mm}\noindent {\bf Proof.} }{\frameqed}

\newtheorem{rmk}[thm]{Remark}

\newcommand{\ra}{\rightarrow}              
\newcommand{\lra}{\longrightarrow}         

\newcommand{\rah}{\hookrightarrow}         

\renewcommand{\mid}{{\ \ \big| \ \ }}

\newcommand{\cal}{\mathcal}
\renewcommand{\Bbb}{\mathbb}

\newcommand{\w}{\omega}

\newcommand{\bZ}{{\Bbb Z}}      
\newcommand{\bQ}{{\Bbb Q}}      
\newcommand{\bR}{{\Bbb R}}      
\newcommand{\bC}{{\Bbb C}}      
\newcommand{\bA}{{\Bbb A}}      

\newcommand{\cM}{{\cal M}}

\newcommand{\cO}{{\cal O}}

\newcommand{\cH}{{\cal H}}

\newcommand{\fa}{{\mathfrak a}}
\newcommand{\fb}{{\mathfrak b}}
\newcommand{\fc}{{\mathfrak c}}
\newcommand{\fd}{{\mathfrak d}}

\renewcommand{\Im}{{\rm Im}}

\numberwithin{equation}{subsection}
\begin{document}

\title{Notes On Hilbert's 12th Problem}
\footnote{This paper is dedicated to my late friend Huang Yu, who was killed in an auto accident
on December 24, 2004. Huang Yu was the one who first drew my attention to the Hilbert 12th
problem. My own investigation of the matter, as outlined in this note, has been largely
stimulated from numerous discussions I had with him. I write this note for the memory
of these happy times.}

\author{Sixin Zeng}
\email{sixinz@yahoo.com}
\date{April,2006}

\begin{abstract}
In this note we will study the Hilbert's 12th problem for a primitive CM
field, and the corresponding Stark's conjectures. Using the idea of
``Mirror Symmetry'', we will show how to generate all the class fields
of a given primitive CM field, thus complete the work of Shimura-Taniyama-Weil.
\end{abstract}
\maketitle


\section*{\bf Introduction}
\subsection{}
Let $K$ be a number field, $H_K$ be the ideal class group of $K$, and let
$K_0$ be the Hilbert class field of $K$. The class field theory tells us
there is a canonical isomorphism $Gal(K_0/K) \simeq H_K$. In general given 
any integral ideal $\fd$, let $K_{\fd}$ be the maximal abelian extension
of $K$ unramified outside $\fd$, and let $H_{\fd}$ be the generalized
ideal class group relative to $\fd$, then a similar isomorphism holds
as well: $Gal(K_{\fd}/K) \simeq H_{\fd}$. On the other hand, Hilbert's
12th problem asking for an explicit generation of all the abelian extension
of $K$, more precisely it is asking for finding a special transcendental
function, whose values at some special points would generate all the abelian
extension of $K$.

When $K$ is the rational field $\bQ$, the transcendental function is the 
exponential function and the special points are the division points on the
unit circle. Indeed by the classical Kronecker-Weber theorem, all the
abelian extension of $\bQ$ can be obtained by adding the roots of unity to $\bQ$.

When  $K$ is an imaginery quadratic field, this problem is answered by the 
theory of complex multiplications. As this theory has very much influenced
the later thinking about this problem, let's recall some details. So let
$K$ be such a field, consider the set of elliptic curves $E$ satisfying
$End(E)\otimes \bQ \simeq K$, i.e., the set of elliptic curves with complex multiplications
by $K$. All such elliptic curves can be constructed by the following way:
take the caononical embedding $\iota: K \lra \bC$, let $\fa$ be an integral
ideal of $K$, then $\fa$ is a rank 2 $\bZ$-module in $\bC$, i.e., $\iota(\fa)$ is
a lattice in $\bC$, so we can take the quotient $\bC / \iota(\fa)$ which is an
elliptic curve with complex multiplication by $K$. From the construction it is
clear that the ideal group of $K$ operates on this set of elliptic
curves, indeed for any representative $\fb$ of an ideal class, we have the
isogeny: $\bC /\iota(\fa) \ra \bC/ \iota(\fa \fb)$. On the other hand all such
elliptic curves are defined over some number field, and they are represented
by the moduli points on the moduli space of all elliptic curves $M_1=\cH / PSL_2(\bZ)$.
In particular the Galois group acts naturally on the moduli points of this set of elliptic
curves. So if $E=\bC/\iota(\fa)$, let $\wp$ be a prime of $K$, $E(\wp)=\bC /\iota(\fa \wp)$,
$s(\wp) \in H_K \simeq Gal(K_0/K)$ the Artin symbol, then the main result
of the complex multiplication is $E^{s(\wp)} \simeq E(\wp)$. In more concrete
terms if $j$ is the classical modular function, and let $p_E$ be the moduli
point of $E$ on $M_1=\cH/PSL_2(\bZ)$, then $j(p_E)$ is an algebraic number and
generates the Hilbert class field of $K$. Similar statements hold for all the
ray class fields of $K$.

\vspace{1mm}
The idea behind the theory of complex multiplication can be summerized
as the following:

\begin{itemize}
\item
Given the number field $K$, in order to solve the Hilbert's 12th problem, first
we need to find a suitable class of algebro-geometric objects $X$, say varieties;

\item
$X$ should be closely related to $K$ such that the ideal classes of $K$ can 
naturally act on them;

\item
All such $X$ should be defined over some number field, all such $X$ should be living 
on some natural moduli space so their field of
moduli are given by the coordinate function evaluated at the moduli points.
The Galois group can naturally acts on them;

\item
Moreover the action of the ideal classes of $K$ and the action of Galois group
on the moduli points should be related by the reciprocity law of class field theory and  
the Kronecker  congruence relations;

\end{itemize}

If we can find such a class of $X$,
then the field of moduli of $X$ give the answer to Hilbert's 12th problem.

For a general number field besides  imaginery quadratic, this philosophy is
difficult to apply. So from now on we will concentrate on a special case.
From now on let $K$ be a primitive CM field of degree $2n$, 
$n \geq 1$. $K$
is then an imaginery quadratic extension of a totally real field $F$, with
$[F: \bQ]=n$. Further let $K^*$ be the reflexive field of $K$, $F^*$ be the
maximal totally real subfield of $K^*$, then $K^*$ is also a primitive CM
field of degree $2n$, $[K^* :F^*]=2$. Let $\{\sigma_1,\sigma_2,\cdots,\sigma_n \}$
be the set of archemidean primes of $F$, lifted to $K$ as the CM types, and let
$\rho$ be the complex conjugate.

When $n > 1$, according to the above philosophy, we naturally
consider the set of abelian varieties of dimension $n$
with complex multiplication by $K$ with CM type $\{\sigma_1,\sigma_2,\cdots,\sigma_n \}$ ,
 i.e., those abelian varieties satisfying
$End(X) \otimes \bQ \simeq K$, and moreover for any $\alpha \in K$ the action
of $\alpha$ on the space of holomorphic differentials of $X$ is given
as the diagonal action of $diag(\alpha^{\sigma_1},\cdots,\alpha^{\sigma_n})$.
All such $X$ can be constructed in the following way: let $\iota: K \ra \bC^n$
be the embedding $\iota(x)=(x^{\sigma_1},\cdots,x^{\sigma_n})$, then for
any integral ideal $\fa$, $\iota(\fa) \subset \bC^n$ as a $\bZ$-lattice
is of rank $2n$, hence the quotient $X=\bC^n / \iota(\fa)$ is an abelian
variety, one see that such $X$ satisfying the above conditions.

From the construction  such set of abelian varieties are
naturally associated to $K$, with the integral ideals of $K$ act on them as 
isogenies. This is exactly the same as elliptic curves.

So next we shall consider the moduli space for these $X$ and their arithmetic
properties, it is here some complication arised, let's be careful.

Let $\cO_K$ and $\cO_F$ be the rings of integers of $K$ and $F$, and let
$\fa$ be an integral ideal of $\cO_K$. Regarding $\fa$ as a module
of $\cO_F$, it is of rank 2. More precisely we have the following(see for
example, Yoshita's book \cite{Yoshita}):

\begin{lem}
For any ideal $\fa$ of $\cO_K$ we have the isomorphism as $\cO_F$-module:
$\fa \simeq \fb \cdot \w_1 \oplus \cO_F \cdot \w_2$, where $\w_1,\w_2 \in K$,
$\fb$ a fractional ideal of $\cO_F$. Moreover the ideal class of $\fb$ is
$\fc \cdot N_{K/F}(\fa)=\fc \cdot \fa \fa^{\rho}$, where $\fc$ is a fractional
ideal of $F$, independent of $\fa$, and $\fc^2 = D_{K/F}$.
\end{lem}

In particular $\cO_K=\fc \cdot \w_1 \oplus \cO_F \cdot \w_2$, and if $\fa$
satisfying $\fa \cdot \fa^{\rho}=(\mu)$, then 
$\fa = \fc \cdot \w_1 \oplus \cO_F \cdot \w_2$.

Now let $X=\bC^n / \iota(\fa)$, by Shimura the polarization
of $X$ is given by the Riemann form $E(x,y)$ which is $E(x,y)=Tr_{K/{\bQ}}(\zeta x y^{\rho})$,
with $\zeta \in K$ satisfying $\zeta^{\rho}=-\zeta$, $\Im(\zeta^{\sigma_i}) >0$, any $i$.
The ``type'' of this polarization, by definition,  is the ideal class of 
$\zeta \fd_{K/F} \fa \fa^{\rho}$, 
where $\fd_{K/F}$ is the relative different of $K$ over $F$. From the above lemma
it is clear that $\fb$ is the type of $X$.

For any moduli problem, in order to have a good moduli space we need to
fix a polarization of $X$. For our class of $X$, we observe that they all have
a large ample cone, in fact the dimension of the ample cone are all of $n$,
hence   it makes sense to consider the moduli problem with all the ample
classes fixed. This is the moduli space of ``type polarized abelian varieties''.

Concretely for any fractional ideal $\fb$ of $\cO_F$, we can constructed
a moduli space $M_n(\fb)$ which parametrizes families of abelian varieties
with the form $\bC^n /{ v_n \cdot\iota(\fb) \oplus \iota(\cO_F)}$, where
$v_n$ is a vector in the product of upper half plane $\cH^n$, and the dot product is 
defined component-wise. It is
well-known that $M_n(\fb)=\cH^n /\Gamma(\fb)$ with 
$\Gamma(\fb)=\{\alpha \in SL_2(\cO_F) | \alpha \equiv 1 \mod (\fb) \}$.
If $X$ is the above abelian variety of type $\fb$, the 
$X=\bC^n / \iota(\fa)$ with $\fa = \fb \cdot \w_1 \oplus \cO_F \cdot \w_2$,
up to isomorphism we can choose $\w_1,\w_2 \in K$ such that $\w_X=\w_1/\w_2$
satisfying $\Im(\w_X^{\sigma_i}) >0$ for any $i$, i.e. 
$(\w_X^{\sigma_1},\cdots,\w_X^{\sigma_n})$ is a vector in $\cH^n$, hence
defines a moduli point of $X$ in $M_n(\fb)$.

So in summary given the class of CM abelian varieties with a fixed CM type,
we can construct a natural moduli space $M_n(\fb)$ for a fixed type of $X$,
i.e., if $X_1$, $X_2$ are of the same type, they live on the same moduli
space, but different types can produce different moduli spaces. So
although the ideal classes of $K$ can act on all of the $X$, the Galois
group can only act on the subclass of $X$ with the same type.

On the other hand it is easy to show that all $X$ are defined over a 
finite extension of the reflex field $K^*$. In fact for an ideal
$\fa^*$ of $K^*$, we define $\fa=g(\fa^*)=\prod_{\tau_i} {\fa^*}^{\tau_i}$,
where $\{\tau_1,\cdots,\tau_n \}$ is a CM type of $K^*$, then $g(\fa^*)$
acting on $X$ will not change the type of $X$. It is from this point
of view that in 1955, Shimura-Taniyama-Weil(see Shimura's book \cite{S1}) established that the
field of moduli of $X$ will generate part of the class field of $K^*$,
precisely it is the class field corresponding to the subgroup of
the ideal class group $H^0_{K^*}=\{\fa \in H_{K^*} |g(\fa)=(\mu),\mu \in K\}$.

\subsection{}

In this note  we will try to extend the work of Shimura-Taniyama-Weil
to cover all the class fields of $K$. There are actually two problems
here. First, the appearance of the reflex field $K^*$ is quite inconvenient,
since our abelian varieties $X$, their moduli spaces $M_n(\fa)$,
and the CM points on the moduli are all constructed naturally from the
given CM field $K$, it is natural for us to look for the invariants from
these abelian varieties that directly generate the class fields of $K$, instead
of the reflex $K^*$. Second and more important problem, is how can we
deal with the class fields that corresponding to the isogenies of
abelian varieties that live on the different type of moduli spaces.

The first problem can be solved in the following way. Since
the fields $K$ and $K^*$ are reflex to each other, by Shimura-Taniyama-Weil theory,
the class fields of $K$ is generated by the abelian varieties associated to $K^*$. So
the problem is to find invariants on $X$ that somehow related to Shimura varieties
associated to $K^*$. For this purpose I find the following notion of ``cone polarized
Hodge Structure'' quite useful.

First we observed the for all the abelian varieties of the CM type
the Kahler cone is very large. Indeed as the polarizations are
determined by the Riemann form $E(x,y)=Tr_{K/{\bQ}}(\zeta x y^{\rho})$,
so the polarizations are determined by $\zeta \in K$ such that
$\zeta^{\rho}=-\zeta$, $\Im(\zeta^{\sigma_i}) >0$, any $i$.
Such $\zeta$ form a cone of dimension $n$,
hence the Kahler cone is of $n$ dimensional. We can also see this by the
fact that since $End(X)$ is an order in $K$, the automorphism group
of $X$ is then the unit group of the order, which is of the form
${\bZ}^{n-1} \oplus Torsion$. So the Kahler cone of $X$ is necessarily $n$
dimensional.

Given such Kahler cone $C_X$ of $X$
we can consider the primitive classes of $X$ in the middle dimension
relateive to this cone, i.e., those classes of dimension $n$ such
that are annilated by any element in the Kahler cone.

$$\{\alpha \in H^n(X,\bC) \mid \alpha \cdot x=0, \forall x \in C_X \}$$

This is the generalization of the notion ``transcendental lattice'' in the
theory of $K3$ surface. 
The Kahler cone
and these primitive classes relative to this cone should be considered
as the most basic Hodge theoretical invariants of our abelian varieties.
So to understand these abelian varieties, 
we need to construct the 
appropriate moduli spaces for these invariants.

First let's fix the Kahler cone and consider the 
primitive classes of the middle degree relative to this cone. Such primitive
classes carries a natural Hodge structure, indeed for abelian varieties with $F$-multiplication
the primitive classes can be characterized as the invariant classes of the automorphism
group $U_F$, it carries a Hodge structure of weight $n$, with the $i$-th Hodge number to be
${i \choose n}=\frac{n!}{i! (n-i)!}$. 

Let's consider the
classifying space of such Hodge structures. As is well-known, for the Hodge structures of weight
higher than 2, the classifying space of Hodge structures is usually
much bigger than the geometric moduli space. Standard examples are the
Calabi-Yau threefolds. On the other hand, in our situation, when fix the Kahler cone 
$C=C_X$, we can define the notion of $C$ polarized abelian varieties. These
are the abelian varieties $X^{\prime}$ such that we have a natural
inclusion $C \ra C_{X^{\prime}}$ into the Kahler cone of $X^{\prime}$.
This notion is a generalization of the lattice polarized K3 surfaces in
the theory of K3 (I think this is due to Dolgachev \cite{Dol})

We can consider the moduli space of cone $C$ polarized abelian varieties.
Let's denote it as $M_C$.
This is well-defined, and it turns out to be a finite quotient of $M_n(\fa)$, i.e., fix
the Kahler cone is the same as fix the type of $X$, all the $X^{\prime}$
of the same type are all $C_X$ polarized. Now for any $C$ polarized abelian
variety of dimension n, we can define the notion of primitive classes
of the middle degree relative to $C$. So it's natural to consider
the mapping from the moduli of $C$ polarized abelian variety to the classifying
space of primitive classes relative to $C$. 
The image of this mapping, denoted as $M_{PH}$, can be regarded as the classifying
space of Hodge structures that coming from geometry.

The key observation is then
this image of moduli space $M_{PH}$
turns out to be the standard Hilbert moduli varieties of the reflexive $K^*$.
In other words given the type of the $X$ we can construct a moduli
space which is a Hilbert moduli variety of $K^*$. 

To see this we note that all these CM abelian varieties have a large automorphism
group $U_K$, the unit group of $\cO_K$, and the primitive classes relative to this cone
is simply the invariant classes under this action. Moreover for $\alpha \in K=End(X)_{\bQ}$, 
we have $\alpha \in End(H^1(X))_{\bQ}$,
but $H^n(X)=\wedge H^1(X)$, so we have in fact
$g(\alpha)=\prod_{i}\alpha^{\sigma_i} \in K^*$ acting on $PH^n(X)$, i.e. 
the primitive classes admit
multiplication by $K^*$. So we have

\begin{thm}
\begin{enumerate}
\item
We have a natural isomorphism $M_{PH} \simeq \cH^n /SL_2(\cO_{F^*})$;
\item
The natural morphisms $M_n(\fa) \ra M_C \ra M_{PH}$ are all finite.
\end{enumerate}
\end{thm}

By the theory
of Shimura-Taniyama-Weil, the moduli points
of $X$ on this moduli space will generate the class field of $K$.
The natural modular function on the Hilbert modular variety of $K^*$,
when pulled back to $M_n(\fa)$, will become a natural modular function on 
$M_n(\fa)$. In this way, by directly considering the geometry of $X$,
we can have the class fields of $K$ generated. This answers our first
question.

Note that in our approach we do not emphasis on the use of ``field of
moduli'', indeed since our moduli space of primitive Hodge structures
can actually be regarded as the moduli space of abelian varieties
with multiplication by $F^*$, the natural ``field of moduli'' somehow
lose its meaning on this moduli space
(\footnote{If we insist to use ``field of moduli'', it should be ``field of moduli of the
primitive Hodge Structures''. I am not sure how ``motivic'' this notion of  ``field of moduli''
would be.}).
We only need the natural
modular functions on the moduli spaces, which give the coordinates of the
CM points. In any case our approach is
more natural in view of Hodge theory, and serve the purpose of
generating class fields of $K$ well.

\subsection{}

The second question is more difficult. We need to find a natural way
to interpolate the moduli spaces of different types. For this purpose we
will use the idea of Mirror Symmetry. Precisely if $R_1$ and $R_2$ are two
ideals such that $R_1R_2^{-1}$ is a real ideal, then although the abelian
varieties $X_{R_1}$ and $X_{R_2}$ defined by $R_1$ and $R_2$ are living on
the different moduli space, we will show their ``Mirror partners'' $X_{R_1}^{\prime}$
and $X_{R_2}^{\prime}$ are living on a single moduli, and the Mirrors' field of moduli
can be used to generate the class fields. 

To motivate our
idea, let's consider another approach to the Hilbert's 12th problem,
namely the Stark's conjectures(\cite{Stark}).

The point of departure is to consider the Dedekind zeta function and
Hecke L-functions for the number fields. Hilbert's problem asks for 
a natural transcendental function, indeed for the abelian extension
point of view, nothing can be more natural than these L-functions, as they
transformed under the Galois group explicitly.
Usually the Stark's conjectures are formulated and studied for a totally
real field, but as we shall see, it is more natural and simple to study
it for the CM field, because all the finite extension of CM fields
are necessarily CM, i.e., any extension of $K$ has no real infinity.

The remarkable fact about the Stark's conjectures is that it can be
formulated on any class fields of $K$ uniformly, not just those
class fields in the Shimura-Taniyama-Weil theory, so in the explicit
form it can be used  as a guide for us to search for the 
solution of missing class fields of the old theory.

So let $K$ as above, let $R$ be an ideal class of $K$, recall that the Dedekind
zeta function is defined as $\zeta_K(s)=\sum_{\fa} \frac{1}{N(\fa)^s}$, where
the sum is over all the integral ideals of $K$. Dirichlet's formula gives:

$$\zeta_K(s)=\frac{h_K}{s-1} \kappa_K + \rho_K + O(s-1)$$

where $\kappa_K=\frac{(2\pi)^n R_K}{w_K d_K^{1/2}}$, $\rho_K$ is a constant,
when $K$ is $\bQ$, $\rho_K$ is Euler's constant.

We can also define the partial zeta function associated to an ideal class as
$\zeta_K(s,R)=\sum_{\fa \in R} \frac{1}{N(\fa)^s}$.
Natually $\zeta_K(s)=\sum_R \zeta_K(s,R)$, with the limit formula:

$$\zeta_K(s,R)=\frac{\kappa}{s-1} + \rho_K(R) + O(s-1)$$

The constant $\rho_K(R)$ depends on the ideal class of $R$, and holds
vital information about the class fields of $K$.

Now let $K_0$ be the Hilbert class field of $K$, $[K_0 :K]=m$. By
class field theory $Gal(K_0 /K) \simeq H_K$. For any character
$\chi: Gal(K_0 /K) \ra \bC^{\times}$ we can define the Artin
L-function $L_K(s,\chi)=\sum_{R \in H_K} \chi(R) \zeta_K(s,R)$.
It is well-known that $\zeta_{K_0}(s)=\prod_{\chi} L_K(s,\chi)$. By
the above Dirichlet's formula, by studying the behavier of $s \ra 1$,
we can have some relations between the regulators of $R_K$, $R_{K_0}$,
and $\rho_K(R)$.

To make the formula simple it is more convenient to consider $s \ra 0$,
and by functional equation this is equivalent to $s \ra 1$. In fact if
we write
$$\zeta_K(s,R)=-\frac{R_K}{w_K} s^{n-1}(1+\delta_K(R)s) + O(s^{n+1})$$
then
$$\delta_K(R)=n \gamma + nlog{2\pi} -log|D_K| -\frac{w_K |D_K|^{1/2}\rho_K(R^{-1})}
{2^n \pi^n R_K}$$

From $L_K(s,\chi)=\sum_R \chi(R)\zeta_K(s,R)$ we have

\begin{itemize}
\item
if $\chi = \chi_0$ is trivial, then 
$$L_K(s,\chi_0)=\zeta_K(s)=-\frac{h_K R_K}{w_K} s^{n-1}(1+\delta_K s) + O(s^{n+1})$$
with $\delta_K = \sum_R \delta_K(R)$

\item
if $\chi \neq \chi_0$ is not trivial, then
$$L_K(s,\chi)=\zeta_K(s,R)=-\frac{R_K}{w_K} (\sum_R \chi(R) \delta_K(R))s^n + O(s^{n+1})$$

\end{itemize}

Since $L_K(s,\chi)=\sum_{R \in H_K} \chi(R) \zeta_K(s,R)$, comparing the leading
coefficient of the both side we have

$$-\frac{h_{K_0}R_{K_0}}{w_{K_0}} = (-1)^m (\frac{R_K}{w_K})^m h_K 
\prod_{\chi \neq \chi_0}(\sum_R \chi(R) \delta_K(R))$$

The Stark's conjecture predicts that if we write $L_K(s,\chi)=R_K(\chi)s^n + O(s^{n+1})$,
then $R_K(\chi)$ also has the form of regulators, i.e. it is a determinent of
a matrix whose entries are linear combination of logarithm of units of $K_0$. 
Moreover we should expect $R_K(\chi)^{\sigma}=R_K(\chi^{\sigma})$ for any
$\sigma \in Aut(\bC)$.

In our present case $K$ is a primitive CM field, in this case Stark's conjecture
is somehow simple in the following sense: the regulator $R_K$ of $K$ is a determinent
of $(n-1) \times (n-1)$ matrix, and the regulator $R_{K_0}$ of $K_0$ is a determinent of
$(mn-1)\times (mn-1)$ matrix. Stark's conjecture in fact says that the quotient $R_{K_0}/R_K$
is a determinent of $(mn-n)\times (mn-n)$ matrix, and this matrix should be diagonalized
into $m-1$ blocks, with each block a $n \times n$ matrix, further in this case we can use
the $R_K$ matrix to simplify these blocks, so in this extension there are only $m-1$ essential
new units. In other words, $U_{K_0}$ as a module of $U_K$ is free of rank $m-1$.
This is very much similar to the classical case of imaginery quadratic fields.

By some elementry argument we can show that
$R_{K_0}=(R_K)^m \cdot vol(S)$, where $vol(S)$ is a determinent
of $(m-1) \times (m-1)$ matrix whose entries are linear combinations
of logarithm of units in $K_0$. These are the basis of $U_{K_0}$ as the module of $U_K$. 
 On the other hand by the Frobenius
determinent formula,

$$\prod_{\chi \neq \chi_0}(\sum_R \chi(R) \delta_K(R))=
\det_{R_1,R_2 \neq 1}(\delta_K(R_1)-\delta_K(R_2))$$

Hence we should expect :

\begin{itemize}
\item
$\delta_K(R)=log |\eta_K(R)|$;
\item
$\delta_K(R_1) -\delta_K(R_2)=log |\frac{\eta_K(R_1)}{\eta_K(R_2)}|$,
 and $\frac{\eta_K(R_1)}{\eta_K(R_2)}$ is a unit in $K_0$. More
over these units should transform under the Galois group according to the
reciprocity law.

\end{itemize}

\subsection{}

To prove things like this we need to have a good expression for $\delta_K(R)$,
this can be done by using the theory of $GL(2)$ Eisenstein series over
the totally real field $F$. 
The Eisenstein series is defined as

$$E(w,s;\fa)=\sum_{(c,d) \in (\fa \oplus \cO_F )/U_F, (c,d) \neq (0)}\prod_{i=1}^{n}
 y_i^{s}|c^{\sigma_i}z_i + d^{\sigma_i}|^{-2s}
$$

The idea is a classical one, since $\zeta_K(s,R)=\sum_{\fa \in R}\frac{1}{N(\fa)}$,
fix an ideal $\fa_1 \in R^{-1}$. Regarding $\fa_1$ as a module over $\cO_F$,
$\fa_1 \simeq \fa \cdot w_1 + \cO_F \cdot w_2$, where $\fa$ is an fractional ideal
of $\cO_F$, and we choose $w_1, w_2 \in K$ such that $w=\frac{w_1}{w_2}$ satisfying
$\Im w^{\sigma_i} >0$, $\forall i$. Then we have

$$
\begin{array}{lll}
\zeta_K(s,R) & = & N(\fa_1)^s \sum_{\alpha \in \fa_1/U_F, \alpha \neq 0} \frac{1}{N(\alpha)^s}\\
 &=& N(\fa_1)^s \sum_{(c,d) \in (\fa \oplus \cO_F)/U_F,(c,d) \neq (0)} N(cw_1 + d w_2)^{-s} \\
 &=& N(\fa_1)^s \sum_{(c,d) \in (\fa \oplus \cO_F)/U_F,(c,d) \neq (0)}N(w_2)^{-s}\prod_i|c^{\sigma_i}w^{\sigma_i} + d^{\sigma_i}|^{-2s} \\
&=& N(\fa_1)^sN(w_2)^{-s}\prod_{i}\Im(w^{\sigma_i})^{-s}E(w,s;\fa)\\
\end{array}
$$ 

The advantage of using Eisenstain series is that it has an explicit Fourier expansion
at infinity of $\cH^n$, 
the calculation is long and can be found in Yoshita's book(\cite{Yoshita}),
the final formula is :

$$
E(w,s;\fa)=-2^{n-2}h_FR_Fs^{n-1}[1 + (CONST + logN(\fa) + log(\prod_i \Im(w_i)) -h(w;\fa))s]+O(s^{n+1})
$$

where CONST is some universal constant, and
we write $w_i = w^{\sigma_i}$ in the case of no confusion, $w=(w_1,\cdots,w_n)$. And 
$$
h(w;\fa)=\sum_{\chi_F} \chi_F(\fa)h_{\chi_F}(w;\fa)
$$

where 
\begin{itemize}
\item
$\chi_F$ is a character of the ideal class group of $F$;
\item
$\fa$ is the type of $R$;
\item 
$\w$  is the CM point defined by $R$ (with the type $\fa$);
\item
$h_{\chi_F}(\w;\fa)$ is a function on the product of
upper half planes who has the following Fourier expansion:

$$
\begin{array}{lll}
h_{\chi_F}(\w;\fa)&=&\frac{D_F N(\fa)}{2^{n-2}\pi^n h_F R_F}[\chi_F({\fd_F})L_F(2,\chi_F^{-1})\prod_i\Im(\w_i)\\
&+&\pi^nD_F^{-3/2}\sum_{0\neq b \in \fd_F^{-1}\fa} \sigma_{1,\chi}(bda)|N(b)|^{-1}
exp(2\pi i (\sum_{j=1}^n b_j\Re(w_j) + i|b_j\Im(w_j)|))]\\
\end{array}
$$
where
   
\begin{enumerate}
\item
$a,d \in \bA_F^{\times}$ such that $div(a^{-1})=\fa$ and $div(d)=\fd_F$;

\item
$\sigma_{s,\chi}(x)$ is a function defined as
\begin{displaymath}
\sigma_{s,\chi}(s)=\prod_{v\in (finite\ primes)} \left\{ \begin{array}{ll}
1+\chi_v(w_v)q_v^s + \cdots + (\chi_v(w_v)q_v^s)^{ord_v(x_v)} & if \ x_v \in \fd_v,\\
0 & if \ x_v \notin \fd_v.
\end{array} \right.
\end{displaymath}
Here $\fd_v$ denotes the ring of integers of $F_v$, $w_v$ is the prime element of
$F_v$ and $q_v=|\fd_v/w_v \fd_v|$.

\end{enumerate} 
\end{itemize}

So we have
$$
\zeta_K(s,R)=-\frac{R_K}{\w_K}s^{n-1}[1 +(CONST+ logN(\fa) +log\prod_i\Im(w_i)-h(w;\fa))s]+O(s^{n+1})
$$

Since in the end we are dealing with the difference $\delta_K(R_1) - \delta_K(R_2)$, the
universal constant CONST would be cancelled out, so essentially we have

$$\delta_K(R)=h(\w;\fa) -log\prod_{i}
\Im(w_i) -logN(\fa)$$

We reminded that $\fa$ is the type of $R^{-1}$, $w=(w_1,\cdots,w_n)\in \cH^n$ is the CM point
defined by $R$, and $h(w,\fa)$ is the complicated function defined above. So our basic task 
is to understand this function.

\subsection{}

We notice that $h_{\chi_F}$ satisfying the following modular properties(\cite{Yoshita}):
for $\gamma \in \Gamma_{\fa}$ we have
\begin{enumerate}
\item
$h_{\chi_F}(\gamma \w; \fa)=h_{\chi_F}(\w,\fa)$ if $\chi_F$ is not trivial;
\item
$h_1(\gamma \w;\fa) -log(\prod_i \Im(\gamma \w)_i)=h_1(\w,\fa) -log(\prod_i \Im\w_i)$
\end{enumerate}
So in particular we have

$$\sum_{\chi_F}\chi_F(\fa)h_{\chi_F}(\gamma \w;\fa) -log(\prod_i \Im(\gamma \w)_i)
=\sum_{\chi_F}\chi_F(\fa)h_{\chi_F}(\w,\fa) -log(\prod_i \Im\w_i)$$

which suggests that we may actually have a Hilbert modular form 
$\eta_K(w;\fa)$ of parallel weight  
such that $h(w;\fa)=log|\eta_K(w;\fa)|$.
Classically in case $F=\bQ$ this is indeed the case as $\eta_K$ is
the classical Dedekind eta function $\eta$, it is well-known that $\eta$ has
an infinite product expression, which when taking the logarithm translated
into a Fourier expansion, which is exactly the above Fourier series.

In the higher dimensional case 
it is not that easy 
\footnote{Actually for a long time I believe we can have an analogue infinite
product formula for any $\eta_K$, just like the one of $\eta$ has. But after many unsuccessful try, 
I came to the conclusion that such infinite product can not exist.}.
Besides the fact that the Fourier expansion is too complicated and difficult to work with,
we can not expect by directly exponenciate the above expression of $h$
to get a meaningful function, precisely because of the infinite unit
group of $F$, as shown as the regulator term $R_F$ in the leading coefficients
of the Fourier expansion. The regulator $R_F$ is not a rational number, in fact
we expect it to be transcendental, so even if we exponenciate $h$ we
can not get anything useful for the arithmetic purpose. In particular we can
not expect to get an infinite product expression as the classical $\eta$ function.

\subsection{}
What should I do? It turns out although we can not exponenciate the $h$ function, 
can not get the explicit formula for $\eta_K(w;\fa)$, we still can say something 
quantitatively about it. 

The idea is to consider a twisted version of Eisentein series: 

$$E_{u,v}(w,s;\fa)=\sum_{(c,d) \in (\fa \oplus \cO_F )/U_F, (c,d) \neq (0)}\prod_{i=1}^{n}
 e^{2\pi i(c^{\sigma_i}u_i + d^{\sigma_i}v_i)}(\Im w_i)^{s}|c^{\sigma_i}w_i + d^{\sigma_i}|^{-2s}
$$

where $(u,v) \in \bR^n \oplus \bR^n$.

This is entirely adopted from the classical Kronecker's second limit formula(see \cite{Lang}).

To explain why we need to develop the twisted Eisenstein series, let's recall the classical
situation, i.e., when $F=\bQ$.

In this case the Eisenstein series is 
$$
E(w,s)=\sum_{m,n \in \bZ,(m,n) \neq 0}\frac{(\Im(w))^s}{|mw+n|^{2s}}=\frac{1}{2}[1+(CONST+log\Im(w)-4log|\eta(w)|)s] +O(s^2)
$$
 The twisted Eisenstein series is:
$$
E_{u,v}(w,s)=\sum_{m,n \in \bZ,(m,n) \neq 0} e^{2\pi i(mu+nv)}\frac{(\Im(w))^s}{|mw+n|^{2s}}=-2log|g_{-v,u}(w)|s +O(s^2)
$$

where $\eta(w)$ is the Dedekind $\eta$ function:
$$
\eta(w)=q_w^{\frac{1}{24}}\prod_{n=1}^{\infty}(1-q_w^n)
$$

and $g_{u,v}$ is the Siegel's function:
$$
g_{u,v}=-q_w^{\frac{1}{2}B_2(u)}e^{\pi i v(u-1)}(1-q_z)\prod_{n=1}^{\infty}(1-q_w^nq_z)(1-q_w^n /q_z)
$$
with $z=u-vw$, $B_2(u)=u^2-u+1/6$ the Bernoulli polynomial, and for any variable $z$, we write
$q_z=e^{2 \pi i z}$.

Let's note several facts:

\begin{enumerate}
\item
For $(u,v) \notin \bZ \oplus \bZ$, $E_{u,v}(w,s)$ doesn't have pole at $s=1$, hence by functional equation, the
vanishing order at $s=0$ is 1. This is different from the Eisenstein series.

\item
The advantage of using Siegel function $g_{u,v}$ is that we may regard it is a function of a new variable
$z$, an in this respect it is very close to the theta function
$$
\phi(w,z)=(q_z -1)\prod_{n=1}^{\infty}(1-q_w^nq_z)(1-q_w^n /q_z)
$$

Recall that such theta function is characterized up to a constant by

$$
\phi(w,z+1)=\phi(w,z); \phi(w,z+w)=-\frac{1}{q_z}\phi(w,z)
$$

and we have
$$
g_{u,v}=q_w^{\frac{1}{2}B_2(u)}e^{\pi i v(u-1)} \phi(w,z)
$$

\item
$g_{u,v}$ is also closely related to Dedekind $\eta$. In fact as a function of $z$ $g_{u,v}$ has a simple
zero at $z=0$ with $\eta(w)$ as the coefficient, i.e., 
$$
|g_{u,v}|=|\eta(w)|^2 |q_z -1| + O(z^2)
$$

This shows that the absolute value $|\eta(w)|$ is not a theta null, but rather a ``derivative theta null''.
But we note our theta functions are normalized at $z=0$ to be zero, $\phi(w,0)=0$, 
if the theta function normalized in this way, their derivatives can also be used as the
coordinates on the moduli space.
By abuse of notation we still call $|\eta(w)|$ a theta null, hence it gives a modular function
on the moduli space.

\end{enumerate}

From the definition, $E_{u+1,v}=E_{u,v+1}=E_{u,v}$, hence $g_{u+1,v}=g_{u,v+1}=g_{u,v}$, but since
$|g_{-v,u}|=|q_w^{\frac{1}{2}B_2(-v)}|\cdot |\phi(w,z)|$ with $z=u-vw$, we see that this periodic
condition is actually the same as the periodic condition that characterize the theta function $\phi(w,z)$.
This suggests that we may start from this periodic condition to get the theta function directly.

Precisely, assuming that we don't know the explicit infinite product form of $\eta$ and $g_{u,v}$,
only starting from the Eisenstein series:
$$
E_{u,v}(w,s)=-2log|g_{-v,u}|s + O(s^2)
$$
Define
$$
\phi(w,s)=q_w^{-\frac{1}{2}B_2(-v)}g_{-v,u}(w)
$$

Then from the periodic condition of $g_{u,v}$ we immediately see that
$$
|\phi(w,z+1)=\phi(w,z); |\phi(w,z+w)|=|q_z^{-1}|\cdot |\phi(w,z)|
$$

That is, $|\phi(w,z)|$ satisfying the characterization of a theta function,moreover since
we know it is an analytic function, it then has to be a theta function itself.

This is the idea we would follow in the higher dimensional case, as we observed before,
since in the higher dimension we can not expect any explicit infinite product formula
for the function $\eta_K(w;\fa)$, but we still have all the periodic properties as the
1-dimensional case.

Now we go back to the higher dimensional case, the Eisenstein series is:

$$E(w,s;\fa)=\sum_{(c,d) \in (\fa \oplus \cO_F ))/U_F, (c,d) \neq (0)}\prod_{i=1}^{n}
 (\Im(w_i))^{s}|c^{\sigma_i}w_i + d^{\sigma_i}|^{-2s}
$$

We have the limit formula:
$$
E(w,s;\fa)=-2^{n-2}h_FR_Fs^{n-1}[1 + (CONST + logN(\fa) + log(\prod_i \Im(w_i)) -h(w;\fa))s]+O(s^{n+1})
$$

The twisted Eisenstein series is:

$$E_{u,v}(w,s;\fa)=\sum_{(c,d) \in (\fa \oplus \cO_F )/U_F,(c,d) \neq (0)}\prod_{i=1}^{n}
 e^{2\pi i(c^{\sigma_i}u_i + d^{\sigma_i}v_i)}(\Im w_i)^{s}|c^{\sigma_i}w_i + d^{\sigma_i}|^{-2s}
$$

and we have the limit formula:

$$
E_{u,v}(w,s;\fa)=-2^{n-2}h_FR_F log|g_{-v,u}(w;\fa)|s^n + O(s^{n+1})
$$
where $log|g_{-v,u}(w;\fa)|$ has an explicit Fourier expansion, just like $h(w;\fa)$.

Then we argue as the following:
\begin{enumerate}
\item
Recall that $u=(u_1,\cdots,u_n) \in \bR^n$, $v=(v_1,\cdots,v_n) \in \bR^n$, and $\cO_F \subset \bR^n$
as a lattice. From the definition, for any $\alpha \in \cO_F$, we have
$E_{u+\alpha, v}=E_{u, v+\alpha}=E_{u,v}$, i.e., translate invariant under $\cO_F$,
hence $|g_{-v + \alpha,u}|=|g_{-v,u+\alpha}|=|g_{-v,u}|$.

\item
Now write $z=u-vw$, i.e., $z=(z_1,\cdots,z_n)$, $z_i=u_i-v_i w_i$,$\forall i$. We try to write
$g_{-v,u}$ as a function of $(w,z)$, so we define $\phi(w,z)=q_w^{-\frac{1}{2}B_2(-v)}g_{-v,u}(w)$,
where $q_w^{-\frac{1}{2}B_2(-v)}=\prod_i q_{w_i}^{-\frac{1}{2}B_2(-v_i)}$.

Then from the periodic property of $g_{-v,u}$ we immediately have
$$
|\phi(w,z+\alpha)|=|\phi(w,z)|;\ |\phi(w,z+\alpha w)|=|q_z^{-\alpha}||\phi(w,z)|
$$

Recall that on our abelian variety $X=\bC^n /(w \fa \oplus  \cO_F)$ the theta function is
characterized by 
$$
\theta(z+\alpha)=\theta(z); \ \theta(z+\alpha w)=q_z^{-\alpha}\theta(z)
$$

Since $g_{-v,u}$ is an analytic function, we conclude that $|g_{-v,u}(w,z)|$ as a function
of $z$ is the absolute value of a theta function.

\item
From the explicite Fourier expansion of $log|g_{-v,u}|$ and $h(w;\fa)$ we see they
are closely related, in fact if we write $h(w;\fa)=2log|\eta_K(w;\fa)|$,we can verify directly
$$
\lim_{z\ra 0}\{log|q_w^{\frac{1}{12}}\phi(w,z)| -log(|\eta_K(w;\fa)|^2\prod_i |z_i|) \}=0
$$
that is 
$$
\lim_{z \ra 0}\frac{|q_w^{\frac{1}{12}} \phi(w,z)|}{|\eta_K(w;\fa)|^2\prod_i|z_i|}=1
$$

Also we verify that our theta function is normalized at $z=0$ to be 0.

This implies that $\eta_K(w;\fa)$ is a theta null.

\end{enumerate}
 
By the classical theory of theta function, theta null naturally gives rise to the 
modular forms on the moduli space. In fact this should be  more or less expected.
By Mumford's theory of algebraic theta function, we may further conclude that these
theta nulls in fact defines the moduli space as an integral scheme over $\bZ$, hence
have all the expected integral properties. 

\subsection{}

In summary we have found the explicit form of the function $\delta_K(R)$.
$$
\begin{array}{lll}
\delta_K(R)&=& logN(\fa) + log\prod_i \Im(w_i) -h(w;\fa)\\
&=& logN(\fa) + log \prod_i \Im(w_i) -log|\eta_K(w;\fa)|^2\\
&=& log[N(\fa)\prod_i \Im(w_i) |\eta_K(w;\fa)^{-2}|]\\
\end{array}
$$

where $R$ is an ideal class of $\cO_K$, $\fa$ its type, $w$ its CM point on $M_n(\fa)$,
and $\eta_k(w;\fa)$ is the theta null:
$$
\frac{\partial}{\partial z_1}\cdots \frac{\partial}{\partial z_n}\phi(w,z)|_{z=0}=\eta_K(w;\fa)^2
$$

By the Stark's conjecture, we need to understand $\delta_K(R_1) - \delta_K(R_2)$, i.e., we need to
understand the quotient $\frac{\eta_K(w_1;\fa_1)}{\eta_K(w_2;\fa_2)}$. When $R_1$ and $R_2$
are of the same type, $\eta_K(w;\fa)$ are the modular forms on the same moduli space $M_n(\fa)$,
hence $\frac{\eta_K(w_1;\fa)}{\eta_K(w_2;\fa)}$ is meaningful, as the modular function evaluating
at the CM points, and the natural action of Galois group on them is prescribed by the reciprocity law.
But when $R_1$ and $R_2$ are of the different type, for example, when $R_1R_2^{-1}$ is a real ideal,
then $\eta_K(w;\fa_1)$ and $\eta_K(w;\fa_2)$ are on the different moduli spaces, thus their
quotient becomes meaningless. This is precisely the limit of Shimura-Taniyama-Weil's theory.

So what can we do? To go further we need to find a natural way to interpolate the different
moduli spaces, and it is here the idea of Mirror symmetry comes. As we shall see, in this case this 
quotient will have a 
meaning similar to the classical one if we consider the complexified Kahler
moduli of the abelian varieties.

Why should we interpolate the different moduli spaces?
The function $\delta_K(R)$  is defined not on a single moduli space of the
fixed type, but rather automatically been defined on all the moduli spaces,
as the type $\fa$ can vary accordingly. Likewise the modular form
$\eta_K$ is a Hilbert modular form on all the type-fixed
Hilbert modular varieties, and when the type vary, can be regarded as a modular form
on all the moduli spaces. This strongly suggests that we should have
a natural way to interpolate all these moduli spaces of different types,
such that these functions can naturally defined. In other words,
when we fix the type $\fa$, we get the Hilbert modular forms, what then
happens if we fix the CM points $\w$ and let the type vary?

When we look the explicit form of $\delta_K(R)$ and $h$, we note the
apparent symmetric roles played by the quantities $\prod_{\sigma}\Im(\w^{\sigma_i})$
and $N(\fa)$. Indeed if we fix the CM points $\w$ and let $\fa$ vary, we should
have a meaning for the quantity $N(\fa)$. This can be achieved by considering
the Kahler moduli of our abelian varieties.

\subsection{}
Mirror Symmetry is usually formulated for the Calabi-Yau varieties, roughly it asserts
that Calabi-Yau always come in pairs, $X$ and $X^{\prime}$, with the ``complex moduli''
and ``Kahler moduli'' exchanged. In terms of the Hodge number, it means the Hodge diamond
of $X^{\prime}$ is a rotation of $X$.

For abelian varieties, Mirror Symmetry is generally regarded as ``trivial'', as the underlying
topological type would not change. Nevertheless we can still talk about it. There are
several constructions of Mirror manifold for the abelian varieties, the simplest one I
believe, is given by Manin(\cite{M}).
It goes as the following: let $k$ be any complete field,
$X$ an abelian variety over $k$, $T$ be the algebraic torus of dimension
$n$ over $k$, $T\simeq (k^{\times})^n$. Then the multiplicative uniformization
is $0 \rightarrow P_X \rightarrow T \rightarrow X \rightarrow 0$, where $P_X$ is a free
abelian group of rank $n$, $P_X$ is called the period of $X$. Under this
uniformization the Mirror partner $X^{\prime}$  is then
$0 \rightarrow P_X \rightarrow T^{\vee} \rightarrow X^{\prime} \rightarrow 0$,
i.e., we explicitly indentify the periods in $T$ and $T^{\vee}$. When $k \simeq \bC$
one verify that the complex moduli and Kahler moduli of the two are exchanged.

In our situation, given the Mirror pair $X$ and $X^{\prime}$ we will be mainly concern
about the relations of theta functions on them. Since the underlying topological
type would not change, we may regard the Mirror transform as a ``rotation'' of complex
structure of $X$. So to compare the theta functions on $X$ and $X^{\prime}$ we have
to fix the underlying real structures.

We begin with $X=\bC^n /(w \cdot \fa \oplus \cO_F) \simeq \bR^{2n}/\bZ^{2n}$, any polarization of $X$
is given by an integral skew-symmetric bilinear form on $\bR^{2n}$. Given such a form $\w$,
we can find an integral basis $\{\lambda_1,\cdots,\lambda_{2n} \}$ of the integral lattice such
that if $\{ x_1, \cdots,x_{2n} \}$ is the dual basis, then $\w=\sum_{i=1}^n \delta_i dx_i \wedge dx_{n+i}$,
with $\delta_1|\delta_2|\cdots $ the elementary divisors.

Note in our case the biliear form is given by the trace $Tr_{K/\bQ}(\zeta x y ^{\rho})$ with the
admissible $\zeta \in K$ such that $\zeta^{\rho}=-\zeta$, $\Im(\zeta^{\sigma_i}) >0$. Thus we
may regard $(x_1,\cdots,x_n)$ as an integral basis of $\cO_F$, and $(x_{n+1},\cdots,x_{2n})$ as an
integral basis of $\fa$. In particular $(x_{n+1},\cdots,x_{2n})$ depends on $\fa$. In the following we will
denote it as $x_{n+i}(\fa)$ if we need to use this dependence.

Next we introduce the complex structure, so $X$ becomes a complex tori, and we can introduce the complex
coordinates. To do this let $e_i=\lambda_i /\delta_i$, $i=1,2,\cdots,n$, and let $\{z_i\}$ be the complex
dual of $\{e_i \}$. Consider the change of coordinates transform:
$$
\Omega \cdot (x_1,\cdots,x_{2n})^T=(z_1,\cdots,z_n)^T
$$

Then $\Omega=(\Delta_{\delta}, Z)$ with $\Delta_{\delta}=diag(\delta_1,\cdots,\delta_n)$ the diagonal
matrix, and $Z$ symmetrical, $\Im(Z)>0$. We recogonize that $\Delta_{\delta}^{-1}Z$ is the period matrix of $X$. 
Note
in our case for abelian varieties with CM by $K$, the peiod martix $Z$ is necessarily diagonal
$\Delta_{\delta}^{-1}Z=diag(w_1,\cdots,w_n)$, with $w=(w_1,\cdots,w_n) \in \cH^n$.

In particular we have $z_i=\delta_i \cdot x_i + \delta_i \cdot w_i \cdot x_{n+i}(\fa)$, we may regard it
as the transform from the underlying real coordinates to the complex coordinates.

Now recall the theta function on $X$ is characterized by the periodic condition:
$$
\theta(z+\lambda_i)=\theta(z); \ \theta(z+\lambda_{n+i})=e^{-2\pi i z_i} \theta(z)
$$

Taking absolute values we have:
$$
|\theta(z+\lambda_i)|=|\theta(z)|; \ |\theta(z+\lambda_{n+i})|=e^{2\pi \Im(z_i)} |\theta(z)|
$$

However from the above coordinates transformation, 
$$\Im(z_i)=\delta_i \cdot \Im(w_i) \cdot x_{n+i}(\fa)$$

The Mirror symmetry transform
says that we can exchange the complex moduli with the Kahler moduli,
while
the coordinates $w=(w_1,\cdots,w_n) \in \cH^n$ can be regarded as the complex moduli of $X$, 
where is the Kahler moduli? Our Kahler moduli coordinates are actually in the variables
$(x_{n+1}(\fa), \cdots,x_{2n}(\fa))$.
Since they are depend on the type $\fa$, we want to write the dependency explicitly, in 
order to understand the transformation of types.

For this end let's write $x_{n+i}=x_{n+i}(\cO_F)$. The type ideal $\fa$ as a $\bZ$ module, is
a submodule of $\cO_F$ of full rank, i.e., if we fix an integral basis of $\cO_F$, then
$\cO_F/{\fa} \simeq \oplus_{i=1}^n \bZ /{t_i \bZ}$, with $x_{n+i}(\fa)=t_i x_{n+i}$ and $\prod_i t_i= N(\fa)/D_F$. Thus
the positive rational numbers $(t_1,\cdots,t_n)$ can be conviniently regarded as the coordinates of the
ideal $\fa$, and under appropriate identification, can be regarded as the coordinates of 
of the Kahler class in the the Kahler moduli. So in particular we have
$$
\Im(z_i)=\delta_i \cdot \Im(w_i) \cdot t_i \cdot x_{n+i}
$$

But from the above formula, when we exchange the complex moduli $(w_1,\cdots,w_n)$ and Kahler moduli
$(t_1,\cdots,t_n)$,
it's not going to change the multiplier 
$e^{2 \pi \Im(z_i)}=e^{2 \pi\delta_i \cdot t_i \cdot \Im(w_i) \cdot x_{n+i}}$.
Since the absolute value of theta functions can be regarded as  a real analytic function on
$\bR^{2n}$, thus we conclude that for the given Mirror pair $X$ and $X^{\prime}$, their theta functions'
absolute values satisfying the same periodic condition, hence must be only differed by a constant!

Recall our previous puzzle, when two ideal classes $R_1$ and $R_2$ of $\cO_K$ satisfying $R_1R_2^{-1}$ is a real
ideal, then $R_1$ and $R_2$ are of the different type, so the corresponding abelian varieties $X_1$ and $X_2$
living on the different moduli spaces, so $\eta_K(w,\fa_1)/\eta_K(w,\fa_2)$ has no meaning. But we
know $\eta_K(w,\fa)$ is the  theta null of $X$, thus by the above formula, $\eta_K(w,\fa)$ is
also the  theta null of the Mirror $X^{\prime}$. So although $X_1$ and $X_2$ living on the
different moduli space, if their Mirror $X_1^{\prime}$ and $X_2^{\prime}$ are in the same moduli space,
then $\eta_K(w,\fa_1)/\eta_K(w,\fa_2)$ would be meaningful! This is indeed the case, as 
$X_1^{\prime}$ and $X_2^{\prime}$ are on the single moduli space, the complexified Kahler moduli space of
$X$. This is the underlying rationale
for us to use the Mirror symmetry.

Note from this relation of theta functions we also see that if $X$ is defines over a number field, then
the Mirror $X^{\prime}$ also is defined over that number field.

\subsection{}
How to construct the Kahler moduli?
Recall that when we consider the Kahler cone of $X$ and the notion
of $C_X$ cone polarized abelian varieties, we know that for a fixed
type $\fa$, all these $X^{\prime}$ are $C_X$ polarized. Hence for
the fix-typed abelian varieties, their Kahler cones are more or less the
same. Changing the type really means changing the Kahler cone. 

To do this first we observe
that all our $X$ are isogenies to each other, so their Kahler cones
are comparable. For a fixed $\zeta \in K$ satisfying 
$\zeta^{\rho}=-\rho$, $\Im(\zeta^{\sigma_i}) >0$, any $i$. Such $\zeta$
defines a polarization for any $X$, let's denoted as $p_X$.
For any two $X_1$ and $X_2$, we have an isogeny
$f:X_1 \ra X_2$, the pull back of the polarization $f^*(p_{X_2})$ is
a polarization on $X$, hence defines a point in the Kahler cone of $X$.
So it makes sense to consider the Kahler moduli space, if $C_X$ is the
Kahler cone of abelian variety $X$, and if $U_X$ is the units group
acts as the automorphisms of $X$, the Kahler moduli is then the
quotient $C_X/U_X$. As indicated above, once we fix a $\zeta \in K$,
then different $X$ can have different moduli points in $C_X/U_X$.
Recall that for the abelian variety $X$ defined by ideal $\fa$, we have
introduced a coordinate $(t_1,\cdots,t_n)$ for $\fa$ such that $\prod_i t_i = N(\fa)/D_F$, 
under this  identification we see that $(t_1,\cdots,t_n)$ can be regarded as the
Kahler coordinate for the variety $X$,   hence makes
the symmetric role of $\prod_{\sigma}\Im(\w^{\sigma_i})$ and $N(\fa)$
transparent.

So far we have only considered the real part of the Kahler moduli, 
to get a reasonable moduli space in the category of algebraic varieties,
we may consider the complexified Kahler moduli space, which means
we consider the tube domain $C_X(\bC)=C_X \otimes \bR \oplus i C_X$,
and we introduce the natural discrete group action $\Gamma$ on
$C_X(\bC)$, which is an extension of the natural unit group action
on $C_X$, and we take the quotient to form $N_X=C_X(\bC) /\Gamma$,
which is the complexified Kahler moduli. Note that this is also consistent
with the view that $\prod_{\sigma}\Im(\w^{\sigma_i})$ and $N(\fa)$ should
be made symmetrical.

The complexified Kahler moduli is again a finite covering 
space of the standard Hilbert modular variety
of $F$. In fact since the unit group $U_F$ acts on the Kahler cone $C_X$
via the square map $U_F \ra U_F^2$, let $\Gamma_2$ be the smallest
congruence subgroup of $SL_2(\cO_F)$ that contains $U_F^2$, $\Gamma_2$ is
of finite index in $SL_2(\cO_F)$, then we have:

$$N_X \simeq \cH^n /\Gamma_2$$

and the Kahler moduli points of $X$ on it are all some algebraic points.

In the spirit of Mirror Symmetry
$N_X$ can be interpretated as the complex moduli space of the Mirror
$X^{\prime}$.
In this way the moduli space of the types are interpolated,
and the moduli points of X on the Hodge moduli
space and Kahler moduli space will generate all the class fields of $K$.
In this way we complete the work of Shimura-Taniyama-Weil. Note that
all our invariants can be constructed from the cohomology of $X$, so
our approach can be simply stated as a cohomological approach.

Some further investigation along this line, similar
to the classical analysis of Shimura-Taniyama(\cite{S1}),
would show that we can have a similar Kronecker congruence properties
for the isogenies of mirror abelian varieties of the different types.

In summary we can state our main theorem as

\begin{thm}
Let $K$ be a primitive CM field, and let $X$ be an abelian variety
of CM type $(K,\{\sigma_i\})$, then there are well-defined complex
moduli point $p_X$ in the geometric moduli space of primitive classes
of middle degree, and Kahler moduli point $q_X$ in the complexified
Kahler moduli, such that together they generate the Hilbert class fields
of $K$.
\end{thm}

Similarly we have the corresponding result for the ray class fields:

\begin{thm}
Let $\fa$ be an integral ideal of $K$. We can define the set of CM
abelian varieties with the $\fa$-level structure as $\{X,\fa \}$,
and from this there are well defined complex moduli points via
primitive Hodge structures $p_{X,\fa}$ and Kahler moduli points
$q_{X,\fa}$, together they generate the ray class fields $K_{\fa}$
of conductor $\fa$.
\end{thm}

We can define a extended class of modular functions similar to
our $\delta_K(R)$. Let $\cM_X$ be the set of function that

\begin{itemize}
\item
It is a finite function defined on the set of all abelian varieties
$X$ of CM type $\{K,\sigma_i \}$;
\item
When the type if fixed, it is a restriction of Hilbert modular function;
\item
When the CM point $z$ is fixed, it is a restriction of modular function
on $N_X$ to the imaginery axis.
\end{itemize}

The quotient $\eta_K(R_1)/\eta_K(R_2)$ is certainly belong to $\cM_X$.
For the function in $\cM_X$ 
we can formulate a Shimura reciprocity law, similar to
the classical approach.

\begin{thm}
Let $s \in \bA_K$ be an adele, $(s^{-1},K)$ be the Artin symbol, and
we define the action of $s$ on the set of $X$ as $s(X)$ by isogenies,
then
\begin{enumerate}
\item
For any $f  \in \cM_X$, $(s^{-1},K)$ acts on the value of $f(X)$;
\item
Precisely we have $f(X)^{(s^{-1},K)}=f(s(X))$.
\end{enumerate}
\end{thm}

\subsection{}

Now we can return to the Stark's conjecture,
if we write $\eta_K(R)=\frac{\eta_K(w;\fa)}
{\prod_i \Im(\w^{\sigma_i}) N(\fa)}$, with $\fa$ the type of $R$
and $w$ the CM point defined by $R$, then we have
$$\delta_K(R)=log|\eta_K(R)|$$

Since $\eta_K(w;\fa)$ is the theta null, by Mumford's algebraic theory of theta functions, $\eta_K(w;\fa)$
can also be interpretated as the algebraic coordinates on the moduli space, hence $\frac{\eta_K(R_1)}{\eta_K(R_2)}$
is necessarily an algebraic number. By the rationality of the abelian variety $X$ and its Mirror $X^{\prime}$,
and by the Shimura reciprocity law, we have

\begin{thm}
For a primitive CM field $K$, with $\eta_K$ defined as above, we have
\begin{enumerate}
\item
$\epsilon_{R_1R^{-1}_2}=\eta_K(R_2) /\eta_K(R_1)$ is an algebraic number in $K_0$;
\item
If $S$ is an ideal class and $(S^{-1},K)$ is the Artin symbol, then
$\epsilon_R^{(S^{-1},K)}=\epsilon_{RS} \epsilon_S^{-1}$
\end{enumerate}
\end{thm}

Stark's conjecture actually predicts $\frac{\eta_K(R_1)}{\eta_K(R_2)}$ is a unit in the Hilbert class field
$K_0$. It is this statement I find quite elusive. For although it's quite clear the quotients
$\frac{\eta_K(R_1)}{\eta_K(R_2)}$ are  algebraic numbers, and 
they are trnasformed under Galois group according to the reciprocity law, it's hard to assert 
this quantity also has to be an algebraic integer. In the classical situation
we have the explicit Fourier expansion of Dedekind's $\eta$ from the infinite product formula,
with all the coefficients integers. So some general theorem somehow garantteed that their quotient
has to be an algebraic integer. But in our case we simply don't know the explicit form of $\eta_K(w,\fa)$,
let alone their Fourier expansions. The formulism so far developed seems not adquate for this problem, we will
comment more on this at the end of the paper.

In this way the Stark's conjectures are partially proved for $K$.
 
\subsection{}

Since the Stark's conjecture can be formulated for a general number field,
we may ask the question that if it is true for the primitive CM field, then
how about a general number field? For example a totally real field?

In a certain sense our results can be used to answer this question, for
example if $F$ as above is a totally real field, and if $F_1$ is a class
field of $F$. Then we can find a imaginery quadratic extension of $F$, $K$, such that
the composition field $F_1 K$ is a class field of $K$, and $F_1$ is realized
as a subfield of $F_1 K$. Since we know how to explicitly generate $F_1K$ over
$K$, we know at least in princeple how to generate $F_1$. The Stark's conjecture
for $F$ should be proved along this line.

For a general number field $F$ we can follow an analogue path, namely for
any class field $F_1$ of $F$ we may find a CM field $K$ which is an extension
of $F$ such that $F_1K$ is a class field of $K$, and $F_1$ can be embedded in 
$F_1K$ as a subfield. If we can find such a $K$ then we may follow the same step 
as above.

This approach however has something unnatural in it, we need to find an
extra CM field to generate the class field of a given field. It's like using
the elliptic functions to write down the roots of unity. To get the truely natural
transcendental function as required by Hilbert, some further works are needed.

\subsection{}

The above is the brief summary of the main results, the plan for the rest of
the notes is the following. From section 1 to section 5 we will study the class field
generation problem, we will more or less follow Shimura's expository(\cite{S1}). Then from
section 6 to section 8 we will study the Stark's conjecture, by applying the Shimura
reciprocity law. 

Section 9 is rather independent, it consists an integral representation formula
for the function $h$ as a singular theta
lifting of Borcherds' type. This is one of my earlier attempts to understand
the forms $\eta_K(w;\fa)$, by looking at its singularities at the boundary divisors of moduli space.
Although it didn't give the infinite product as I wanted, I feel
it may be of some independent interest, so I include it here.

Finally in section 10 we discuss some further possible works.

\section{\bf Cone Polarized Hodge Structures and Moduli Spaces}
\subsection{}

First let's fix the notations.\\
$F$ --- a totally real field of degree $n$ over $\bQ$; \\
$K$ ---  a primitive imaginery quadratic (CM) extension of $F$; \\
$K^*$ --- the reflexive field of $K$, also a CM field; \\
$F^*$ --- the maximal real subfield of $K^*$, $[K^*:F^*]=2$,$[F^*:\bQ]=n$.\\
$\{\sigma_1,\cdots,\sigma_n\}$ --- the set of embeddings of $F$ into $\bC$, lifted
to the embeddings of $K$ as the CM-type;\\
$\rho$ --- the complex conjugate.\\

In general $F$ and $F^*$ are different, $K^*$ is always primitive, and since $K$ is primitive,
$K$ is the reflexive field of $K^*$.
For any ideals $\fa$ of $K$, define 
$f(\fa)=\fa\cdot \fa^{\rho}$, $g(\fa)=\prod_{i=1}^n \fa^{\sigma_i}$.
Then $f(\fa)$ is an ideal in $K$, while $g(\fa)$ is an ideal in $K^*$. We call an ideal
$\fa$ in $K$ the real ideal if $\fa=\fa^{\rho}$,i.e., if
$\fa$ is coming from an ideal in $F$, and the imaginery ideal if 
$\fa=g(\fc)$, for some ideal
$\fc$ in $K^*$. If $\fa$ is imaginery, then $\fa\fa^{\rho}=(\mu)$,
for some $\mu \in F$.

Let $X$ be an abelian variety of dimension $n$ with CM type $(K,\{\sigma_i\})$, 
we have ${\rm End}(X)_{\bQ}\simeq K$. Such $X$ can be constructed in the following way:
let $R$ be an order in $K$, then as a $\bZ$-module $R$ is of rank $2n$. Using
$\{\sigma_i\}$ we define an embedding:
$$\iota: R \longrightarrow \bC^n, \iota(x)=(\sigma_1(x),\cdots,\sigma_n(x));$$
then all such $X$ can be constructed as $X \simeq \bC^n / \iota(\fa)$, with
$\fa$ an ideal in $R$ and
$R \simeq {\rm End}(X)$.  In particular if 
$U_R$ is the units group of $R$, then we have ${\rm Aut}(X)=U_R$. 
In the simplest case when $R=\cO_K$,
we write the units group as $U_K$. By Dirichlet's theorem, we have
$U_R \simeq \bZ^{n-1} \oplus Torsion$. This means our $X$ has a large automorphism
group.

We want to construct natural moduli spaces for the above CM abelian varieties
$X$ to live on. Since $X \simeq \bC^n / \iota(\fa)$, to classify such $X$ we
only need to consider all the ideals in $\cO_K$. Let $\fa$ be such an ideal,
considered as a module of $\cO_F$, $\fa$ is of rank 2, in fact we have
the following:

\begin{lem}
For any ideal $\fa$ of $\cO_K$ we have $\fa \simeq \fb \cdot w_1 \oplus \cO_F \cdot w_2$,
where $w_1$, $w_2 \in K$, $\fb$ a fractional ideal of $F$. Moreover
$\fb$'s ideal class is equal to 
$\fc\cdot N_{K/F}(\fa)=\fc\cdot \fa\fa^{\rho}$, where $\fc$ is a
fractional ideal of $F$, independent of $\fa$, and $\fc^2=D_{K/F}$.
\end{lem}

In particular from this lemma we have $\cO_K \simeq \fc w_1 \oplus \cO_F \w_2$,
and if $\fa$ is an imaginery ideal then $N_{K/F}(\fa)= (\mu)$, then
$\fa \simeq \fb w_1(\fa) \oplus \cO_F w_2(\fa)$.
If $\fa$ is a real ideal then $N_{K/F}(\fa)=\fa_0^2$, $\fb=\fc\fa_0^2$,
so $\fa \simeq \fc\fa_0^2 w_1 \oplus \cO_F w_2$.

If $X \simeq \bC^n / \iota(\fa)$, then by Shimura, the polarization of 
$X$ is given by the Riemann form
$E(x,y)$ which is   
$$E(x,y)={\rm Tr}_{K / \bQ}(\zeta x y^{\rho})$$
with $\zeta \in K$ satisfying
$\zeta^{\rho}=-\zeta$, $\Im(\zeta^{\sigma_i}) > 0$, $\forall i$

The ``type'' of the polarization, by Shimura(\cite{S1}), is the ideal class of
 $\zeta\fd_{K/F}\fa\fa^{\rho}$ where
$\fd_{K/F}$ is the different of $K$ over $F$. This ideal  is actually a real
ideal by the definition of the different, in fact we have

\begin{lem}
$\fb$ is the type of $X \simeq \bC^n / \iota(\fa)$.
\end{lem}
 
Given the fractional ideal $\fb$ we can construct a moduli space 
$M_n(\fb)$ which parametrizing families of abelian varieties 
$\bC^n /(v_n \cdot \iota(\fb) \oplus \iota(\cO_F))$, where the dot product
$v_n \cdot \iota(\fb)$ is the component-wise product. We have
$M_n(\fb) \simeq \cH^n / \Gamma(\fb)$ with 
$\Gamma(\fb)=\{ \alpha \in Sl_2(\cO_F) | \alpha \equiv 1 \mod (\fb) \}$.
If $X$ is of type $\fb$,$X \simeq \bC^n / \iota(\fa)$ with
$\fa \simeq \fb w_1 \oplus \cO_F w_2$, then up to isomorphism
we can choose $w_1$ and $w_2$ such that $w_X=w_1/w_2$ satisfying
$\Im(w_X^{\sigma_i}) >0$ for any $i$. Hence $w_X$ defines a point
on $\cH^n$, which then defines the moduli point of $X$ on $M_n(\fb)$.

These moduli spaces parametrized abelian varieties with real multiplication
by $F$, and these abelian varieties all have a common feature, that
they have a large automorphism group, and have a large Kahler cones.
These properties suggested the following geometric formulation by
Hodge Structures. 

Recall that the ring $R={\rm End}(X)$ naturally acts on
the cohomology groups $H^1(X,\bZ) \simeq \bZ^{2n}$, also naturally acts on
the space of holomorphic 1-forms $H^0(X, \Omega_X^1)$, the action is
compatible with the Hodge decomposition
$H^1(X, \bZ)  \otimes \bC \simeq H^1(X,\bC) \simeq H^{0,1}(X) \oplus H^{1,0}(X)$

Let $H^{1,1}(X,\bR)=H^2(X,\bR) \cap H^{1,1}(X)$, and let 
${\rm NS}(X)=H^2(X,\bZ) \cap H^{1,1}(X,\bR)$ be the Neron-Severi group of X,
i.e.,the group generated by the first Chern classes of line bundles on $X$.
Since $X$ is an abelian variety, $H^*(X,\bC)$ is generated by $H^1(X,\bC)$.
Let $\{\omega_1,\cdots,\omega_n \}$ be an basis of $H^{0,1}(X)$, then
$\{\omega_i \wedge \overline{\omega_j} + \omega_j \wedge \overline{\omega_i} \}_{i,j}$
is a basis for $H^{1,1}(X,\bR)$. Inside ${\rm NS}(X)\otimes \bR$ let $C_X$ be the cone
generated by the Chern classes of ample line bundle over $X$, $C_X$ is called
the Kahler cone of $X$. Let $C_X(\bZ) \subset {\rm NS}(X)$ be the integral points
in $C_X$ that generate $C_X$.

\begin{lem}
We have ${\rm End}(C_X(\bZ))\otimes \bQ \simeq {\rm End}_0(X)_{\bQ}={\rm End}_0(X) \otimes \bQ$,
where ${\rm End}_0(X)$ is the complex conjugate invariant of ${\rm End}(X)$.
\end{lem}
\begin{proof}
For any element $\alpha \in C_X(\bZ)$, $\alpha$ can be regarded as an integral
bilinear Riemann form on $H^1(X,\bZ)$. Now for any $f \in {\rm End}(C_X(\bZ))$, $f$ is
induced by a $f^{\prime} \in {\rm End}(H^1(X,\bZ)$ such that $f^{\prime}$ preserves
the Hodge decomposition and maps a Riemann form to another Riemann form. 
Such $f^{\prime}$ is always real. 
If such $f^{\prime}$ preserves a Riemann form, then $f^{\prime}$ is a torsion,
as the projective automorphism group is always finite. Hence we must have
${\rm End}(C_X(\bZ))\otimes \bQ \simeq {\rm End}_0(X) \otimes \bQ$.
\end{proof}

\begin{prop}
we have $\dim_{\bR}{\rm NS}(X) \otimes \bR = \dim_{\bR}C_X =n$. i.e., the Picard number
of $X$ is $n$.
\end{prop}
\begin{proof}
Let $r$ be the Picard number of $X$, by Hodge index theorem, when we fix
a polarization on $X$, the Signature of ${\rm NS}(X)\otimes \bR$ is $(1,r-1)$,
so since ${\rm Aut}(C_X(\bZ))$ is abelian, then ${\rm Aut}(C_X(\bZ))\otimes \bQ$ is of
rank $r-1$. On the other hand, by the above lemma, since
${\rm Aut}(C_X(\bZ))\otimes \bQ \simeq {\rm Aut}(X) \otimes \bQ \simeq U_R \otimes \bQ$,
by the Dirichlet unit theorem, $U_R \otimes \bQ$ is of rank $n-1$, hence  $r=n$.
\end{proof}

\begin{rmk}
\begin{enumerate}
\item
This proposition is inspired from a similar result on $K3$ surfaces (\cite{PS}).
\item
It's actually easy to find a basis for $C_X$, in fact let 
$\{\epsilon_1,\cdots,\epsilon_n \}$ be a set of fundamental units of the
order $R$ with $\epsilon_1=1$, then for any chosen basic polarization $c \in C_X(\bZ)$, the
elements $\{ \epsilon^*_1(c),\cdots,\epsilon^*_n(c) \}$ form a n-dimensional
subspace of ${\rm NS}(X)\otimes \bR$, all in the
cone $C_X(\bZ)$, it can be regarded as the basis of $C_X$.
\end{enumerate}
\end{rmk}

\subsection{}

\begin{Def}
Given a cone $C=C(\bZ)\otimes \bR$, a cone polarized abelian variety $X$
is such $X$ that we have a inclusion $C(\bZ) \subset C_X(\bZ)$.
\end{Def}
Intuitively such $X$ has a lot of polarizations.

Next we consider the moduli space of cone polarized abelian varieties.
Let $M_X$ denotes the coarse moduli space of the corresponding functor,
for any $ t\in M_X$ let $X_t$ denote the corresponding abelian variety,
since $C_{X_t} \supset C_X$, by the lemma before we have
${\rm End}_0(X_t)\otimes \bQ \simeq {\rm End}(C_{X_t}(\bZ))\otimes \bQ \supset 
{\rm End}_0(C_X(\bZ)) \otimes \bQ = F$, i.e., $X_t$ always has a multiplication
by $F$.

It's well-known that for abelian vartieties $X$ with ${\rm End}(X) \otimes \bQ \supset F$
have a moduli space $M_F$, with $M_F \simeq \cH^n / SL_2(\cO_F)$ the
Hilbert modular varieties(\cite{S2}). In Shimura's theory $M_F$ can be also identifies with
the moduli space of abelian varieties $X^{\prime}$ with $\dim X^{\prime}=2n$ and
${\rm End}(X^{\prime}) \supset M_2(F)$, where $M_2(F)$ is the $2 \times 2$ matrix
algebra over $F$, in fact in this case $X^{\prime} \simeq X \times X$. 

\begin{prop}
There is a finite morphism $M_F \longrightarrow M_X$, that is, $M_F$ is the
finite covering space of $M_X$.
\end{prop}
\begin{proof}
We need to recall Shimura's theory of PEL and weak PEL structures(\cite{S2}).
(we only need a simple case, Shimura's theory is far more general).
Let $B$ is a quaternion algebra over $F$, let $X$ be the abelian variety,
and let $\theta: B \hookrightarrow {\rm End}(X)\otimes \bQ$ be an embedding
of algebras, $C$ is a polarization on $X$, then the collection
$\{B,\theta,C\}$ is a PEL structure of $X$. Let $F^+$ be the set of
totally positive elements in $F$, then the collection
$\{B,\theta,\beta C, \forall \beta \in F^+ \}$ is a weak PEL structure of $X$.

In our situation a PEL structure is nothing but a self product of polarized abelian variety
$X \times X$ with ${\rm End}(X)\otimes \bQ \supset F$, and with $B=M_2(F)$. A
weak PEL structure is a self product of $X$ which is a cone polarized abelian
variety.

By Shimura's result(\cite{S2}), if we interpretate the moduli spaces $M_F$ and $M_X$
in this way, then we have a finite covering map: $M_F \longrightarrow M_X$.
\end{proof}

\begin{rmk}
From the infinitesimal point of view it's easy to see that $M_F$ and $M_X$
have the same tangent spaces.

\end{rmk}

\subsection{}

Next we recall some relations between $K$ and its reflex $K^*$. Given a CM
type of $K$, $\{\sigma_i \}_i$, there is a canonical CM type of $K^*$, say
$\{\tau_i \}_i$, and given an ideal $\fa$ of $K$, by transformation theory
we  have a natural ideal $\fb=g(\fa)=\prod_i \fa^{\sigma_i}$ in $K^*$. Thus
given an abelian variety $X=\bC^n / \iota(\fa)$, there is a natural
corresponding $X^*=\bC^n / \iota(g(\fa))$. We easily check that all such
$g(\fa)$ are of the same type, hence all these $X^*$ lie on the same moduli
space $M_n^*$. Now fix the type of $X$, these $X$ with the same type also
lie on a single moduli space $M_n(\fa)$, one may wonder if there is a natural
morphism $M_n(\fa) \ra M_n^*$ such that extend the correspondence $X \ra X^*$.
Such morphism indeed exists, and it can best be seen by the following 
concept.

\begin{Def}
Given a cone polarized abelian variety $X$ of dimension $n$, the primitive
classes of $X$ respective to to the Kahler cone $C_X$ is 
$$\{\alpha \in H^n(X,\bC) \mid \alpha \cdot x=0, \forall x \in C_X \}$$
i.e., the classes of middle degree that annilated by all the elements in the
Kahler cone. Let's denote it as $P_CH^n(X)$.
\end{Def}
This notion is the generalization of the transecendental lattice in the theory of $K3$
surfaces(\cite{PS}).

Let $U_K$, $U_F$ be the units group of $\cO_K$, $\cO_F$ respectively, then
by Dirichlet's unit theorem $U_K / U_F$ is a finite group. Actually
for a generic CM extension of $F$ we have $U_K / U_F=1$, both of the
groups acting on $H^*(X,\bZ)$.

\begin{prop}
$P_CH^n(X)$ is the invariant classes of the action $U_F$ on $H^*(X,\bC)$.
\end{prop}

\begin{proof}
From the definition of $X$ we know we can choose a basis of $H^{0,1}(X)$
as $\{\omega_1, \cdots, \omega_n \}$ such that for any 
$\alpha \in K \simeq {\rm End}(X) \otimes \bQ$, we have 
$\alpha^*(\omega_i) = \alpha^{\sigma_i} \cdot \omega_i$.
Consider the action of $U_F$ on $H^*(X,\bC)$, since $H^*(X,\bC)$ is generated
by $H^1(X,\bC)=H^{0,1}(X)\oplus H^{1,0}(X)$, it's easy to see that
the only invariant classes are those 
$\omega_1 \wedge \cdots \wedge \omega_i \wedge \overline{\omega}_{i+1}\wedge 
\cdots \wedge \overline{\omega}_n$
any $1 \leq i \leq n$. On the other hand any Kahler classes can be written
as the linear combination of
$\{\omega_1 \wedge \overline{\omega}_1,\cdots,\omega_n \wedge \overline{\omega}_n \}$,
hence only those 
$\omega_1 \wedge \cdots \wedge \omega_i \wedge \overline{\omega}_{i+1}\wedge 
\cdots \wedge \overline{\omega}_n$
can be annilated by all elements in $C_X$.
\end{proof}

\begin{corr}
$P_CH^n(X)$ carries a canonical Hodge structure, with 
$P_CH^n(X,\bZ)={\rm Inv}_{U_F}H^n(X,\bZ)$, and
$P_CH^n(X) \simeq P_CH^{0,n} \oplus \cdots \oplus P_CH^{n,0}$, with the
Hodge number $\dim P_CH^{i,n-i}= {i \choose n}=\frac{n!}{i! (n-i)!}$
\end{corr}

\begin{corr}
For $\alpha \in End(X)_{\bQ}=K$, the action of $\alpha$ on the primitive
class $PH^n(X)$ is given by the multiplication of $g(\alpha)$.
\end{corr}

\begin{proof}
This follows from the same proof of the proposition.
\end{proof}

\begin{rmk}
The reason to study the cone polarized variety and primitive classes can be explained 
as the following:
\begin{enumerate}
\item
In the classical theory of Shimura-Taniyama-Weil(\cite{S1}), when considering the notion of 
``field of moduli'' we need to fix a polarization, but it does not specify
any particular choice of the polarization, which implies that any polarization
will give the same field of moduli, so it makes sense to consider all
the possible polarization at the same time, hence this notion of cone
polarized abelian varieties.

\item
Also in the classical theory we need to quotient off the finite automorphism
group of $X$ to get the precise generating functions of the class fields.
Here however the automorphism group of $X$ becomes infinite, so directly dividing
$X$ by the automorphism group doesn't make sense in the category of
algebraic varieties, but it still makes sense to consider the invariant
classes in $H^*(X,\bC)$. $P_CH^n(X)$ can be interpretated as the equivariant
cohomology of $X$ under the action of $U_F$.
\end{enumerate}
\end{rmk}

\subsection{}

Now let's consider the classifying space for the primitive Hodge structures.
As is well known, this classifying space $CL(PH^n)$ is a product of Grassmanians and in general
when the weight is $n \geq 3$, it is much bigger than the geometric moduli
space of the underlying manifolds. One such example is the Calabi-Yau threefolds.
One can also see this by the infinitesimal point of view. Indeed in our
situation the infinitesimal geometric deformation of $X$ is given by
$H^1(X,T_X) \simeq H^{1,n-1}(X)^{\vee}$, in the case of keeping the Kahler cone,
this is given as $PH^{1,n-1}(X)^{\vee}$, which is $n$ dimensional, while the
dimension of $CL(PH^n)$ is much bigger.

But since we already know the moduli space $M_X$, we nevertheless
have a classifying map: $M_X \ra CL(PH^n)$, and we have the
following:
       
\begin{prop}
This map is factoring as $M_X \ra M_n^* \ra CL(PH^n)$, i.e., the image
of $M_X$ in $CL(PH^n)$ is coming from $M_n^*$, where $M_n^* \simeq \cH^n /SL_2(\cO_{F^*})$
is the standard Hilbert modular variety of the reflexive field $F^*$.
\end{prop}

\begin{proof}
The key fact is that for $\alpha \in K=End(X)_{\bQ}$, we have $\alpha \in End(H^1(X))_{\bQ}$,
but $H^n(X)=\wedge^n H^1(X)$, by the above corrolary we have in fact
$g(\alpha) \in K^*$ acting on $PH^n(X)$, i.e. the primitive classes admit
multiplication by $K^*$. Likewise for the ideal $\fa$ of $K$, the
ideal transformation of $\fa : X \ra Y$ can be seen as the ideal
tranformation $H^1(Y) \ra H^1(X)$, which then becomes the ideal transformation
$g(\fa):PH^n(Y) \ra PH^n(X)$. Since $g(\fa)$ is an ideal in $K^*$, and since
$M_n^*$ is the moduli spaces for the corresponding abelian varieties, we must
then have the map $M_n(\fa) \ra M_n^*$.
\end{proof}

Denote $M_{PH}$ as the image of $M_X$ in $CL(PH^n)$, then the above proof shows
we actually have isomorphism $M_n^* \simeq M_{PH}$.

In summary we have shown:

\begin{thm}
\begin{enumerate}
\item
We have a natural isomorphism $M_{PH} \simeq \cH^n /SL_2(\cO_{F^*})$;
\item
The natural morphisms $M_n(\fa) \ra M_C \ra M_{PH}$ are all finite.
\end{enumerate}
\end{thm}

\begin{rmk}
\begin{enumerate}
\item
So from the cohomological point of view the most essential invariants of the
total cohomology of these $X$ are the Kahler cones and primitive Hodge structures.
In the next section we will start to consider the Kahler moduli space.
\item
The morphism $M_n(\fa) \ra M_n^*$ means that we can use the moduli space
$M_n(\fa)$ for the purpose of class field generation of $K$. Indeed in this
view the modular coordinate function of the CM points on $M_n(\fa)$ can be
used for both the class field generations of $K$ and $K^*$, although clearly
there is no direct geometric relations between $X$ and $X^*$, their moduli
points are closely related.
\end{enumerate}
\end{rmk}

\section{\bf Mirror Symmetry for Abelian Varieties}
\subsection{}

Mirror symmetry is usually formulated for the Calabi-Yau manifold, roughly
for a Calabi-Yau manifold of dimension $n$, mirror symmetry predicts there
exists another Calabi-Yau manifold $Y$ of the same dimension such that
we have  natural isomorphisms $H^{i,j}(X) \simeq H^{i,n-j}(Y)$ for any
$i,j$, that is, the Hodge diamond of $Y$ is a rotation from $X$.

Here we will study the case when $X$ is an abelian variety, precisely
let $X$ be as before, ${\rm End}(X) \otimes \bQ \simeq K$, we will try to
find out its mirror partner $Y$.

Let $C_X$ be the kahler cone of $X$, by the last section we know the 
primitive decomposition
$H^*(X)=P_CH^n(X) \oplus C_X \cdot H^*(X)$, now let
$$HC_X=\{\alpha \in H^{i,i}(X), \forall i \mid 
\alpha \ is \ generated \ by \ the \ elements \ in \ C_X \otimes \bC \}$$
i.e., $HC_X$ is the classes generated by the $C_X \otimes \bC$. It can
be shown that $HC_X$ is the invariant cycles in any degenerations
of $X$ on the moduli space $M_X$. We have the formal decomposition
$H^*(X) = P_CH^n(X) \oplus HC_X \oplus C_X \cdot (H^*(X) \setminus HC_X)$.
If we put a formal filtration on $HC_X$ by the degree, we may
regard $HC_X$ carrying a formal Hodge structure, although this
Hodge structure is not coming from the geometry.

Now the mirror partner $Y$ of $X$ is another abelian variety of dimension
$n$ such that we have the natural isomorphisms $H^{i,j}(X) \simeq H^{i,n-j}(Y)$,
further we expect $P_CH^n(X) \simeq HC_Y$,$P_CH^n(Y) \simeq HC_X$  as Hodge structures, 
in particular we should have $P_CH^{1,n-1}(X) \simeq HC_Y^{1,1}$,
$P_CH^{1,n-1}(Y) \simeq HC_X^{1,1}$. From the infinitesimal deformation
of $X$ we know $P_CH^{1,n-1}(X)$ is the tangent space of the complex
moduli space $M_X$, while $HC_X^{1,1}\simeq C_X \otimes \bC$ can be regarded
as the tangent space of the complexified Kahler moduli space, hence mirror
symmetry exchange the complex moduli with the Kahler moduli.

A general framework to construct the Kahler moduli spaces has been suggested 
by Golyshev-Lunts-Orlov(\cite{GLO}), it goes as the following: given our $X$ we can define two
algebraic groups from $H^*(X)$, first let $J_X: H^1(X,\bR) \rightarrow H^1(X,\bR)$
be the complex structureof $X$, $J_X^2=-1$. For any $\omega \in C_X$ we
have the positive definite bilinear form
$$H^1(X,\bR) \times H^1(X,\bR) \longrightarrow \bR,
\  (x,y) \longrightarrow w(x,J_X \cdot y)$$
Consider the morphism of $\bR$ algebraic group $h_X: S^1 \rightarrow GL_{2n}(H^1(X,\bR))$,
with $h_X(e^{i\theta})=\cos\theta \cdot Id + \sin \theta \cdot J_X$.

\begin{Def}
$Hdg{X,\bQ}=$smallest $\bQ$-algebraic subgroup of $GL_{2n}(H^1(X,\bR))$ such that
$h_X(S^1)$ is in and the quadratic form defined by any $\omega \in C_X$ remains
positive definite.
\end{Def}

$Hdg_{X,\bQ}$ is called the Hodge group, note that in our situation it is nothing
but $GL_2(F)$ regarding as an algebraic group over $\bQ$. From it we can
construct the complex moduli space $M_X$.

On the other hand, let $D^b(X)$ be the derived category of coherent sheaves on
$X$, and let ${\rm Auteq}(D^b(X))$ be the group of exact autoequivalence of $D^b(X)$.
There is a natural representation ${\rm Auteq}(D^b(X)) \rightarrow GL(H^*(X,\bZ))$,
let ${\rm Spin}(X)_{\bZ}$ be the image. This is a discrete group whose elements are realizable 
as algebraic correspondences. Let ${\rm Spin}(X)_{\bQ}$ be the Zariski closure of
${\rm Spin}(X)_{\bZ}$ in $GL(H^*(X,\bQ))$. It can be shown that ${\rm Spin}(X)_{\bQ}$ is
in fact a semi-simple algebraic group, whose Lie algebra over $\bC$ is isomorphic
to the Neron-Severi algebra $g_{{\rm NS}}(X)$ of $X$. ($g_{{\rm NS}}(X)$ is 
defined as the following(\cite{LL}):
for each element in $C_X$, by Lefschetz theory we can define a representation
of $sl_2$ on the total cohomology $H^*(X, \bC)$, $g_{{\rm NS}}(X)$ is then the
algebra generated by all such $sl_2$, for all the elements in $C_X$.)

Now let $C_X(\bC)=C_X \otimes \bR + iC_X$ be the complexified Kahler cone, 
the real part $C_X \otimes \bR$ is usually called the B field.
The algebraic $\bR$ group ${\rm Spin}(X)_{\bR}={\rm Spin}(X)_{\bQ} \otimes \bR$ acts on
$C_X(\bC)$ transitively, and for any $\omega \in C_X(\bC)$ the fixed point subgroup
$K_{\omega}$ is a maximal compact subgroup of ${\rm Spin}(X)_{\bR}$. Hence
$C_X(\bC) \simeq {\rm Spin}(X)_{\bR} / K_{\omega}$.

\begin{Def}
The Kahler moduli space $N_X$ is defined as $N_X = C_X(\bC)/ {\rm Spin}(X)_{\bZ}$.
\end{Def}

In this framework the mirror symmetry then is to exchange the two algebraic groups
$Hdg_{X,\bQ}$ and ${\rm Spin}(X)_{\bQ}$, and consequently have the two moduli spaces
$M_X$ and $N_X$ exchanged.

To have a concrete picture of what the Kahler moduli should look like, we need to know more
about the group ${\rm Spin}(X)_{\bZ}$ acting on $C_X(\bC)$.  For this  it's 
instructive to recall the Mumford's description of the Hilbert modular
varieties $M_n(\fb)$. Following Mumford(\cite{AMRT}) we proceed in the following
steps:

\begin{enumerate}
\item
Define $W_d=\{(z_1,z_2,\cdots,z_n) | \Im z_1 \cdot \Im z_2 \cdots \Im z_n \geq d \}
\subset \cH^n$, which is a set in the upper half space,
let 
$\Gamma_2 = \Gamma(\fb) \cap 
\left( \begin{array}{cc}
a & b \\
0 & d
\end{array} \right) $, 
the upper triangler
matrix, and let 
$\Gamma_1=\{\alpha \in \Gamma_2 | \alpha = 
\left( \begin{array}{cc}
1 & a \\
0 & 1
\end{array} \right)\} $, 

the matrix with diagonal being 1. 
Then the action of $\Gamma(\fb)$ on $W_d$ reduced to the action of 
$\Gamma_2$ on $W_d$. We regard $W_d$ as  a neighborhood of infinity, hence
to compactify $M_n(\fb)$ we only need to consider $W_d / \Gamma_2$.

\item
$\Gamma_1 \subset \Gamma_2$ acts on $\cH^n$ as translation by the lattice 
$\iota(\fb)$, hence $\cH^n / \Gamma_1 \subset \bC^n / \iota(\fb) \simeq (\bC^*)^n$.

\item
$\Gamma_2/\Gamma_1 \simeq U_{\fb}$, where $U_{\fb}$ is the unit group
respect to $\fb$, i.e. the unit $u$ satisfying $u \equiv 1 \mod ({\fb})$.By
the Dirichlet theorem we know $U_{\fb} \simeq \bZ^{n-1}$ up to a finite
torsion, and $U_{\fb}$ acts on $\cH^n / \Gamma_1 \simeq (\bC^*)^n$.

\item
take the imaginery part of the tori $(\bC^*)^n \simeq i (\bR^+)^n \times T^n$,
which is the image of the imaginery part of $\cH^n \simeq i (\bR^+)^n \times \bR^n$,
as $\fb$ is a real ideal, embedded as a real lattice. Consider the action
of $U_{\fb}$ on $(\bR^+)^n$, the action is actually the same as
$U_{\fb}$ acts on $F^+ \rah (\bR^+)^n$.

\item
So by the theory of torus embedding, if we can find a rational polyhedral 
decomposition of $(\bR^+)^n$ respect to the action of $U_{\fb}$, we
can form a complete variety $\bar{M_n(\fb)}$.

This decomposition however has already been shown by Shintani(\cite{Shintani1}), indeed he
showed that there is a finite collection of rational polyhehral cones
$V_1, V_2, \cdots, V_{\ell}$ in $(\bR^+)^n$, such that
$(\bR^+)^n= \bigcup_{\epsilon \in U_{\fb}} \bigcup_{i=1}^{\ell} V_i \cdot \epsilon$,
which is a rational simplicial cone decomposition. By the theory of
torus embedding we have a complete toroidal variety, with each
face of dimension $i$ gives rise to a toric orbit of dimension $n-i$.

\end{enumerate}

Now back to the construction of Kahler moduli space, we start from the 
Kahler cone, for the CM abelian variety by Shimura the Kahler cone is
exactly $(\bR^+)^n$, with the action of automorphism group $Aut(X)$, which
is precisely a finite indexed subgroup of the units group $U_F$. Complexifying
the Kahler cone is really adding the real component, and by considering the
derived categoty $D^b(X)$ means we can tensor all the elements of
$NS(X)$, hence the lattice $\bZ^n$ action on the $\cH^n \simeq (\bR^+)^n \times \bR^n$.
So from these consideration we conclude:

\begin{prop}
The Kahler Moduli space $N_X$ is a finite quotient of the Hilbert
modular varity $M_n$.
\end{prop}


Now given our $X$, let $\omega \in C_X$, and let $\Omega$ be the non-vanishing
holomorphic n-form on $X$, how can we construct the mirror partner $Y$?
There is a general method suggested by Strominger-Yau-Zaslow(\cite{Yau}), 
the idea is to use the special
Lagrangian real tori fibration. A real submanifold of real dimension $n$
in $X$ is called Lagrangian if the restriction of $\omega$ on it is zero,
it is special Lagrangian if in addition the $n$-form $\Omega$ restricted
to a volumn form on it. In our situation let $X \rightarrow S$ be a torus
fibration such that each fiber is a special Lagrangian, then the mirror partner $Y$ can be
constructed as the dual torus fibration on the same base $S$. One shows
we can put a natural complex structure and Kahler form on $Y$.

If we realize $X$ as $X=\bC^n / L$ where $L$ is a rank $2n$ lattice,
and with the $\omega$ and $\Omega$ given as above, then the special Lagrangian
construction can be rephrased as the following, we have an exact sequence
$0 \ra L_2 \ra L \ra L_1 \ra 0$ such that $\omega$ restricted to zero on $L_2$ and
$\Omega$ restricted to a non-vanishing $n$ form on $L_2$, here both
$L_1$ and $L_2$ are of rank $n$ lattices. Then $\omega$ defines an 
isomorphism $L_1(\bR) \simeq L_2(\bR)^{\vee}$, 
this defines a natural complex structure on  $L_1(\bR) \oplus L_2(\bR)^{\vee}$
as $J_{\omega}(x,y)=(y,-x)$, moreover the volume form of restriction 
$\Omega$ on $L_2(\bR)$ defines
an isomorphism $L_2(\bR)^{\vee} \simeq \wedge^{n-1} L_2(\bR)$, 
let $L_2^{\vee}\simeq \wedge^{n-1} L_2 \subset L_2(\bR)^{\vee}$ be 
the image lattice, then
$$Y=\frac{L_1(\bR) \oplus L_2(\bR)^{\vee}}{L_1 \oplus L_2^{\vee}}$$

Since $\Omega$ is a $n$ form, it defines a pairing between  $L_1(\bR)$
and $\wedge^{n-1} L_2(\bR)$, hence a pairing between
$L_1(\bR)$ and $L_2(\bR)^{\vee}$, we can check this is a Kahler form
under the complex structure $J_{\omega}$.

\begin{rmk}
All the constructions of Mirror symmetry are more or less pure transecendental, i.e., only
valid over $\bC$, but as we shall see, in our situation it would be far
better to have an algebraic construction, as we will need a mirror
partner defined over a number field.
\end{rmk}

\subsection{}

Now for our $X$ since ${\rm End}(X)_{\bQ} =K$, there is a rank $2n$ lattice $L$ in $K$
such that $X=\bC^n / \iota(L)$, $\omega$ defines a Riemann
form $E(x,y)$ which is given on $L$ as 
$$E(x,y)={\rm Tr}_{K / \bQ}(\zeta x y^{\rho})$$
with $\zeta \in K$ satisfying
$\zeta^{\rho}=-\zeta$, $\Im(\zeta^{\sigma_i}) > 0$, $\forall i$

On the other hand the $n$ form $\Omega$ as a $n$ form on $L$
is given as
$$G(x_1,x_2,\cdots, x_n)={ \rm Tr}_{K^* / \bQ} \det (x_i^{\sigma_j})$$
where $(x_i^{\sigma_j})$ is an $n \times n$ matrix, every term
of its determinent is in $K^*$, so it makes sense to take the trace.

For the lattice $L$ arised in this context the special Lagrangian
fibration  $0 \ra L_2 \ra L \ra L_1 \ra 0$ 
always exists.
Indeed since all our $X$ can be given as $X = \bC^n /(v_n \cdot \iota(\fb) \oplus \iota(\cO_F))$,
we have $\iota(L) \simeq v_n \cdot \iota(\fb) \oplus \iota(\cO_F)$,
then $L_2$ can be taken as the real elements $\cO_F$
in $L$, which is of rank $n$, and one checks that $\omega$ is restricted
to zero on $L_2$, and the form $G$ gives the square root of fundamental discriminent
of $L_2$, in particular, it's non-degenerate.
So the conditions
of being special Lagrangian are satisfied.

Concretely let $X=\bC^n / \iota(L)$ and $0 \ra L_2 \ra L \ra L_1 \ra 0$ be the special
Lagrangian fibration.
Since $\omega$ defines an isomorphism $L_1(\bR) \simeq L_2(\bR)^{\vee}$, 
let $\{e_1,\cdots,e_n\}$ be a basis of $L_1$, let $\{f_1,\cdots,f_n \}$
be the image of $\{e_i\}$ in $L_2(\bR)^{\vee}$, write
$V=L_1(\bR) \oplus L_2(\bR)^{\vee}$, then $\{e_i,f_j\}$ is a basis
of $V$, and the complex structure $J_{\omega}$ of $Y_{\omega}$ 
is defined as $J_{\omega}(e_i)=f_i$, $J_{\omega}(f_i)=-e_i$.
Since $J_{\omega}^2 = -1$, we extend it to $V_{\bC} \rightarrow V_{\bC}$,
let $V^{1,0}$ and $V^{0,1}$ be the eigenspace of $J_{\omega}$ with
the eigenvalue $i$ and $-i$ respectively, then $V_{\bC} = V^{1,0} \oplus V^{0,1}$,
$V^{1,0} = \overline{V^{0,1}}$. Write $\omega_i = e_i +i f_i$,
it's well-known $\{\omega_i \}$ is a
basis of $V^{1,0}$.

If we choose the special Lagrangian fibration  $0 \ra L_2 \ra L \ra L_1 \ra 0$
as above, then since ${\rm End}(X)_{\bQ} \simeq K \supset F$,
the action of $F$ on $L$ keeps $L_2$, hence induces an action on $Y$.

In fact write $L_Y=L_1 \oplus L_2^{\vee}$ as the integral lattice, 
$Y_{\omega}=V / L_Y$, then 
$L_Y \subset L_Y \otimes \bR \subset L_Y \otimes \bC \simeq V_{\bC}$.
If we write the holomorphic forms $\{\omega_i \}$ in terms of the
integral basis of $L_Y$, by the explicit forms of $G$, we see that all
the coefficients are in $F$. Consider the action
of $F \subset {\rm End}(Y_{\omega})_{\bQ}$ on $\{\omega_i \}$, since
$F$ acts on $L_Y$ rationally, $F: L_Y \otimes \bQ \ra L_Y \otimes \bQ$,
we get that $F$ acts on $\{\omega_i \}$ as matrices with entries in $F$. 
Conversely such representation of $F$
determines $Y_{\omega}$ uniquely.

So we have:
\begin{prop}
If we construct the mirror partner $Y$ as above, then $F \subset {\rm End}(Y)_{\bQ}$,
i.e., $Y$ has a $F$ multiplication.
\end{prop}

Since $Y$ has the $F$ multiplication, its automorphism group are of infinite order, and 
it's Kahler cone is $n$ dimensional. From this it's easy to see that if $X$ is algebraic, so
does $Y$, i.e., we always have an integral element in its Kahler cone.

\subsection{}
But for our purpose of generating the class fields, we really need to know the rationality of
the Mirror $Y$, i.e., if $X$ has a $K$ multiplication, we want to be sure it's Mirror $Y$ is defined
over a number field. In so far the above construction shows little about this. What can we do?

Note that in the Mirror construction the complex moduli of $Y$ is given by the complexified Kahler
moduli of $X$. But any given $X$ only has a Kahler class, which is only the imaginery part of
the complexified Kahler class, so we have the freedom of specifying a B-field. We shall show that
an appropriate choice of this B-field will garanttee the rationality of $Y$.

In general for any variety, to study it's rationality we need to construct appropriate algebraic
invariants. If the variety has a good moduli space, then these invariants should be regarded as the
algebraic coordinates on the moduli space, and the rationality of the variety then is given by the
rationality of these coordinates on the moduli points.

For abelian varieties Mumford has developed an algebraic theory of theta functions(\cite{Mumford}),
by this theory the algebraic invariants of an abelian variety is in general given by the ``theta-null'',
i.e., by the theta function on $X$ evaluated at the identity element. So to understand the rationality
of the Mirror abelian variety $Y$, we need to know the theta-null of $Y$. We shall see that for a
good choice of the B-field, the theta function of $X$ and $Y$ are related, hence the theta-null
are ralted as well. So the rationality of $Y$ would follow from the rationality of $X$.

To compare the theta functions on $X$ and $Y$, we need to identify their underlying real manifold, and
regard the Mirror transform as a ``rotation'' of complex structures. For this purpose we first need to
introduce appropriate coordinates.

We begin with $X=\bC^n /(w \cdot \fa \oplus \cO_F) \simeq \bR^{2n}/\bZ^{2n}$, any polarization of $X$
is given by an integral skew-symmetric bilinear form on $\bR^{2n}$. Given such a form $\w$,
we can find an integral basis $\{\lambda_1,\cdots,\lambda_{2n} \}$ of the integral lattice such
that if $\{ x_1, \cdots,x_{2n} \}$ is the dual basis, then $\w=\sum_{i=1}^n \delta_i dx_i \wedge dx_{n+i}$,
with $\delta_1|\delta_2|\cdots $ the elementary divisors.

In our case the biliear form is given by the trace $Tr_{K/\bQ}(\zeta x y ^{\rho})$ with the
admissible $\zeta \in K$ such that $\zeta^{\rho}=-\zeta$, $\Im(\zeta^{\sigma_i}) >0$. Thus we
may regard $(x_1,\cdots,x_n)$ as an integral basis of $\cO_F$, and $(x_{n+1},\cdots,x_{2n})$ as an
integral basis of $\fa$. In particular $(x_{n+1},\cdots,x_{2n})$ depends on $\fa$. In the following we will
denote it as $x_{n+i}(\fa)$.

Next we introduce the complex structure, so $X$ becomes a complex tori, and we can introduce the complex
coordinates. To do this let $e_i=\lambda_i /\delta_i$, $i=1,2,\cdots,n$, and let $\{z_i\}$ be the complex
dual of $\{e_i \}$. Consider the change of coordinates transform:
$$
\Omega \cdot (x_1,\cdots,x_{2n})^T=(z_1,\cdots,z_n)^T
$$

Then $\Omega=(\Delta_{\delta}, Z)$ with $\Delta_{\delta}=diag(\delta_1,\cdots,\delta_n)$ the diagonal
matrix, and $Z$ symmetrical, $\Im(Z)>0$. We recogonize that $\Delta_{\delta}^{-1}Z$ is the period matrix of $X$. 
Indeed
in our case for abelian varieties with CM by $K$, the peiod martix $Z$ is necessarily diagonal
$\Delta_{\delta}^{-1}Z=diag(w_1,\cdots,w_n)$, with $w=(w_1,\cdots,w_n) \in \cH^n$.
(Recall that our moduli space $M_n(\fa)$ is of $n$-dimensional, as opposed to the general
abelian varieties moduli spaces, which is $\frac{n(n+1)}{2}$ dimensional.)

In particular we have $z_i=\delta_i \cdot x_i + \delta_i \cdot w_i \cdot x_{n+i}(\fa)$, we may regard it
as the transform from the underlying real coordinates to the complex coordinates.

Now recall the theta function on $X$ is characterized by the periodic conditions:
$$
\theta(z+\lambda_i)=\theta(z); \ \theta(z+\lambda_{n+i})=e^{-2\pi i z_i} \theta(z)
$$

Taking absolute values we have:
$$
|\theta(z+\lambda_i)|=|\theta(z)|; \ |\theta(z+\lambda_{n+i})|=e^{2\pi \Im(z_i)} |\theta(z)|
$$

From the above coordinates transformation we have:
$$\Im(z_i)=\delta_i \cdot \Im(w_i) \cdot x_{n+i}(\fa)$$

The Mirror symmetry transform
says that we can exchange the complex moduli with the Kahler moduli,
while
the coordinates $w=(w_1,\cdots,w_n) \in \cH^n$ can be regarded as the complex moduli of $X$, 
the Kahler moduli coordinates are actually in the variables
$(x_{n+1}(\fa), \cdots,x_{2n}(\fa))$.
Since they are depend on the type $\fa$, we want to write the dependency explicitly, in 
order to understand the transformation of types.

For this end let's write $x_{n+i}=x_{n+i}(\cO_F)$, these $(x_{n+1},\cdots,x_{2n})$ can be
regarded as a fixed integral basis of $\cO_F$. Introduce a new coordinate $(t_1,\cdots,t_n)$ by
$x_{n+i}(\fa)=t_i \cdot x_{n+i}$, $\forall i$. Since
the type ideal $\fa$ as a $\bZ$ module, is
a submodule of $\cO_F$ of full rank, i.e., if we fix an integral basis of $\cO_F$, then
$\cO_F/{\fa} \simeq \oplus_{i=1}^n \bZ /{t_i \bZ}$, with $\prod_i t_i= N(\fa)/D_F$. Thus
the positive rational numbers $(t_1,\cdots,t_n)$ can be conviniently regarded as the coordinates of the
ideal $\fa$, and under appropriate identification, can be regarded as the coordinates of 
of the Kahler class in the the Kahler moduli.

These $(t_1,\cdots,t_n)$ are only the real coordinates of the Kahler cone, to introduce the complxified
Kahler cone, we need the B-fields, where are they? In our case to make the formula of theta function
consistent we can use the information from the complex moduli $(w_1,\cdots,w_n)$. In fact if we write
$w_j=x_j + i y_j$, then 
$$
\begin{array}{lll}
z_j & = & \delta_j \cdot x_j + \delta_j \cdot t_j \cdot w_j \cdot x_{n+j} \\
    & = & \delta_j \cdot x_j + \delta_j \cdot t_j \cdot (x_j + i y_j ) \cdot x_{n+j} \\
    & = & \delta_j \cdot x_j + \delta_j \cdot y_j \cdot (\frac{t_j \cdot x_j}{y_j} + i t_j ) \cdot x_{n+j} \\
\end{array}
$$

Thus the complexified Kahler coordinate is $(s_1,\cdots,s_n)$ with $s_j=\frac{t_j \cdot x_j}{y_j} + i t_j$, $\forall j$.
If we choose the B-field in this way, the multiplier of the theta function would not change. In the following we
will always use this choice, and for a given Kahler class $\omega$, we will denote the Mirror partner constructed
in this way by $Y_{\omega}$.

But from the above formula, when we exchange the complex moduli $(w_1,\cdots,w_n)$ and Kahler moduli
$(t_1,\cdots,t_n)$,
it's not going to change the multiplier 
$e^{2 \pi \Im(z_i)}=e^{2 \pi\delta_i \cdot t_i \cdot \Im(w_i) \cdot x_{n+i}}$.
Since the absolute value of theta functions can be regarded as  a real analytic function on
$\bR^{2n}$, thus we conclude that for the given Mirror pair $X$ and $X^{\prime}$, their theta functions'
absolute values satisfying the same periodic condition, hence must be only differed by a constant!

Hence by the theory of theta function we have:

\begin{prop}
Let $X$ as before, ${\rm End}(X)_{\bQ} =K$. If $\omega$ is an integral point in $C_X$, 
i.e., it comes from
an ample line bundle on $X$, then $Y_{\omega}$ is defines over
a number field.
\end{prop}

\begin{prop}
The complexified Kahler cone can be identified as 
$$\cH^n \simeq i \fb^+ \otimes \bR \oplus \cO_F \otimes \bR$$
where $\fb^+$ is the set of totally positive elements in $\fb$. Moreover this
identification is compatible with the action of unit groups $U_F$ on both side.
\end{prop}
\begin{proof}
The problem really is to find a an appropriate identification for the Kahler
cone $C_X$. By definition $C_X \subset H^2(X,\bR)$, we know $C_X \simeq (\bR^+)^n$.
But we really want to know the position of $C_X(\bZ) \subset H^2(X,\bZ)$, the 
Neron-Severi cone. Because $X= \bC^n/ \iota(\fa)$, so $H^1(X,\bZ)$ can be identified
with $\fa$ as a $\cO_F$ module. 

Recall that as a $\cO_F$ module we have $\fa \simeq \fb \cdot \w_1 \oplus \cO_F \cdot \w_2$,
with $\fb$ the type of $\fa$. Since $H^2(X,\bZ)\simeq \wedge^2 H^1(X,\bZ)$, so
as a $\cO_F$ module, $\wedge^2 \fa$ is of rank 1, in fact $\wedge^2 \fa \simeq \fb$.
Hence we have as a $\cO_F$ module, $H^2(X,\bZ) \supset \fb$, and the Neron-Severi
cone $C_X(\bZ)$ can be identified as the totally positive elements in $\fb$.

\end{proof}

\begin{prop}
Assume $Y_{\omega}$ is defined over some field $k$. Then if all the elements of
action $F \subset {\rm End}(Y_{\omega})_{\bQ}$ are also defined over $k$, we have
$k \supset F$. Conversely if $k \supset F$, then all then elements of
$F \subset {\rm End}(Y_{\omega})_{\bQ}$ are defined over $k$.
\end{prop}
\begin{proof}
recall that we have constructed a basis of holomorphic differential
forms $\{\omega_i\}$ such that the action of $F \subset {\rm End}(Y_{\omega})_{\bQ}$
on them is a representation over $F$, so if all the elements of
$F \subset {\rm End}(Y_{\omega})_{\bQ}$ are defined over $k$, by taking the trace we
get $k \supset F$. Conversely if $\sigma$ is any automorphism of $\bC$ that keeps
$k$ fixed, then since $k \supset F$, $\sigma$ fixes $F$, hence
fixes all the elements of $F \subset {\rm End}(Y_{\omega})_{\bQ}$.
\end{proof}

\begin{corr}
$Y_{\omega}$ and all the elements of $F \subset {\rm End}(Y_{\omega})_{\bQ}$
are defined over a finite extension of $F$.
\end{corr}

\subsection{}
Now for the set of abelian varieties of CM type $\{K,\sigma_i\}$, we want
to give an explicit Kahler moduli point for each of them.

To do this since we have for each $X$, $X=\bC^n /\iota(\fa)$, let
$X_0=\bC^n / \iota(\cO_K)$, and let $N_X$ be the complexified Kahler moduli
of $X_0$. We need to assigne a point for each $X$ on this space.

Recall all the polarization of any $X$ is given as the Riemann form
$$E(x,y)={\rm Tr}_{K / \bQ}(\zeta x y^{\rho})$$
with $\zeta \in K$ satisfying
$\zeta^{\rho}=-\zeta$, $\Im(\zeta^{\sigma_i}) > 0$, $\forall i$

So let's from now on fix such a $\zeta$, this then defines a polarization
$p_X$ on any $X$. If $X=\bC^n /\iota(\fa)$, then we have an isogeny

$$f: X_0 \ra X$$

Let $f^*(p_X)$ be the pull back of the polarization $p_X$, it's easy to see
that $f^*(p_X)$ is also an polarization on $X_0$, hence defines a point in $C_{X_0}$,
hence a point $q_X \in N_X$.
This will be our Kahler moduli point of $X$ in $N_X$.

\begin{prop}
\begin{enumerate}
\item
If $X_1$ and $X_2$ are of the same type, then there is a positive rational $\xi$
such that $q_{X_1}=\xi \cdot q_{X_2}$;
\item
If $X=\bC^n / \iota(\fa)$, and if $(t_1,\cdots,t_n) \in \cH^n$ is the coordinates
of $q_X$, then $\prod_i t_i =N(\fa)/D_F$.
\end{enumerate}
\end{prop}
\begin{proof}
The first one is obvious. For the second one, note that in our identification
of the complexified Kahler cone, 
$$C_X(\bC) \simeq  i \fb^+ \otimes \bR \oplus \cO_F \otimes \bR$$
with $\fb$ the type of $\cO_K$, which we recall satisfying $\fc^2=D_{K/F}$.
On the other hand the type of $\fa$ is $\fc \cdot \fa \fa^{\rho}$, under
the pull back map $f^*$, $\fc \cdot \fa \fa^{\rho}$ becomes a subset of
$\fc$. By the naturality of the mapping and construction, and by the
definition of Norm, we then have $\prod_i t_i =N(\fa)/D_F$.

\end{proof}
\begin{rmk}
After we define the Kahler moduli point of $X$, since we know the complexified Kahler moduli space $N_X$
is a Hilbert modular variety, the natural question is, what is the nature of these Kahler moduli points? 
Are they CM points?

These Kahelr moduli points are explicitly given in our coordinates as $\{\frac{t_j \cdot x_j}{y_j} + i t_j\}_{j=1}^n$,
they are certailnly some kind of algebraic points, but unlikely to be CM points. Because although $\frac{x_j}{y_j}$ comes from
some CM points, the quantity $t_j$ only depends on the real ideal $\fa$, and we simply put them in the imaginery
axis. There seems no reason to believe that they come from a simple CM extension.
\end{rmk}

\begin{rmk}
This explicit form of the coordinates will be useful when we consider the
Stark's conjecture. Roughly it means that our coordinates of Kahler moduli
behave similarly with the imaginery part of a CM point.
\end{rmk}

\section{\bf Transformations and Isogenies}

\subsection{}

Let $X$ be the abelian variety as before, ${\rm End}(X)_{\bQ} = K$, we know
there exists a rank $2n$ lattice in  $K$ such that $X=\bC^n / \iota(L)$.
Let $X^{\prime}$ be another such abelian variety with 
$X^{\prime}=\bC^n / \iota(L^{\prime})$.
Let $f: X \rightarrow X^{\prime}$ be a surjective morphism, such $f$ can be lifted to
a lattice map $f^*: \iota(L) \rightarrow \iota(L^{\prime})$. Let $\fc$ be an ideal in $K$,
if the above lattice map is induced by an ideal multiplication 
$\fc \cdot L \rightarrow L^{\prime}$, we call $f$ a $\fc$ transformation. Note that in
this case $f$ is an isogeny with 
$\ker(f) = \{ x \in K \mid \fc \cdot x \in L^{\prime} \} / L$.

A particular simple case is when both $L$ and $L^{\prime}$ are ideals in $K$,
in this case if $\gamma \in L^{-1}L^{\prime}$, $\gamma \neq 0$, then the diagonal
$n \times n$ matrix of elements $\{\sigma_1(\gamma), \cdots , \sigma_n(\gamma) \}$ acts
on $\bC^n$, transforms $\iota(L)$ to $\iota(L^{\prime})$, hence every $\gamma$ represents
a $(\gamma (L^{\prime})^{-1}L)$ transform of $X$ to $X^{\prime}$. Conversely every
homomorphism of $X \rightarrow X^{\prime}$ is induced by such $\gamma$. If $h_K$
is the class number of $K$, then there are exactly $h_K$ such abelian varieties,
all the isogenies between them are ideal transformations.

On $X = \bC^n / \iota(L)$, recall that we can define two canonical forms,
the Kahler form $\omega$ and the  $n$ form $\Omega$. $\omega$ is defined
by a Riemann form on $L \subset K$ as
$$E(x,y)={\rm Tr}_{K / \bQ}(\zeta x y^{\rho})$$
with $\zeta \in K$ satisfying:
$$\zeta^{\rho} = -\zeta, \Im (\zeta^{\sigma_i}) > 0, \forall i$$

On the other hand $\Omega$ is defined by a degree $n$ alternating form on $L \subset K$
as:
$$G(x_1,\cdots,x_n)={\rm Tr}_{K^* / \bQ}(\det (x_i^{\sigma_j}))$$

Assume we have two such abelian varieties $X$ and $X^{\prime}$, with the
forms $\omega$, $\Omega$, and $\omega^{\prime}$, $\Omega^{\prime}$ on them
respectively. Assume there is an isogeny $f : X \rightarrow X^{\prime}$,
on the form level the pull-backs $f^*(\omega^{\prime})$ and $f^*(\Omega^{\prime})$
are also Kahler form and $n$ form, so it makes sense to compare
them with $\omega$ and $\Omega$ on $X$. We say $f$ preserves $\omega$
if $f^*(\omega^{\prime})$ is a constant rational multiple of $\omega$, likewise
we define the same for $\Omega$.

\begin{prop}
If $\fc$ is an ideal of $K$ and  $f : X \rightarrow X^{\prime}$ is a $\fc$ transform,
then $f$ preserves $\omega$ if and only if $\fc$ is an imaginery ideal, i.e.,
$\fc=g(\fb)$ for some ideal $\fb$ in $K^*$. Likewise $f$ preserves $\Omega$ if and
only if $\fc$ is a real ideal, i.e, $\fc=\fc^{\rho}$.
\end{prop}

\begin{proof}
If $f$ preserves $\omega$ then by the definition of Riemann form $E$ that defines
$\omega$ we must have $\fc \fc^{\rho}=(\mu)$ for some $\mu \in F$, by the property
of reflexive field this implies $\fc=g(\fb)$ for some idelas $\fb$ in $K^*$.

Similarly if $f$ preserves $\Omega$ then from the defining form $G$ we must
have $g(\fc)=(\nu)$ with $\nu \in K^*$, again by the property of reflexive field
we have $\fc=\fc^{\rho}$.
\end{proof}

\subsection{}

The Kahler form $\omega$ defines a polarization on $X$, and the $n$ form
$\Omega$ defines a polarization on the mirror partner $Y_{\omega}$. Classically
$\omega$ also defines an isogeny $X \rightarrow X^{\vee}$ which is
an ideal $\fa$ transformation, where $\fa=\zeta \fd (L L^{\rho})$, $\fd$ is the different
from $K$ to $\bQ$, and $(L L^{\rho})$ is the smallest ideal that contains
$L L^{\rho}$. In particular such $\fa$ is a real ideal, and it is called the
type of $\omega$.

Similarly for the form $\Omega$ we can define the ideal $\fb=(g(L))$ in $K^*$,
where $(g(L))$ is the smallest ideal in $K^*$ that contains $g(L)$. $\fb$ is
an imaginery ideal in $K^*$, we call it the type of $\Omega$.

\begin{prop}
If $\fc$ is an ideal of $K$ and  $f : X \rightarrow X^{\prime}$ is a $\fc$ transform,
and if $\fa$,$\fb$, $\fa^{\prime}$, $\fb^{\prime}$ are the respective types,
then $\fa^{\prime} = \fa \fc \fc^{\rho}$, 
$\fb^{\prime}=\fb \cdot g(\fc)$.
\end{prop}

\begin{proof}
This is clear since $\fc \cdot L = L^{\prime}$.
\end{proof}

Note since for our $X$ we have ${\rm End}(X)_{\bQ}=K$, the lattice $L \subset K$
is actually an order in $K$. If $X^{\prime}$ is another abelian variety which is
isogenious to $X$, then
we can arrange to make $L^{\prime}$ and $L$ to be
in the same maximal order, so any isogenies between $X$ and $X^{\prime}$ is
an ideal transformation, so we can conclude:

\begin{corr}
Let $f : X \rightarrow X^{\prime}$ be an isogeny. If $f$ keeps the type of
$\omega$, then $f$ is an imaginery ideal transform. If $f$ keeps the type
of $\Omega$, then $f$ is a real ideal transform.
\end{corr}

\begin{prop}
Let $f: X \rightarrow X^{\prime}$ be an isogeny who is induced from the
lattice map $L \rightarrow L^{\prime}$, if the lattice map preserves the special
Lagrangian fibration, then $f$ induces an isogeny 
$Y_{\omega} \rightarrow Y^{\prime}_{\omega}$ between the mirror partners.
\end{prop}

Recall that the special Lagrangian fibration can be chosen as taking $L_2$
as the real sublattice in $L$, in this case the multiplication by the real
ideal will preserve the special Lagrangian fibration, so we have
\begin{corr}
Let $f : X \rightarrow X^{\prime}$ be an isogeny given by the real ideal
multiplication, then $f$ induces an isogeny between the mirror partners.
\end{corr}

\begin{rmk}
In the theory of Shimura-Taniyama-Weil, when considering the class fields
generation, i.e., when considering the Galois action on $X$, the form $\omega$
or the type of $\omega$ needs to be preserved. This is the main reason why
they only get part of the class fields. Our idea roughly can be stated as
to use the type of
$\Omega$ to catch those missing class fields.
\end{rmk}

\section{\bf Congruence Relations and Class Fields Generation}
\subsection{}

Our idea is to use the field of moduli of the mirror partner $Y_{\omega}$
to generate the class field of $K$. Recall that we already know $Y_{\omega}$
is defined over a finite extension of $K$, so let $k$ be such a field of
definition, $k \supset K$. Let $\wp$ be a prime in $K$, $\beta$ be a prime
in $k$ above $\wp$, $q=N(\wp)$. Let $\widetilde{Y}_{\omega}$ be the reduction of
$Y_{\omega}$ by $\beta$, and $\widetilde{Y}_{\omega}^q$ the Frobenius transformation
that takes each points'coordinates to their $q$-th power.

\begin{prop}
Let $X$ and $Y_{\omega}$ as above, $\wp$ a prime in $K$, and $X^{\prime}$ be the
$\wp$ transform of $X$ by regarding $\wp$ as an ideal in $K$, and $Y^{\prime}_{\omega}$ be its
mirror partner. Then
$\widetilde{Y}_{\omega}^q $ is isomorphic to $\widetilde{Y}^{\prime}_{\omega}$.
\end{prop}

\begin{proof}
Let $\sigma$ be a lifting of the Frobenius transform:
$\widetilde{Y}_{\omega} \rightarrow \widetilde{Y}_{\omega}^q$
to $\sigma: k \rightarrow k$ such that $\sigma(\beta)=\beta$, 
$\sigma(x)=x^q (mod \ \beta)$, such lifting always exists. Let $Y_{\omega}^{\sigma}$
be the transform of $Y_{\omega}$ by $\sigma$, and let $X^{\sigma}$ be the mirror of
$Y_{\omega}^{\sigma}$, then since $\sigma$ doesn't change the polarization
of $Y_{\omega}$, then $X$ and $X^{\sigma}$ are of the same type of $\Omega$,
by our theory there is a real ideal $\fa$ transform
$X \rightarrow X^{\sigma}$ that induces $\sigma : Y_{\omega} \rightarrow Y_{\omega}^{\sigma}$.

Let $r=[k_{\beta} : K_{\wp}]$, then $\sigma^r =1$, so by iteration we
have the $\fa^r$-tansform $X \rightarrow X$, this means there is an element
$\gamma \in K$ such that $\fa^r =(\gamma)$, it's easy to see $\gamma$ is actually
real, $\gamma \in F$.

On the other hand since $\sigma: Y_{\omega} \rightarrow Y_{\omega}^{\sigma}$
is a lifting of Fronbenius tranform, the order of the isogeny $\sigma$ then
is $N(\beta)^n$, this then implies the order of the $\gamma$ multiplication
$a^r : X \rightarrow X$ is $N_{K / \bQ}(\gamma)=N(\beta)^{nr}$. In particilar
if $\wp^{\prime}$ is any prime ideal in $\fa$, 
$\wp^{\prime}$ must be dividing $\wp$, further inspection on
the order of $\wp$ in the class group shows that we must have 
$(\gamma)={\wp}^r$, i.e.,
$\fa= \wp$ as ideals in $K$.
\end{proof}

\begin{rmk}
This should be regarded as an analogue of the Kronecker congruence relation.
\end{rmk}

\subsection{}

Recall that for the polarized abelian variety $Y_{\omega}$ the field of moduli is
a field $k_0$ such that if $Y_{\omega}$ is defined over any $k \supset k_0$, and
$\sigma$ be any isomorphism of $k$ into some other field, then
$Y_{\omega} \simeq Y_{\omega}^{\sigma}$ if and only if $\sigma$ restricted on
$k_0$ is identity. Our objective is to show $k_0$ is a class field of $K$.

For this let $X$ as before, $\fb=g(\fa)$ be the type of $\Omega$, and write
$P=[X,\fb]$. If $X$ and $X^{\prime}$ are of the same type of $\Omega$, then
there is a real ideal transformation from $X$ to $X^{\prime}$, so for
$P$ and $P^{\prime}$ if they have the same $\fb$, we write $[P:P^{\prime}]$
as the isogeny from $X$ to $X^{\prime}$, and it's easy to check that
$[P:P]=1$ and all such $[P:P^{\prime}]$ forms an finite abelian group $G$.

On the other hand, as $k_0$ is the field of moduli of $Y_{\omega}$,
if $\sigma \in Gal(k_0 / F)$, since $\sigma$ doesn't
change the polarization of $Y_{\omega}$, let $X^{\sigma}$ be the mirror
partner of $Y_{\omega}^{\sigma}$, put $P^{\sigma}=[X^{\sigma},b]$, then
it makes sense to define $[P^{\sigma}:P]$, in this way we get homomorphism
$Gal(k_0 / K) \rightarrow G$, one checks that this is in fact
an isomorphism. So the field of moduli of $Y_{\w}$ will generate some
class fields of $K$.

To determine precisely which class fields it is generating, we note the following:

\begin{prop}
The class group $C_K$ is generated by the real and imaginery ideals.
\end{prop}

\begin{proof}
First of all note the norm map : $N_{K / F}: C_K \rightarrow C_F$ is surjective.
In fact if $H \subset C_F$ is the image and $k$ is the class field corresponding
to $H$, then it's easy to see that $k \subset K$. But $k$ is an unramified extension
of $F$, hence we must have $k=F$, so $H=C_F$.

Now assume $\fa \in \ker(N_{K / F})$, then $\fa \cdot \fa^{\rho}=(\mu)$, we are
going to show $\fa=g(\fb)$ for some ideals in $K^*$. 

For this let $k$ be a Galois extension of $\bQ$ that contains both $K$ and $K^*$,
let $\{\sigma_i \}$, $\{\tau_i\}$ be the CM-types of $K$ and $K^*$ respectively,
and let $S$, $S^*$ be the sets of complex embeddings of $k$ that induce
$\{\sigma_i \}$ and $\{\tau_i\}$ respectively, it's well-known (\cite{S1}) that
$S^* = \{ g \in Gal(k / \bQ \mid g^{-1} \in S \}$. Let
$H = \{ g \in Gal( k / \bQ) \mid g \cdot S = S \}$, then $H$'s fixed field is
$K^*$, the action of $H$ on $S$ has the decomposition $S=\sum_{i=1}^n \alpha_i H$.

Still using $\fa$ to denote the ideal in $k$ generated by $\fa$, then since
$\fa \cdot \fa^{\rho} = (\mu)$, and since $k$ is Galois, we have
$\fa=\prod_{\sigma \in S} \fc^{\sigma}$, and for any $\sigma \in S$, we have
$\fc^{\sigma} \cap K = \fa$. So we get 
$\fa=\prod_{\sigma \in S} \fc^{\sigma}=\prod_{i=1}^n \prod_{\sigma \in \alpha_i H} 
\fc^{\sigma}$
write ${\fb}_i=\prod_{\sigma \in \alpha_i H} \fc^{\sigma}$, then since $H$ is the
fix subgroup of $K^*$, we have $\fb_i$ are the ideals in $K^*$. Now it's straight
forward to check that $\fb_i = \fb^{\tau_i}$, for some ideal $\fb$ in $K^*$, hence
$\fa=g(\fb)$.
\end{proof}

If $\wp$ is a prime in $K$ and consider the reduction of  $Y_{\omega}$ over
a prime $\beta$ in $k_0$ over $\wp$, by the above congruence relation
the Frobenius transform of 
$\widetilde{Y}_{\omega} \rightarrow \widetilde{Y}_{\omega}^q$ is induced by
the $\wp$ tranform on $X$, hence by class field theory and by the above proposition
we get:

\begin{thm}
Let $H_1$ be the subgroup of class group of $K$ that is generated by
all the imaginery ideals $\fa$, $\fa \cdot \fa^{\rho}=(\mu)$, then the field
of moduli of $Y_{\omega}$ composed with $K$ is the class field corresponding to the
fixed field of $H_1$.
\end{thm}

\subsection{}

Let $\fa$ be a real integral ideal in $K$, consider $\fa$'s action on $X$, since
$\fa$ is real, the action of $\fa$ will preserve the special Lagrangian fibration
of $X$, so in particular $\fa$ acts on the mirror partner $Y_{\omega}$, and we
can define the cross-section $t(\fa)=\{x \in Y_{\omega} \mid \fa \cdot x =0 \}$ on 
$Y_{\omega}$, and there is a field of moduli for $(Y_{\omega}, t(\fa))$, let's
denote it as $k_0(\fa)$. Now analogous to the above theory we want to show
$k_0(\fa)$ is in the ray class fields of $K$ with conductor $\fa$.

If $\lambda: X \rightarrow X^{\prime}$ is a real ideal $\fc$ transform, with $\fc$ is
prime to $\fa$, then the induced isogeny $Y_{\omega} \rightarrow Y_{\omega}^{\prime}$
maps $t(\fa)$ to $t^{\prime}(\fa)$. If $\lambda_1: X \rightarrow X^{\prime}$ is
another isogeny that keeps the special Lagrangian fibration, then there is
a $\mu \in F$ such that $\lambda_1=(\mu)\lambda$, $\lambda_1$ maps $t(\fa)$ to
$t^{\prime}(\fa)$ as well, so $\lambda t(\fa)=t^{\prime}(\fa)$, 
$\lambda \mu t(\fa)=t^{\prime}(\fa)$, this implies 
$(\lambda -\lambda_1)t(\fa)=\lambda (1-\mu)t(\fa)=0$, by the definition of
$t(\fa)$, $(1-\mu)\fc \subset \fa$, hence $\mu \equiv 1 (\mod^{\times}  \fa )$.

So for a fixed integral real ideal $\fa$, it's natural to consider the pair
$\{\fb,\nu\}$ with $\fb$ a real ideal in $K$ prime to $\fa$, and $g(\fb)=(\nu)$ with
$\nu \in K^*$. Given another pair $\{\fb_1,\nu_1\}$ we call
it equivalent to $\{\fb,\nu\}$ if there is a $\mu \in K$ such that
$\fb_1=\mu \fb$, $\nu_1=g(\mu)\nu$, and $\mu \equiv 1 (\mod^{\times}\fa)$.
The equivalent classes of these pairs form a finite abelian  group
$G(K,\fa)$. 

Now by the exactly same argument as in the  unramified case we can show
that $Gal(k_0(\fa) / K) \simeq G(K,\fa)$, and further by the congruence relation
and class field theory we get

\begin{thm}
Let $H_1(\fa)=\{ideals \ \fb \ of \ K, prime \ to \ \fa 
\mid \fb \fb^{\rho}=(\mu), N(\fb)=g(\mu),
\mu \equiv 1 (\mod^{\times} \fa ) \}$,
then $H_1(\fa)$ is a subgroup of ideal group modulo $\fa$, and the field of moduli
$k_0(\fa)$ composed with $K$ is a subfield of ray class field of conductor $\fa$, 
precisely it is the fixed field of $H_1(\fa)$.
\end{thm}

\subsection{}

By Shimura-Taniyama-Weil theory, the field of moduli of $X^*$ is the class field
of $K$ corresponding to the subgroup of real ideals in the ideal class group.
By our previous results, the field of moduli of $X^*$ can be obtained from the
geometric classifying space of primitive classes of weight n, 
hence from the moduli space $M_n(\fa)$. Now combined with
the above theorems, we get:

\begin{thm}
Let $K$ be a primitive CM field, and let $X$ be an abelian variety
of CM type $(K,\{\sigma_i\})$, then there are well-defined complex
moduli point $p_X$ in the geometric moduli space of primitive classes
of middle degree, and Kahler moduli point $q_X$ in the complexified
Kahler moduli, such that together they generate the Hilbert class fields
of $K$.
\end{thm}

Similarly we have the corresponding result for the ray class fields:

\begin{thm}
Let $\fa$ be an integral ideal of $K$. We can define the set of CM
abelian varieties with the $\fa$-level structure as $\{X,\fa \}$,
and from this there are well defined complex moduli points via
primitive Hodge structures $p_{X,\fa}$ and Kahler moduli points
$q_{X,\fa}$, together they generate the ray class fields $K_{\fa}$
of conductor $\fa$.
\end{thm}

{\bf Example}
Let's consider the calssical example, $F=\bQ$. In this case the Kahler moduli
is again the modular curve $\cH/ \Gamma$. By our choice of the B-field,
the Kahler moduli point is inside the set $\{ \zeta \in K | \iota(\Im\zeta) > 0 \}$.
We know that $j$ function whenevaluated over them takes algebraic integer values in the
Hilbert class fields of $K$, hence from generation of class fields, the Kahler moduli
is the same as the complex moduli.  This is consistent with the classical theory
of complex multiplications for the elliptic curves.

\section{\bf Explicit Reciprocity Laws}
\subsection{}

We have constructed the Kahler moduli space $N_X = C_X(\bC) / \Gamma_2$.
let $X$ as before, $Y_{\omega}$ the mirror partner of $X$, then $Y_{\omega}$ defines
a moduli point $q_X \in N_X$, we recall that $q_X$ is given as
$(s_1,\cdots,s_n)$ with $s_j=\frac{t_j \cdot x_j}{y_j} + i t_j$.
The Kahler moduli $N_X$ is actually a Shimura variety, since it has a 
finite covering over $\cH^n /SL_2(\cO_F)$.

By the above theory the field of moduli of $Y_{\omega}$ generates the class field of
$K$. By the general idea of Shimura varieties, such field of moduli should be
given as the values of coordinate functions on $N_X$ at $q_X$, so it is natural
to seek a Shimura type reciprocity laws (\cite{S1}) in this context.

Fot this purpose we  define a extended class of modular functions on the moduli spaces.
Let $\cM_X$ be the set of function that

\begin{itemize}
\item
It is a finite function defined on the set of all abelian varieties
$X$ of CM type $\{K,\sigma_i \}$;
\item
When the type if fixed, it is a restriction of Hilbert modular function;
\item
When the CM point $z$ is fixed, it is a restriction of modular function
on $N_X$.
\end{itemize}

As we shall see later, such functions always exist from the formulism of
Stark's conjecture. We may say that these functions provided an interpolation
between the different moduli spaces $M_n(\fa)$ of various types.

The reason to introduce such class of modular function is that they provide
a natural place for us to formulate the Shimura reciprocity law.

To do this we need to define a natural action of the Adele group of $K$ on this
set of functions. Let $\bA_K$ be the Adele ring of the field $K$, for any
$s \in \bA_K$, to define the action of $s$ on these modular functions, we
need to distinguish two cases, as according to the real and imaginery ideals.

Let $\{K, \sigma_i \}$ be the CM type, and let $\{K^*,\sigma_i^*\}$ be the
CM type of the reflexive field $K^*$. These $\{\sigma_i\}$ can be regarded as
the automorphisms on $\bA_{K}$, so in particular we can define
$g(s)=\prod_i s^{\sigma_i}$ for any $s \in \bA_{K}$. We check that
$g(s)$ is a well-defined element in $\bA_K^*$. Likewise for any $s \in \bA_K$,
let $f(s)=s\cdot s^{\rho}$.

\begin{Def}

Let $s \in \bA_K$. If $s=g(t)$ for some  $t \in \bA_{K^*}$, then $s$ is called
an imaginery adele. If $s=f(t)$, then $s$ is called a real adele.

\end{Def}

\begin{prop}
For any $s \in \bA_K^{\times}$ in the multiplicative group, $s$ can
be written (non-uniquely) as $s=t_1 \cdot t_2$, where $t_1$ is an
imaginery adele, $t_2$ is a real adele.
\end{prop}

This is an reformulation of the proposition in the last section.

When $s$ is an imaginery adele, since $K$ acts on $H^1(X,\bQ)$, so
$\bA_K^{\times}$ acts on $\wedge^n H^1(X,\bQ)\otimes \bA_{\bQ}
\simeq H^n(X,\bQ) \otimes \bA_{\bQ}$. Since $s$ is imaginery, the
action of $s$ is trivial on the Kahler cone, hence the action of $s$
preserves the primitive rational Hodge structures. On the other
hand, the action of $s$ on $H^n(X,\bQ) \otimes \bA_{\bQ}$ is factorized
through $s \ra g(s)$, hence defines an action of $g(s)\in \bA_{K^*}$ on
$H^n(X,\bQ)\otimes \bA_{\bQ}$. As we have seen, then geometric moduli
space of $X$ in the classifying space of primitive classes are the
Hilbert modular variety of $K^*$, hence in the case when $s$ is imaginery
we have a well-defined element $s \ra \tau_1(s) \in GL_2(\bA_F)$, 
$\tau_1(s)$ acts on the modular function of $GL_2(\bA_F)$ by acting on
its argument.

When $s$ is a real adele, $s=t \cdot t^{\rho}$. The action of $t$ on the
on $H^2(X,\bQ)\otimes \bA_{\bQ}$ factorized through $t \ra f(t)=s$
and preserving the Kahler cone, hence $s$ acts on $C_X(\bQ) \otimes \bA_{\bQ}$.
By our construction of the complexified Kahler cone, this action can
be extended to an action on $(iC_X \otimes \bQ^+ + C_X \otimes \bQ) \otimes \bA_{\bQ}$.
The Kahler moduli is constructed by $C_X(\bC) / GL_2(F)$, hence the action
of $s$ defines an element $s \ra \tau_2(s) \in GL_2(\bA_F)$. This element
$\tau_2(s)$ acts on the modular function of $N_X$ by acting on the argument.

Combine these two action and the above proposition we have
\begin{prop}
There is a well-defined action of $\bA_K^{\times}$ on the modular function
space $\cM_X$.
\end{prop}

\subsection{}

Now we can state the reciprocity law.

\begin{thm}
Let $s \in \bA_K^{\times}$ be an adele, $(s^{-1},K)$ be the Artin symbol, then for any
$f \in \cM_X$ that is finite at the modular point $(p_x,q_X)$, we have
$f(p_X,q_X)^{(s^{-1},K)} = f^{\tau(s)}(p_X,q_X)$.
\end{thm}

\begin{proof}
We need to check the consistency of our constructions. But this is
quite straightforward since we construct the action of $\tau$ through
the cohomology groups of $X$, which produces the natural moduli
spaces for the modular function in $\cM_X$ to live.

The rest of the theorem then
follows from the general theory of Shimura varieties. 
First we note that $N_X$ is the geometric moduli space of some
cone polarized abelian varieties, here the cone is the Kahler
cone of $Y_{\omega}$, by mirror symmetry.
So in particular the field of moduli of $Y_{\omega}$ is generated
by the algebraic coordinates of $q_X$.

Next since $N_X$ is a Shimura variety, the natural modular functions
can be used as the algebraic coordinates on $N_X$, as for the sufficient
high weight modular forms define a projective embedding of $N_X$.
So for the modular function $f \in \cM_X$, if $f$ is finite at
$q_X$, then $f(q_X)$ is in the class field of $K$.

Now since by the congruence relation the real ideals of $K$ transform
$q_X$ in the Kahler cone, hence in $N_X$, and since $\bA_K^{\times}$ acts
on $\cM(N_X)$ by transforming the argument of the function,
we get the reciprocity law for real ideals. 

The treatment of imaginery ideals are similar.
\end{proof}

\begin{rmk}
This is a rather formal accessment, as we don't know the detail of
particular modular function $f \in \cM_X$. Later in the Stark's conjecture
we shall see some examples of such functions.
\end{rmk}

\section{Stark's Conjecture For CM Fields}
\subsection{}

Start from 1970 a system of conjectures have been purposed by Stark(\cite{Stark}) as an
approach to the Hilbert's 12th problem, which put emphasis on the structures
of regulators and Dedekind zeta functions. Indeed since Hilbert's 12th problem is
asking for some natural transecendental functions to generate the class fields,
nothing could be more natural than the special functions related to the 
zeta functions, as the later behaves nicely for the class field extensions.
These conjectures can be formulated in a very general
form, for any number fields, but here we will only study a special case.

In this subsection (only)
let $K$ be any number field with $[K:\bQ]=n$, let $r_1(K)$,$r_2(K)$ be the number
of real and complex archimedean places of $K$ respectively, so $r_1(K) + 2 r_2(K)=n$.
By Dirichelet unit theorem, the free part of the unit group $U_K$ has rank
$r(K)=r_1(K) + r_2(K) -1$. Now let $K_1$ be a finite abelian extension of $K$,
$[K_1:K]=m >1$. The unit group $U_{K_1}$ has the rank $r(K_1)$. Let's consider
the problem that how fast of $r(K_1)$ can grow relative to $r(K)$. The reason
for this question is that Stark's conjecture really is about the structure of
regulators, so it makes sense to find the information about the size of the
regulator determinent.

There are two extreme cases:

\begin{enumerate}
\item
If we require $r(K_1)$ to grow slowly and uniformly, then we put
$r(K_1) -r(K)=m-1$. In this case all the relative regulators are size 1
determinent, which is of course the slowest possibility. From the
definition it is easly to see that there are only two possibilities in
this case:

\begin{itemize}
\item
When $r_2(K)=0$, i.e., $K$ is a totally real field and $K_1$ is an extension
such that all but 1 infinities of $K$ are ramified;

\item
When $r_2(K)=1$, i.e.,$K$ has only one complex archimedean prime and there is only
one real infinity unramified in $K_1$.
\end{itemize}

\item 
If we require $r(K_1)$ to grow fast and uniformly, then we put $r(K_1)-m \cdot  r(K)=m-1$.
In this case all the relative regulators are of size n determinent, which is the
fastest growth possibilities. Again there are two cases:

\begin{itemize}

\item
When $K$ is a totally real field, and $K_1$ is also a totally real extension, that
is, $K_1$ is unramified over all the $\infty$ places of $K$.

In this case $r_2(K)=r_2(K_1)=0$,  and $r(K)=n-1$, $r(K_1)=mn-1$.

\item
When $K$ is a CM field, and $K_1$ is  any abelian extension of $K$.

In this case $K_1$ is necessarily CM, since any possible real embedding of
$K_1$ would come from a real embedding of $K$. We have $r_1(K)=r_1(K_1)=0$,
and $r(K)=n-1$, $r(K_1)=mn-1$.

\end{itemize}
\end{enumerate}

Usually the Stark's conjecture is formulated and studied for the first
case, i.e., the slowest case, since the structure of relative regulator is simple
and more importantly can be verified numerically. However Stark's
conjecture can be formulated for all the fields, and in this note we
will study the case when $K$ is a primitive CM field, which is,
according to our classification, one of the fastest growth unit group
case. As we shall see, in this case the problem is actually simpler,
because we can use the regulator of $K$ to simplify the relative
regulators for $K_1$.

\subsection{}

So let $K$ be a CM field, the Dedekind zeta function is defined as
$$\zeta_K(s)=\sum_{a}\frac{1}{N(a)^s}$$
where the sum is over all the integral ideals
of $\cO_K$.
It is well-known that $\zeta_K(s)$ has an analytic continuation to the whole
complex plane, satisfying a standard functional equation, and has a simple
pole at $s=1$. Moreover the Dirichilet class number formula tells us that
$$\zeta_K(s)=\frac{h_K}{s-1}\kappa_K + \rho_K + O(s-1)$$
where $h_K$ is the class
number of $K$, $\kappa_K=\frac{(2\pi)^nR_K}{w_K d_K^{1/2}}$, $R_K$ is the
regulator of $K$, $w_K$ is number of roots of unity in $K$, $d_K$ is the fundamental
discriminent of $K$, and $\rho_K$ is a constant.

Let $H_K$ be the ideal class group of $K$, let $R$ be an ideal class,
we can define the partial zeta function associated to $R$ as the following:
$$\zeta_K(s,R)=\sum_{\fa \in R}\frac{1}{N(\fa)^s}$$

Again $\zeta_K(s,R)$ has an analytic continuation to the whole plane, satisfying
a functional equation between $\zeta_K(s,R)$ and $\zeta_K(s,R^{-1})$, has a simple
pole at $s=1$. We have naturally $\zeta_K(s)=\sum_R \zeta_K(s,R)$, moreover
$$\zeta_K(s,R)=\frac{\kappa_K}{s-1} + \rho_K(R) + O(s-1)$$

So the residue of $\zeta_K(s,R)$ at $s=1$ does not depend on $R$, but the contant
term $\rho_K(R)$ does depend on $R$. This $\rho_K(R)$ is of fundamental importance
for our understanding of the class fields of $K$. When $K$ is quadratic imaginery,
$\rho_K(R)$ is explicitly given by the Kronecker limit formula as special values
of modular functions. These kind of limit formulas will be used in the treatment
of Stark's conjectures.

Now let $K_0$ be the Hilbert class field of $K$, $[K_0 :K]=m$. By
class field theory $Gal(K_0 /K) \simeq H_K$. For any character
$\chi: Gal(K_0 /K) \ra \bC^{\times}$ we can define the Artin
L-function $L_K(s,\chi)=\sum_{R \in H_K} \chi(R) \zeta_K(s,R)$.
It is well-known that $\zeta_K(s)=\prod_{\chi} L_K(s,\chi)$. By
the above Dirichlet's formula, by studying the behavier of $s \ra 1$,
we can have some relations between the regulators of $R_K$, $R_{K_0}$,
and $\rho_K(R)$.

To make the formula simple it is more convenient to consider $s \ra 0$,
and by functional equation this is equivalent to $s \ra 1$. In fact if
we write
$$\zeta_K(s,R)=-\frac{R_K}{w_K} s^{n-1}(1+\delta_K(R)s) + O(s^{n+1})$$
then
$$\delta_K(R)=n \gamma + nlog{2\pi} -log|D_K| -\frac{w_K |D_K|^{1/2}\rho_K(R^{-1})}
{2^n \pi^n R_K}$$

From $L_K(s,\chi)=\sum_R \chi(R)\zeta_K(s,R)$ we have

\begin{itemize}
\item
if $\chi = \chi_0$ is trivial, then 
$$L_K(s,\chi_0)=\zeta_K(s)=-\frac{h_K R_K}{w_K} s^{n-1}(1+\delta_K s) + O(s^{n+1})$$
with $\delta_K = \sum_R \delta_K(R)$

\item
if $\chi \neq \chi_0$ is not trivial, then
$$L_K(s,\chi)=\zeta_K(s,R)=-\frac{R_K}{w_K} (\sum_R \chi(R) \delta_K(R))s^n + O(s^{n+1})$$

\end{itemize}

Since $L_K(s,\chi)=\sum_{R \in H_K} \chi(R) \zeta_K(s,R)$, comparing the leading
coefficient of the both side we have

$$-\frac{h_{K_0}R_{K_0}}{w_{K_0}} = (-1)^m (\frac{R_K}{w_K})^n h_K 
\prod_{\chi \neq \chi_0}(\sum_R \chi(R) \delta_K(R))$$

The Stark's conjecture predicts that if we write $L_K(s,\chi)=R_K(\chi)s^n + O(s^{n+1})$,
then $R_K(\chi)$ also has the form of regulators, i.e. it is a determinent of
a matrix whose entries are linear combination of logarithm of units of $K_0$. 
Moreover we should expect $R_K(\chi)^{\sigma}=R_K(\chi^{\sigma})$ for any
$\sigma \in Aut(\bC)$.

\subsection{}

Now let's take a closer look at the structure of regulators. Given the extension 
$k$ of $K$, note that $k$ is also a CM field. By Dirichilet's unit theorem the
rank of unit group $U_k$ is $mn-1$. Recall the construction of Artin units, let 
$\sigma_1,\sigma_2,..,\sigma_n$ be the real embeddings of $F$ lifted to $K$,
$\sigma_1=id$, let $\epsilon_1,\epsilon_2,..,\epsilon_n \in U_k$ be a
chosen set of fundamental units such that $\epsilon_i=\epsilon_1^{\sigma_i}$.
Let $G=Gal(k/K)$, $G=\{\tau_1,...,\tau_m \}$, the the set
$\{\epsilon_i^{\tau_j} \}_{i=1,...,n}^{j=1,..,m}$ is a system of Artin units,
satisfying one relation:
$$\prod_{j=1}^m \prod_{i=1}^n \epsilon_i^{\tau_j} = 1$$

Consider the logarithm embedding map:
$$\phi: k \mapsto \bC^{mn} \mapsto \bR^{mn}$$
$$\phi(x) = (2log|\phi_1(x)|,..., 2log|\phi_{mn}(x)|)$$
where $(\phi_1,...,\phi_{mn})$ are the infinite places in $k$. The
image of $U_k$ is $L_k=\phi(U_k)$, who lies on a hyperplane $H_k$ defined
as $\sum_{i=1}^{mn} x_i = 0$, $L_k$ is a lattice in $H_k$ of rank $mn-1$.
the regulator $R_k$ is defined as the volume of the fundamental domain
of $L_k$. Inside $L_k$ there is a sublattice $L_k$ generated by
$\{\phi(\prod_{j=1}^m \epsilon_i^{\tau_j})\}_{i=1,...,n}$.
From Galois theory $\prod_{j=1}^m \epsilon_i^{\tau_j}$ is a unit in $K$, hence
$H_K=L_K \otimes \bR \subset H_k$, $dimH_K=n-1$, and $rank(L_K)=n-1$. 
the regulator $R_K$ of $k$ is the volumn of fundamental domain of $L_K$ in $H_K$.

Let's note that in the above definition of the regulator we have have taken
a set of fundamental units, i.e., the generators of the unit group. In general
if we take a set of units that generating a coindex finite subgroup of the
unit group, then the volume of the fundamental domain we got would be a
 rational multiple of the regulator, this is clear from the definition.

Now taken the set of Artin units as above, $\{\epsilon_i^{\tau_j} \}_{i=1,...,n}^{j=1,..,m}$,
from the construction we see that $\eta_i=\prod_{j=1}^m \epsilon_i^{\tau_j}$ for any $i$
is a unit in $K$, and $\xi_j=\prod_{i=1}^n \epsilon_i^{\tau_j}$ any $j$ is
a unit in $k$ satisfying $\xi_j = \xi_1^{\tau_j}$. Since $\prod_{i=1}^{n} \eta_i =1$,
we may take $\eta_1,\cdots,\eta_{n-1}$ as the basis of units in $K$. Similarly
we have the relation $\prod_{j=1}^{m} \xi_j =1$.
In summery we have shown:

\begin{prop}
The unit group of $k$ is generated by the units in $K$ as $\{\eta_i\}_{i=1}^{n-1}$ and
m units in $k$ as $\{\xi_j\}_{j=1}^{m}$ satisfying $\xi_j = \xi_1^{\tau_j}$, up to
a finite index. This is also true for all the extension of $K$.
\end{prop}

Using this new basis of units let's calculate the regulator, 
let $T=L_k/L_K$, then $T \otimes \bR = H_k/H_K$, $rank(T)=mn-n=n(m-1)$.
We may take a basis of $T$ as $\{\xi_j \cdot \eta_i \}$ for $i=1,\cdots,n$ and
$j=1,\cdots,m-1$ where we make the convention that $\eta_n=1$, the problem
then is to calculate the volumn of $T$. Geometrically we may see as the following:
for the subspace $T_j=\{\xi_j \cdot \eta_i \}_{i=1}^n$ the volume is 
naturally $R_K \cdot log|\xi_j|$, and since $T_j$ transformed naturally
under the action of Galois group $Gal(k/K)$, let $S=\{\xi_j \}$ be the
rank $m-1$ lattice, so we have

$$R_k/R_K = vol(T)=R_K^{m-1}\cdot vol(S)$$

where $vol(S)$ is a determinent of a $(m-1) \times (m-1)$ matrix, each
entry is a linear combination of logarithms of $\xi_j$.

\subsection{}

Now let's compare this expression with

$$-\frac{h_{K_0}R_{K_0}}{w_{K_0}} = (-1)^m (\frac{R_K}{w_K})^n h_K 
\prod_{\chi \neq \chi_0}(\sum_R \chi(R) \delta_K(R))$$

since by the Frobenius determinent formula,

$$\prod_{\chi \neq \chi_0}(\sum_R \chi(R) \delta_K(R))=
\det_{R_1,R_2 \neq 1}(\delta_K(R_1)-\delta_K(R_2))$$

hence we should expect :

\begin{enumerate}
\item
$\delta_K(R)=log |\eta_K(R)|$;
\item
$\delta_K(R_1) -\delta_K(R_2)=log |\frac{\eta_K(R_1)}{\eta_K(R_2)}|$,
and $\frac{\eta_K(R_1)}{\eta_K(R_2)}$ is a unit in $K_0$. More
over these units should transform under the Galois group naturally.

\end{enumerate}
   
This is the  form of Stark's conjecture that we are going to study. Let's
note that it is very much the same form as the case of imaginery quadratic
fields.

\section{Eisenstein series for a Totally real fields}
\subsection{}

Our idea of appraoching the  above form of the Stark conjecture is similar
to the known case of imaginery quadratic, namely we first need to find a good formula
of $\delta_K(R)$ by using the Eisenstein series over the real field $F$, then
we determine the integral properties of the singular moduli, and finally
we apply the Shimura reciprocity law. In this section we will recall the
theory of $GL_2$ Eisenstein series over $F$ and try to get a formula for
$\delta_K(R)$, classically this kind of formula is called the Kronecker limit
formula.

The Eisenstein series in question is defined as

$$E(w,s;\fa)=\sum_{(c,d) \in (\fa \oplus \cO_F )/U_F, (c,d) \neq (0)}\prod_{i=1}^{n}
 \Im(w_i)^{s}|c^{\sigma_i}w_i + d^{\sigma_i}|^{-2s}
$$

It is related with the zeta function in the following way, since $\zeta_K(s,R)=\sum_{\fa \in R}\frac{1}{N(\fa)}$,
fix an ideal $\fa_1 \in R^{-1}$. Regarding $\fa_1$ as a module over $\cO_F$,
$\fa_1 \simeq \fa \cdot w_1 + \cO_F \cdot w_2$, where $\fa$ is an fractional ideal
of $\cO_F$, and we choose $w_1, w_2 \in K$ such that $w=\frac{w_1}{w_2}$ satisfying
$\Im w^{\sigma_i} >0$, $\forall i$. Then we have

$$
\begin{array}{lll}
\zeta_K(s,R) & = & N(\fa_1)^s \sum_{\alpha \in \fa_1/U_F, \alpha \neq 0} \frac{1}{N(\alpha)^s}\\
 &=& N(\fa_1)^s \sum_{(c,d) \in (\fa \oplus \cO_F)/U_F,(c,d) \neq (0)} N(cw_1 + d w_2)^{-s} \\
 &=& N(\fa_1)^s \sum_{(c,d) \in (\fa \oplus \cO_F)/U_F,(c,d) \neq (0)}N(w_2)^{-s}\prod_i|c^{\sigma_i}w{\sigma_i} + d^{\sigma_i}|^{-2s} \\
&=& N(\fa_1)^sN(w_2)^{-s}\prod_{i}\Im(w^{\sigma_i})^{-s}E(w,s;\fa)\\
\end{array}
$$ 

The advantage of using Eisenstain series is that it has an explicit Fourier expansion
at infinity of $\cH^n$, from which we can get more information about the basic invariant
$\delta_K(R)$. 

The way to get the Fourier expansion is a standard one (see \cite{Yoshita}), we follow the
following steps:

\begin{enumerate}
\item
First we seperate the defining sum over $(c,d) \in (\fa \oplus \cO_F)/U_F,(c,d) \neq (0)$ into three parts:
$c=0$; $d=0$ and $(c,d) \neq (0)$;

\item 
The sum over $c=0$ or $d=0$ then can be treated by the Dirichlet series of the real field $F$;

\item
To treat the sum over $(c,d) \neq (0)$, we first apply the Mellin transfom with $c$ fixed, then we move the action of 
$U_F$ into the action over the integral domain, and by a change of integral domain we make the ``partial'' theta
function into a full theta function inside the integral.(Exactly as the way Hecke derived the functional equation).
Then we apply the Poisson summation to the full theta function, and convert the integral into a product of
standard Bessel functions.
$$K_s(z)=\frac{1}{2} \int_0^{\infty} exp(-\frac{z}{2}(t+1/t))t^{s-1} dt$$
\end{enumerate}

Applying this method we get(\cite{Yoshita}):
$$
\begin{array}{lll}
E(w,s;\fa)&=& -2^{n-2}h_FR_Fs^{n-1}[1 + (2n\gamma + 2nlog2\pi -2logD_F -\frac{D^{1/2}_F \gamma_F}{2^{n-2}h_FR_F}\\
   &+&   logN(\fa) + log(\prod_i \Im(w_i)) -h(w;\fa))s]+O(s^{n+1})\\
\end{array}
$$

where $\gamma$ is Euler's constant, and $\gamma_F$ is generalized Euler's constant for the field $F$, and
we write $w_i = w^{\sigma_i}$ in the case of no confusion, $w=(w_1,\cdots,w_n)$. And 
$$
h(w;\fa)=\sum_{\chi_F} \chi_F(\fa)h_{\chi_F}(w;\fa)
$$

where 
\begin{itemize}
\item
$\chi_F$ is a character of the ideal class group of $F$;
\item
$\fa$ is the type of $R$;
\item 
$\w$  is the CM point defined by $R$ (with the type $\fa$);
\item
$h_{\chi_F}(\w;\fa)$ is a function on the product of
upper half planes who has the following Fourier expansion:

$$
\begin{array}{lll}
h_{\chi_F}(\w;\fa)&=&\frac{D_F N(\fa)}{2^{n-2}\pi^n h_F R_F}[\chi_F({\fd_F})L_F(2,\chi_F^{-1})\prod_i\Im(\w_i)\\
&+&\pi^nD_F^{-3/2}\sum_{0\neq b \in \fd_F^{-1}\fa} \sigma_{1,\chi}(bda)|N(b)|^{-1}
exp(2\pi i (\sum_{j=1}^n b_j\Re(w_j) + i|b_j\Im(w_j)|))]\\
\end{array}
$$
where
   
\begin{enumerate}
\item
$a,d \in \bA_F^{\times}$ such that $div(a^{-1})=\fa$ and $div(d)=\fd_F$;

\item
$\sigma_{s,\chi}(x)$ is a function defined as
\begin{displaymath}
\sigma_{s,\chi}(s)=\prod_{v\in (finite\ primes)} \left\{ \begin{array}{ll}
1+\chi_v(w_v)q_v^s + \cdots + (\chi_v(w_v)q_v^s)^{ord_v(x_v)} & if \ x_v \in \fd_v,\\
0 & if \ x_v \notin \fd_v.
\end{array} \right.
\end{displaymath}
Here $\fd_v$ denotes the ring of integers of $F_v$, $w_v$ is the prime element of
$F_v$ and $q_v=|\fd_v/w_v \fd_v|$.

\end{enumerate} 
\end{itemize}

So we have the expansion of partial zeta function at $s=0$:
$$
\begin{array}{lll}
\zeta_K(s,R) &=& \frac{R_K}{w_F} s^{n-1}\\
 & & \times \big\{1+(2n\gamma + 2nlog2\pi + nlog2 -logD_F\\
 & & -1/2 log|D_K| -\frac{D_F^{1/2}}{2^{n-2}R_F} \gamma_F(d) + logN(\fa)\\
 & & +log \prod_{i}\Im(w^{\sigma_i}) - \sum_{\chi \in X_F} 
     \chi(\fa)^{-1}h_{\chi}(w;\fa))s \big\} + O(s^{n+1}), \\
\end{array}
$$

But by functional equation we have 
$$\zeta_K(s,R)= -\frac{R_K}{w_K} s^{n-1}(1+ \delta_K(R) s) + O(S^{n+1})$$
where
$$\delta_K(R)=n\gamma + nlog2\pi -log|D_K| - 
\frac{w_K|D_K|^{1/2}\rho_K(R^{-1})}{2^n \pi^n R_K}$$

So we may proceed as before, let $X_K$ be the set of all characters of
$\bA_K^{\times}$ which are trivial on 
$K^{\times}\prod_{v}\cO_v^{\times}K_{\infty}^{\times}$. For $\omega \in X_K$,
put:
$$L(s,\omega)=c(\omega)s^{r(\omega)} + O(s^{r(\omega)+1})$$
as the Hecke L functions, we have $r(\omega)=n$ if $\omega \neq 1$ and
$r(1)=n-1$. We also put:
$$H_K(R)=\sum_{\chi \in X_F}\chi(\fa)^{-1}h_{\chi}(w;\fa) 
- log \prod_{i}\Im(w^{\sigma_i}) -logN(\fa)$$
 
then we have:
$$c(\omega)=\frac{R_K}{w_K} \sum_{R} \omega(R) H_K(R)$$

Note that $H_K(R)$ is actually the same as $\delta_K(R)$ by adding
some universal constant, and by the Frobenius detrminent, replacing
$\delta_K(R)$ by $H_K(R)$ would not affect the statements of the
Stark's conjecture.

We may regard $H_K(R)$ as a function of $w$ defined on the $\cH^n$, whose
value at the CM points $w$ has the arithmetic significence. Indeed,
the Stark conjecture precisely predicts that we should expect
$H_K(R_1 R_2^{-1}) - H_K(R_1)$ is a logarithm of a unit in the Hilbert class
field of $K$. We will demonstrate that this is possible by showing we
have
$$H_K(R)=log|\eta_K(R)|^2$$

for some well defined Hilbert modular forms $\eta_K(R)$ on $M_n{\fa}$
which has good integral
properties. These $\eta_K(R)$ can be regarded as the generalization of the
classical Dedekind $\eta$ function.

\subsection{}
The function $h$ is of central importance to us, let's first note that
by the automorphic property of the Eisenstein series we have the following:

We notice that $h_{\chi_F}$ satisfying the following modular properties(\cite{Yoshita}):
for $\gamma \in \Gamma_{\fa}$ we have
\begin{enumerate}
\item
$h_{\chi_F}(\gamma \w; \fa)=h_{\chi_F}(\w,\fa)$ if $\chi_F$ is not trivial;
\item
$h_1(\gamma \w;\fa) -log(\prod_i \Im(\gamma \w)_i)=h_1(\w,\fa) -log(\prod_i \Im\w_i)$
\end{enumerate}
So in particular we have

$$\sum_{\chi_F}\chi_F(\fa)h_{\chi_F}(\gamma \w;\fa) -log(\prod_i \Im(\gamma \w)_i)
=\sum_{\chi_F}\chi_F(\fa)h_{\chi_F}(\w,\fa) -log(\prod_i \Im\w_i)$$

which suggests that we may actually have a Hilbert modular form 
$\eta_K(w;\fa)$ of parallel weight  
such that $h(w;\fa)=log|\eta_K(w;\fa)|$.
Classically in case $F=\bQ$ this is indeed the case as $\eta_K$ is
the classical Dedekind eta function $\eta$, it is well-known that $\eta$ has
an infinite product expression, which when taking the logarithm translated
into a Fourier expansion, which is exactly the above Fourier series.

In the above expression of $H_K(R)$ and $h_{\chi_F}$ let's note the
apparent symmetrical roles played by $\prod_i \Im(w_i)$ and $N(\fa)$. Recall
that $\fa$ is the type of the ideal $R$, and in our formulism of
the Kahler moduli space, $N(\fa)$ has the meaning of product of the
coordinates of the Kahler moduli point $q_X$. Indeed recall that
if $X=\bC^n / \iota(\fa)$, and if $(s_1,\cdots,s_n) \in \cH^n$ is the coordinates
of $q_X$ with $s_j=\frac{t_jx_j}{y_j}+ it_j$, 
then $\prod_i t_i =N(\fa)/D_F$. This suggests that if we fix
the CM points $\w^{\sigma_i}$ and consider the function $h_{\chi_F}$
as a function of $(s_i,\cdots,s_n)$, we still have the automorphy properties.

Explicitly, recall that $\{t_i\}$ is
the coordinates of Kahler moduli in the Kahler cone $C_X$, and
$(s_1,\cdots,s_n)$ is the complex coordinate of $q_X$
in the complexified $C_X(\bC)$, replacing $t_i$ by $s_i$ means
we make an analytic continuation of $h_{\chi_F}$ to the complexified
$C_X(\bC)$. It's easy to see that for the positive $\Im(w_i)$ the Fourier
series still converges, hence we still have $h_{\chi_F}$ as a function of $(s_1,\cdots,s_n)$ an analytic
function defined on $C_X(\bC)$.

In terms of this coordinates $(s_1,\cdots,s_n)$ we now have:
$$
\begin{array}{lll}
h_{\chi_F}(w;s)&=&\frac{D_F \prod_i \Im(w_i)}{2^{n-2}\pi^n h_F R_F}[\chi({\fd})L_F(2,\chi^{-1})\prod_i t_i\\
&+&
\pi^nD_F^{-3/2}\sum_{(0)\neq (\frac{b_1 \Im(w_1)}{t_1},\cdots,\frac{b_n \Im(w_n)}{t_n}) ( \in \fd^{-1}\fa} \sigma_{1,\chi}(bda)|N(b)|^{-1} \\
& \times &
exp(2\pi i (\sum_{j=1}^n(-i \frac{b_j \Im(w_j)}{t_j}\Re(s_j) + |\frac{b_j \Im(w_j)}{t_j}\Im(s_j)|))]\\
\end{array}
$$

So by the automorphy property of $h_{\chi_F}$ under the transformation of $w$, we
get similar property for $h_{\chi_F}$ under the transformation of $s$:

\begin{prop}
For any $\gamma \in \Gamma_2$,
if $\chi_F \neq 1$ then we have $h_{\chi}(w;\gamma s)=h_{\chi_F}(w;s)$, and if
$\chi_F=1$, we have $h_1(w;s)-log(\prod_i \Im(\gamma s)_i)=h_1(w;w)-log(\prod_i \Im(s_i)$.
\end{prop}

Thus if we exponentiate $h_{\chi_F}$, not only we shall expect a modular form on
$M_n(\fa)$, but we shall expect a modular form on $N_X$ as well. This fits our
view that $N_X$ shall be used to generate the class field of $K$, along with
$M_n(\fa)$.

\subsection{}
The importance of the function $h(w;\fa)$ has been long recognized(see for example \cite{Yoshita}), yet
in general it is proved to be quite difficult to understand. This is in direct contrast with the
classical case $F=\bQ$. In that case the $h$ function becomes (write $w=x+iy$)
$$
\begin{array}{lll}
h(w)&=& \frac{2}{\pi}[\zeta(2)y + \pi \sum_{n \in \bZ,n \neq 0} \frac{1}{n}\sigma_1(n) exp(2\pi i(nx + i |ny|))]\\
 &=& \frac{\pi}{3}y + 4 \sum_{n=1}^{\infty}\frac{\sigma_1(n)}{n} cos(2\pi n x)e^{-2\pi y}\\
 &=& -log|e^{\frac{\pi i w}{12}} \prod_{n=1}^{\infty}(1-e^{2\pi i n w})|^4\\
\end{array}
$$ 

Then we ``recognize'' $e^{\frac{\pi i w}{12}} \prod_{n=1}^{\infty}(1-e^{2\pi i n w})$ is actually the Dedekind's $\eta$-function,
so we can move on by studying the arithmetic properties of $\eta$. Note here the explicit infinite product of $\eta$
is the key. 

Besides the field $\bQ$, however, we are not this lucky. For a general $K$ it seems hopeless to try to find an analogue
infinite product for the $\eta_K(w;\fa)$. Actually this should be as expected: for if we look at the Fourier coefficient
of the $h(w;\fa)$ we see there is a factor $\frac{1}{R_F}$. This is not an algebraic number, so even if
we convert the rest of the Fourier series into a logarithm of an infinite product, we would have the problem of moving
$\frac{1}{R_F}$ to the exponential power, which then would make the infinite product meaningless. So from this respective, 
we should not expect an infinite product expression for $\eta_K(w;\fa)$. This problem, I believe, is the major reason
why the higher dimensional case is difficult.

As explained in the introduction, to overcome this difficulty, we should use a twisted Eisenstein series. The point is
that even we don't know the explicit form of $\eta_K(w;\fa)$, we can still study it's quantitative properties.

The twisted Eisenstein series is defined as: 

$$E_{u,v}(w,s;\fa)=\sum_{(c,d) \in (\cO_F \oplus \fa)/U_F,(c,d) \neq (0)}\prod_{i=1}^{n}
 e^{2\pi i(c^{\sigma_i}u_i + d^{\sigma_i}v_i)}(\Im w_i)^{s}|c^{\sigma_i}w_i + d^{\sigma_i}|^{-2s}
$$

\begin{prop}
We have the limit formula:

$$
E_{u,v}(w,s;\fa)=-2^{n-2}h_FR_F log|g_{-v,u}(w;\fa)|s^n + O(s^{n+1})
$$
where $log|g_{-v,u}(w;\fa)|$ has an explicit Fourier expansion, just like $h(w;\fa)$, 
and $u=(u_1,\cdots,u_n) \in \bR^n$, $v=(v_1,\cdots,v_n) \in \bR^n$, $u,v \notin \cO_F \subset \bR^n$.
In particular the vanishing order of $E_{u,v}(w,s;\fa)$ at $s=0$ is $n$.

Further from the Fourier expansion we have:
$$
\lim_{z\ra 0}\{log|q_w^{\frac{1}{12}}\phi(w,z)| -log(|\eta_K(w;\fa)|^2\prod_i |z_i|) \}=0
$$
that is 
$$
\lim_{z \ra 0}\frac{|q_w^{\frac{1}{12}} \phi(w,z)|}{|\eta_K(w;\fa)|^2\prod_i|z_i|}=1
$$
where we $z=u-vw$, i.e., $z=(z_1,\cdots,z_n)$, $z_i=u_i-v_i w_i$,$\forall i$.

\end{prop}

\begin{proof}
We need to give the Fourier expansion for $E_{u,v}(w,s;\fa)$, and we follow the same method as the
untwisted case (\cite{Yoshita}).

For this let's write
$$
E(w,s;\fa)=\frac{1}{h_F}\sum_{\chi_F} \chi_F(\fa) N(\fa)^{-s}E_{\chi}(w,s;\fa)
$$ 
$$
E_{u,v}(w,s;\fa)=\frac{1}{h_F}\sum_{\chi_F} \chi_F(\fa) N(\fa)^{-s}E_{u,v,\chi}(w,s;\fa)
$$ 

where $\chi_F$ runs through all the characters of the class group of $F$. These $E_{\chi}(w,s;\fa)$ and 
$E_{u,v,\chi}(w,s;\fa)$ are well defined, and they are the similar sum over $(c,d) \in (\fa \oplus \cO_F )/U_F,(c,d) \neq (0)$.

The idea is to split the sum into three pieces: $c=0$,$d=0$,and $c \neq 0, d \neq 0$, which we denoted the corresponding
sum as (I), (II), (III) respectively.

For $E(w,s,\fa)$ these three sums are:
$$
(I)=\chi_F(\fa)^{-1}N(\fa)^s \prod_i \Im(w_i) L_F(2s,\chi_F)
$$

$$
(II)=D_F^{-1}N(\fa)^{1-s} \prod_i (\Im(w_i))^{1-s}(\frac{\pi^{1/2}\Gamma(s-1/2)}{\Gamma(s)})^n L_F(2s-1,\chi_F)
$$

$$
(III)=2^nD_F^{-\frac{1}{2}}(\frac{\pi^s}{\Gamma(s)})^n N(\fa)^{1-s} \prod_i(\Im(w_i))^{\frac{1}{2}}
\sum_{0\neq b \in \fd_F^{-1}\fa}\sigma_{1,\chi}(bda)|N(b)|^{-1}exp(2\pi i(\sum_j(b_j x_j + i|b_j y_j|)))
$$
 where we write $w_j=x_j + i y_j$, and $L_F(s,\chi)$ is the Dirichlet series for the real field $F$.

We note that for the trivial character $\chi_1$, (II) would produce a simple pole at $s=1$, 
this is the only place that we have a pole, hence by functional
equation this will make the vanishing order at $s=0$ to be $n-1$.

We now split the twisted Eisenstein series into similar three sums $(I)^{\prime}$,$(II)^{\prime}$,$(III)^{\prime}$ 
respectively.
Since $E_{u,v}(w,s;\fa)$ is a twist of $E(w,s;\fa)$, we will see the effect of this twist on each of the sums seperately.

\begin{enumerate}
\item
First we note the character $\chi_F$ can be lifted to a character on $\cO_F$, and by the definition of $E_{u,v}(w,s;\fa)$,
$e^{2\pi i(c^{\sigma_i}u_i + d^{\sigma_i}v_i)}$ can be regarded as a continuous character $\chi_{u,v}$ on $\cO_F \oplus \cO_F$,
hence the effect of twisting on (I) is a modification of $\chi_F$: $\chi_F^{\prime}=\chi_F \cdot \chi_{u,v}$,
hence
$$
(I)^{\prime}=\chi_F(\fa)^{-1}N(\fa)^s \prod_i \Im(w_i) L_F(2s,\chi_F^{\prime})
$$

\item
Similarly the effect on (II) is a modification of characters:

$$
(II)^{\prime}=D_F^{-1}N(\fa)^{1-s} \prod_i (\Im(w_i))^{1-s}(\frac{\pi^{1/2}\Gamma(s-1/2)}{\Gamma(s)})^n L_F(2s-1,\chi_F^{\prime})
$$

We note that since $(u,v) \notin (\cO_F \oplus \cO_F)$, the new character $\chi_F^{\prime}$ will never be trivial,
hence the series is regular at $s=0$, since $(I)^{\prime}$,$(III)^{\prime}$  are all regular at $s=1$,
by functional equation, the vanishing order at $s=0$ is $n$. This is
different from the untwisted case.

\item
Finally the effect on the (III) is a shift of variables:
$$
\begin{array}{lll}
(III)^{\prime}&=&2^nD_F^{-\frac{1}{2}}(\frac{\pi^s}{\Gamma(s)})^n N(\fa)^{1-s} \prod_i(\Im(w_i))^{\frac{1}{2}}\\
 & \times &
\sum_{0\neq b \in \fd_F^{-1}\fa}\sigma_{1,\chi}(bda)|N(b)|^{-1}exp(2\pi i(\sum_j(b_j x_j + i|b_j y_j|+ b_j z_j)))\\
\end{array}
$$
\end{enumerate} 

From these expressions we see that if we take $z$ to 0, then $(I)^{\prime} \ra (I)$, and
$(III)^{\prime} \ra (III)$, but for the second sum we need to take out an extra factor:
$$log \prod_i|z_i|$$

i.e., we have
$$
\lim_{z\ra 0}\{log|q_w^{\frac{1}{12}}\phi(w,z)| -log(|\eta_K(w;\fa)|^2\prod_i |z_i|) \}=0
$$ 

\end{proof}

The advantage of using the twisted Eisenstein series is that it is closely related to the theta function
of the ablelian varieties:

\begin{prop}
Regarding $g_{-v,u}(w,z)$ as a function of $(w,z)$ with $z$ defined as above, 
if we define  $\phi(w,z)=q_w^{-\frac{1}{2}B_2(-v)}g_{-v,u}(w)$, 
where $q_w^{-\frac{1}{2}B_2(-v)}=\prod_i q_{w_i}^{-\frac{1}{2}B_2(-v_i)}$,
then  $\phi(w,z)$ is a theta function in the
variable $z$

\end{prop}

\begin{proof}
Recall that $u=(u_1,\cdots,u_n) \in \bR^n$, $v=(v_1,\cdots,v_n) \in \bR^n$, and $\cO_F \subset \bR^n$
as a lattice. From the definition, for any $\alpha \in \cO_F$, we have
$E_{u+\alpha, v}=E_{u, v+\alpha}=E_{u,v}$, i.e., translate invariant under $\cO_F$,
hence $|g_{-v + \alpha,u}|=|g_{-v,u+\alpha}|=|g_{-v,u}|$.

Then from the periodic property of $g_{-v,u}$ we immediately have
$$
|\phi(w,z+\alpha)|=|\phi(w,z)|;\ |\phi(w,z+\alpha w)|=|q_z^{-\alpha}||\phi(w,z)|
$$

Recall that on our abelian variety $X=\bC^n /(w \fa \oplus  \cO_F)$ the theta function is
characterized by 
$$
\theta(z+\alpha)=\theta(z); \ \theta(z+\alpha w)=q_z^{-\alpha}\theta(z)
$$

Since $g_{-v,u}$ is an analytic function, we conclude that $|g_{-v,u}(w,z)|$ as a function
of $z$ is the absolute value of a theta function.

\end{proof}

\begin{corr}
$\eta_K(w;\fa)$ is a derivative theta null.
\end{corr}

\begin{proof}
We only need to notice our theta function is normalized at $z=0$ to be 0, and apply the above propositions.
\end{proof}

By the classical theory of theta function, theta null naturally gives rise to the 
modular forms on the moduli space. In fact this should be  more or less expected.
By Mumford's theory of algebraic theta function, we may further conclude that these
theta nulls in fact defines the moduli space as an integral scheme over $\bZ$, hence
have all the expected integral properties.

In summary we have found the explicit form of the function $\delta_K(R)$.
$$
\begin{array}{lll}
\delta_K(R)&=& logN(\fa) + log\prod_i \Im(w_i) -h(w;\fa)\\
&=& logN(\fa) + log \prod_i \Im(w_i) -log|\eta_K(w;\fa)|^2\\
&=& log[N(\fa)\prod_i \Im(w_i) |\eta_K(w;\fa)^{-2}|]\\
\end{array}
$$

where $R$ is an ideal class of $\cO_K$, $\fa$ its type, $w$ its CM point on $M_n(\fa)$,
and $\eta_k(w;\fa)$ is the theta null:
$$
\frac{\partial}{\partial z_1}\cdots \frac{\partial}{\partial z_n}\phi(w,z)|_{z=0}=\eta_K(w;\fa)^2
$$

\begin{rmk}
In the above identification of the $\eta_K(w;\fa)$ with the theta null, we only use the multiplier
conditions on the theta functions. Since all the theta functions satisfy this condition, and in particular
in our setting the abelian varieties $X$ have large automorphism group $U_F$, so given any theta
function $\phi(w,z)$, the translation of $\phi(w,z)$ under $U_F$ will also be a theta function, thus we
are facing an inevitable question: which theta function shall we use?

From the algebraic point of view, since by Mumford, the theta null give rise to the algebraic coordinates on the
moduli space, so even we have different theta functions, their null values as a function on the moduli space,
must be algebraically related. This means any choice of the theta function would be fine for the rationality
of our $X$. It is in this sense that this identification is consistent with our view that $\eta_K(w;\fa)$ is
a canonical quantity.

Also from this view we see that since any choice of the theta would give the same canonical oordinates, this
caonoical coordinates must be corresponding to the objects that respective to all the theta functions, in other
words, repective to all the polarizations at the same time. We already know these objects, they are the
primitive Hodge structures of $X$ that relative to the Kahler cone. Thus we see that these theta nulls really 
give rise to the algebraic coordinates of this primitive Hodge structures.

\end{rmk}

\section{Rational Properties}
\subsection{}

For an integral ideal $R$ of $\cO_K$ with type $\fa$ and CM point $w$ we define

$$
\eta(R)=\frac{\eta_K(w;\fa)}{N(\fa)\prod_{i=1}^n \Im(w_i)}
$$
and 
$$
\epsilon(R)=\frac{\eta(R)}{\eta(\cO_K)}
$$

we want to show that for two ideals $R_1$ and $R_2$ with different ideal types $\fa_1$,$\fa_2$, and different
CM points $w_1$ ,$w_2$, the quantity $\frac{\eta(R_1)}{\eta(R_2)}$
is always an algebraic number in the class field of $K$. We may actually consider two cases,
one with the type ideal fixed, and $w_1$,$w_2$ lies on the moduli $M_n(\fa)$,
another case with the CM points fixed, and the type ideals are related by
$\fa_1=\fa^2 \cdot \fa_2$.

Recall that for the abelian variety $X_0$ defined by the integer ring $\cO_K$,
the type is $\fc$ with $\fc^2=\fd_{K/F}$. For the imaginery ideal multiplication we
only need to consider the CM points on $M_n(\fc)$. On this moduli space,
let $w_0$ denotes the CM point defined by the ideal $\cO_K$.

Now if $w_1 \in M_n(\fc)$ be another CM point, then there exists a $2 \times 2$
matrix with totally positive integral entries $A \in M_2(\cO_F)$ such that $w_1=A \cdot w_0$.
From the construction we see $N(det(A))=\frac{N(\fa)}{D_K}$ where $\fa$ is the ideal defines
the CM point $w_1$.

The action of $A$ on the CM points can be simplified by considering the action of
$SL_2(\cO_F)$, in fact we have the following:
\begin{lem}
Let $M_2^+(\alpha,\cO_F)$ be the set of $2 \times 2$ matrix with totally positive
entries and determinent $\alpha$. Then $M_2^+(\alpha,\cO_F)/SL_2(\cO_F)$ has the
following repesentatives 
$
\left(\begin{array}{ll}
a & b \\
0 & d \end{array} \right)
$
with $ad=\alpha$, $(a,d)/U_K$,$b \in \cO_F/(d)$.

\end{lem}

In particular, $M_2^+(\alpha,\cO_F)/SL_2(\cO_F)$ is a finite set. So we can 
identify the CM points according to the norm of their defining ideal with
these set.

For real ideal multiplications we have the similar matrix representation as well. Let's
fix an admissible $\zeta \in K$ that defines polarization for all our abelian varieties,
in particular it defines a point $\zeta_0 \in N_{X_0}$ in the complexified Kahler moduli
space of $X_0$. Then by our construction of the Kahler moduli, any real ideal $R_1$ defines
a point $\zeta_1 \in N_{X_0}$ by the isogenies. Since the algebraic group acting on the
Kahler moduli is $Spin(X)_{\bR}$ which is locally isomorphic to $GL_2(F)_{\bR}$, there exists
a $2 \times 2$ matrix with the total positive entries $B \in M_2(\cO_F)$ such that $\zeta_1=B \cdot \zeta_0$.
Matrix $B$ has similar properties as $A$.

\begin{thm}
The number $\epsilon(R)$ is an algebraic number in the class fields of $K$. Moreover
for any ideal class $\fa$ of $K$, if $(\fa^{-1},K)$ is the Artin symbol, then
$\epsilon(R)^{(\fa^{-1},K)}=\epsilon(R \cdot \fa)\epsilon(\fa)^{-1}$.
\end{thm}
\begin{proof}
We only need to consider two extreme cases, when $R$ is an imaginery ideal or a real ideal.

When $R$ is an imaginery ideal, then $R$ is of the same type as $\cO_K$, hence
by definition, $\epsilon(R)=\frac{\eta(R)}{\eta(\cO_K)}$ is the algebraic coordinate function of
the moduli space $M_n(\fc)$ evaluating at a CM point. The set of CM points under the action of
$SL_2(\cO_F)$ can be identified with the quotient space $M_2^+(\alpha,\cO_F)/SL_2(\cO_F)$, so we
can form an auxilliary polynomial
$$
\Phi_{\alpha}=\prod_{A_R \in M_2^+(\alpha,\cO_F)/SL_2(\cO_F)}(x-\epsilon(R))
$$
where $A_R$ is the matrix representative of the CM points defined by the ideal $R$.

By the class fields generation we know $\epsilon(R)$ is an algebraic number, and by the
Shimura reciprocity law, the above polynomial $\Phi_{\alpha}$ is invariant under the
Galois group $Gal(K_0/K)$, hence the coefficients of $\Phi_{\alpha}$ are all in $K$.
Hence $\epsilon(R)$ is in the Hilbert class fields $K_0$. The reciprocity law on $\epsilon(R)$
then follows from the reciprocity law of Shimura on the CM points.

When $R$ is a real ideal, by the theory of Mirror symmetry and construction of Kahler moduli
space, $\epsilon(R)$ is the algebraic coordinate function of the moduli space $N_X$ that evaluated
at the Kahler moduli points. Since we know the similar Shimura reciprocity law on this space,
with the similar matrix representation of the Moduli points, the same arguments as above would apply.

\end{proof}

\begin{rmk}
Here we only treat the Hilbert class field case. Similar ideas can be used to treat more general ray class fields'
Stark conjectures.
\end{rmk}

This theorem can be regarded as a ``rational'' version of Stark's conjecture. In fact Stark's conjecture really
predicts $\epsilon(R)$ is actually an unit, not just an algebraic number. So far we only developed the analytic
theory of the function $\eta_K(w;\fa)$, so for this kind of integral statement we can not say much. We will give
some comments on this matter at the last section.

Compare with the classical case is also interesting. 
Classically when $F=\bQ$ the auxilliary polynomial 
$$
\Phi_{\alpha}=\prod_{A_R \in M_2^+(\alpha,\cO_F)/SL_2(\cO_F)}(x-\epsilon(R))
$$

becomes a $\bZ$ coefficient polynomial, as contrast with our assertion that it's only a $\bQ$ polynomial.
So when $F=\bQ$ we can be sure $\epsilon(R)$ is an algebraic integer, and as we can also consider the
inverse ideal $\epsilon(R^{-1})$ at the same time, we see that it is necessarily an unit.

Why in that case the auxilliary polynomial is of $\bZ$ coefficients? Because there we are dealing
the Dedekind $\eta$ function, who has an infinite product expression, and from which we see $\eta$ has
a Fourier expansion at infinity with integral coefficients. This fact plus the integral property of the
elliptic modular function $j$ then garantteed that auxilliary polynomial is a $\bZ$ polynomial
(see \cite{Lang}).

The idea of Fourier expansion can be generalized formally into our setting as the following:
give the totally real field $F$, let
$\overline{M_n(\fa)}=M_n(\fa) \cup_j D_j$ be the Mumford compactification of $M_n(\fa)$,
by adding smooth normal crossing divisors $D_j$ at the infinity. It is well
known that the smooth projective variety $\overline{M_n(\fa)}$ has a model over
$Spec(\cO_F)$ as a regular scheme, and the boundary divisors $D_j$
are all defines over $Spec(\cO_F)$, each irreducible component has multiplicity 1
over $Spec(\cO_F)$ (see \cite{Faltings}).

The CM points on $M_n(\fa)$ then is rational over
the class fields of $K$, and defines a regular section of 
$Spec(K_0) \times_{Spec(\cO_F)}\overline{M_n(\fa)}$, where $K_0$ is the Hilbert class field of $K$.
The integral property we
need can be summerized as the following:

\begin{prop}
Let $f$ and $g$ be the Hilbert modular forms on $M_n(\fa)$ of the same weight $m$,
assume that both have integral Fourier coefficients, and near the boundary divisors
both have integral coefficients expansion. Then for any $A \in M_2^+(\alpha,\cO_F)/SL_2(\cO_F)$,
consider the function $\phi(z)=|det(A)|^{-m}\frac{f(Az)}{g(z)}$, the value
$\phi(z_0)$ then is an algebraic integer. In other words, $\phi$ is in the
integral ring of the function fields of $M_n(\fa)$, regarding as a scheme over
$Spec(K_0)$.
\end{prop}

We may regard this as a formal analogue of the classical q-expansion princeples.

Now we can see why in our case it is difficult to argue along this line, precisely because we don't
know the exact form of the function $\eta_K(w;\fa)$, let alone its Fourier coefficients at infinity.
It is this problem that causes our lack of understanding of the integral properties.

\section{H function and Singular Theta Lifting}
\subsection{}

We want to understand more about the Hilbert modular form 
$\eta_K(w;\fa)$ of parallel weight 2 
defined by 
$$\sum_{\chi_F} \chi_F(\fa)h_{\chi_F}(\w;\fa)=log|\eta_K(w;\fa)|$$

What we know is the Fourier expansion of $h_{\chi_F}$ near the infinity:
$$
\begin{array}{lll}
h_{\chi_F}(w;\fa)&=&\frac{D_F N(\fa)}{2^{n-2}\pi^n h_F R_F}[\chi({\fd})L_F(2,\chi^{-1})\prod_i \Im(w_i)\\
&+&\pi^nD_F^{-3/2}\sum_{0\neq b \in \fd^{-1}\fa} \sigma_{1,\chi}(bda)|N(b)|^{-1}
exp(2\pi i (\sum_{j=1}^n b_j\Re(w_j) + i|b_j\Im(w_j)|))]\\
\end{array}
$$

$\eta_K(w;\fa)$ can be considered as a generalization of
the classical Dedekind eta function $\eta$, who has
an infinite product expansion

$$\eta(w)=e^{2\pi i w/24}\prod_{n=1}^{\infty}(1-e^{2\pi i n w})$$

which when taking the logarithm translated
into a Fourier expansion of $log|\eta(w)|^2$,  which by some direct
inspection we find it is exactly the above Fourier series.

As we have explained before, 
when $K$ is of higher degree we can not expect $\eta_K(w;\fa)$ has a similar
infinite product expansion, which makes its study quite difficult. Since the Dedekind
$\eta$ is such a classical object, we may ask, what kind of other properties of $\eta$
can be generalized to $\eta_K(w;\fa)$?

In this section we will show that the singular theta lifting representation of $\eta$ indeed
can be generalized to $\eta_K(w;\fa)$.

By the singular theta lifting representation of $\eta$ we mean the following integral representation:

$$log|\eta(w)|^2=\int_{\cH/{\Gamma}} \theta(\tau,w) \theta_1(\tau) y \frac{dxdy}{y^2}$$

where $\theta(\tau,w)$ is the Siegel theta function on SO(1,2), and
$\theta_1(q)=\sum_{-\infty}^{\infty}q^{n^2}$ is the classical Jacobi theta
function.

This is called  a singular theta lifting because the integral on the right hand
side is divergent, so we need to regulate it by the following way, consider

$$\int_{\cH/{\Gamma}} \theta(\tau,w) \theta_1(\tau) y \frac{dxdy}{y^{(2+s)}}$$

as a function of the complex variable $s$, for $\Re(s) >> 0$ this integral
is convergent and is analytic respect to $s$. So we can analytically continue
this integral to the whole $s$ plane, and then we define the
regulated integral as the constant term at $s=0$. 
In the following we always understand the divergent integral
in this sense.

If we replace $\theta_1(q)$ by any other 1/2 weight modular forms $G$, the regulated
integral still makes sense. In fact
singular lifting gives a correspondence between the 1/2 weight modular forms to
the weight $c_0(G)/2$ modular forms, where $c_0(G)$ is the constant Fourier
ceofficient of $G$,
one can show that this correspondence is rational,
preserving the Hecke eigenforms, and L-functions.

The advantage of the theta lifting is that we can explicitly calculate the
Fourier coefficients of the lifting from the Fourier coefficients of $G$.
Moreover as observed by Harvey-Moore(\cite{Harvey}), we can read off the singularities of
the lifting more or less directly. Hopefully this can be used to study the geometrical property of the
form $\eta_K(w;\fa)$ near the boundary divisors of the moduli space.

\subsection{}

We want to generalized this integral lifting to the general $GL_2(F)$, which
means for any given Hilbert modular form of parallel weight 1/2, we want to
define a similar singular integral lifting to another Hilbert modular form of
parallel weight 1/2. For this purpose we need to define an appropriate kernel
function.

The Siegel theta function $\theta(w,\tau)$ is a modular form defined by
the group $SO(1,2)$. Recall that $SO(1,2)$ can be defined as the linear
tranformation group on $\bR^3$ such that leaves the quadratic form
$x_1 x_3 -x_2^2$ invariant. There is an isomorphism $SO(1,2) \simeq SL_2(\bR)$,
which can be seen as the following: for any $g \in SL_2(\bR)$, let
$g$ acts on $\bR^3$ as $g(x_1,x_2,x_3)=g 
\left( \begin{array}{cc}
x_1 & x_2 \\
x_2 & x_3
\end{array} \right)g^T $, 
it is obvious that this action keeps the quadratic form, hence defines a
map $SL_2(\bR) \ra SO(1,2)$, which is an isomorphism.

So the upper half plane $\cH$ is the symmetric space of $SO(1,2)$ as well,
in particular, $\cH$ can be identified as a open set in the Grassmanian
$Gr(3,2)$ such that the two-dimensional plane is negative definite. For
any $w \in \cH$, let $w_+$ and $w_-$ be the projection of $\bR^3$ to this
plane $w$ and its orthogonal complement (which is then positive definite)
respectively, then the Siegel theta function is defined as:

$$
\theta(w,\tau)=\sum_{x \in \bZ^3} exp(2\pi i (\frac{(w_+(x),w_+(x))}{2} \cdot \tau
+ \frac{(w_-(x),w_-(x))}{2} \cdot \bar{\tau})
$$

One can verify that it is invariant under $SO(1,2)$ and transformed under
$SL_2(\bR)$ on the variable $\tau$ as a modular form of weight $(1,1/2)$.

Now in our current situation the modular forms are defined on the product
of upper half planes $\cH^n$. $\cH^n$ is the symmetric space of the group
$SL_2(\bR)^n \simeq SO(1,2)^n$, so in appearence the generalization is straight
forward. However since our arithmetic group is $SL_2(\cO_F)$, we need to
be careful on the integral lattice in order to get the modular form
respect to $SL_2(\cO_F)$.

By the archemedean primes of $F$ we have a natural embedding $\cO_F \ra \bR^n$,
let $L_F \subset \bR^n$ be the image, which is a rank n lattice.
So $L_F^3 \subset (\bR^n)^3 \simeq (\bR^3)^n$ is a rank 3n lattice in $(\bR^3)^n$.
We can conveniently write the elements of this lattice $L_F^3$ as
$(x_1,x_2,x_3) \in \cO_F^3$, which represents the real lattice point by
the real embeddings. Moreover the action of $SL_2(\cO_F)$ on $(x_1,x_2,x_3)$ is given as before:

$$g(x_1,x_2,x_3)=g 
\left( \begin{array}{cc}
x_1 & x_2 \\
x_2 & x_3
\end{array} \right)g^T $$

If we define the $\bR^n$-valued quadratic form on $\bR^3$ as $x_1x_3-x_2^2$, then this action keeps
the quadratic form. We write the inner product under this quadratic form as

$$(x,y)=((x,y)^{\sigma_1},\cdots,(x,y)^{\sigma_n}) \in \bR^n$$

Likewise for any integral ideal $\fa$ of $F$, we identify $\fa$ with its image
in $\iota(\fa) \subset \bR^n$, and we denot $\fa^3$ as the lattice in $\bR^{3n}$.

If $x=(x_1,x_2,\cdots,x_n) \in (\bR^+)^n$ and $w=(w_1,\cdots,w_n)\in \cH^n$, then
we use the dot product to denote
$$x \cdot w = x_1w_1+x_2w_2 +\cdots + x_nw_n$$

Now for any $w=(w_1,\cdots,w_n) \in \cH^n$, $w$ defines a two dimensional planes
in each of the component of $(\bR^3)^n$ which is negatively definite, hence
we can define the projection $w_+=((w_1)_+,\cdots,(w_n)_+)$ and
$w_-=((w_1)_-,\cdots,(w_n)_-)$ for $x=(x_1,\cdots,x_n) \in (\bR^3)^n$, just
as before.

Now for any integral ideal $\fa$ we can define the generalization of 
the Siegel theta function as

$$
\theta_{\fa}(w,\tau)=\sum_{x \in \fa^3} exp(2\pi i (\frac{w_+(x)^2}{2} \cdot \tau
+ \frac{w_-(x)^2}{2} \cdot \bar{\tau})
$$
Here both $w$ and $\tau$ are in $\cH^n$, the dot product with $\tau$ is the
dot product of two vectors in the above sense.

\begin{prop}
$\theta_{\fa}(w,\tau)$ is invariant under $SO(1,2)^n$ respect to $w$, and
transformed under $SL_2(\cO_F)$ on $\tau$ is a Hilbert modular form
of parallel weight $(1,1/2)$.
\end{prop}
\begin{proof}
The idea is exactly the same as the classical case. The invariance under
$SO(1,2)^n$ is easy to see. To check the modular property of $SL_2(\cO_F)$
we use the generators of $SL_2(\cO_F)$. 

Write $\Gamma=SL_2(\cO_F)$, and
let 
$\Gamma_2 = \Gamma \cap 
\left( \begin{array}{cc}
a & b \\
0 & d
\end{array} \right) $, 
the upper triangler
matrix, and let 
$\Gamma_1=\{\alpha \in \Gamma_2 | \alpha = 
\left( \begin{array}{cc}
1 & a \\
0 & 1
\end{array} \right) \} $, 
the matrix with diagonal being 1. 

\begin{lem}
The group $SL_2(\cO_F)$ can be generated by three types of elements:
\begin{itemize}
\item 
Elements in $\Gamma_1=
\left( \begin{array}{cc}
1 & a \\
0 & 1
\end{array} \right)$ with $a \in \cO_F$.

\item
Elements in $\Gamma_2 / \Gamma_1 =
\left( \begin{array}{cc}
u^n & 0 \\
0 & u^{-n}
\end{array} \right) $ with $u \in \cO_F$ is a unit.

\item
The inverse element $
\left( \begin{array}{cc}
0 & 1 \\
-1 & 0
\end{array} \right) $ 
\end{itemize}
\end{lem}

By this lemma we just need to check case by case.

\begin{enumerate}
\item
if $g=
\left( \begin{array}{cc}
1 & a \\
0 & 1
\end{array} \right) $ with $a \in \cO_F$. In this case 
$$
\begin{array}{lll}
\theta_{\fa}(w,g(\tau))& = &
\sum_{x \in \fa^3} exp(2\pi i (\frac{(w_+(x),w_+(x))}{2} \cdot (\tau+a)
+ \frac{(w_-(x),w_-(x))}{2} \cdot (\bar{\tau} + a))) \\
                   & = & 
\sum_{x \in \fa^3} exp(2\pi i (\frac{(w_+(x),w_+(x))}{2} \cdot \tau + (\frac{(w_+(x),w_+(x))}{2} \cdot a
+ \frac{(w_-(x),w_-(x))}{2} \cdot \bar{\tau }+\frac{(w_-(x),w_-(x))}{2} \cdot \bar{a} ))
\end{array}
$$

But by the definition of the dot product and by the definition of the lattice $\fa$, for any
$x$,  we see $(\frac{(w_+(x),w_+(x))}{2} \cdot a$ can be written as $Tr_{F/{\bQ}}(a^{\prime})$
with $a^{\prime} \in \cO_F$, hence it is a rational integer, hence can be dropped from the
expression. We then have $\theta_{\fa}(w,g(\tau))=\theta_{\fa}(w,\tau)$.

\item
When $g=
\left( \begin{array}{cc}
u^n & 0 \\
0 & u^{-n}
\end{array} \right) $ with $u \in \cO_F$ is a unit.

In this case the action of $g$ on $\tau$ is $g(\tau)=\tau \cdot u^{2n}$, but then
$u^{2n}$ can be absorbed into $(\frac{(w_+(x),w_+(x))}{2}$ and $\frac{(w_-(x),w_-(x))}{2}$,
hence into the natural action of $u$ on the lattice $\fa^3$. But $u$ acts on $\fa^3$ as
an automorphism, hence the action of $g$ on the sum is just rearrange the terms, hence
in this case $\theta_{\fa}(w,g(\tau))=\theta_{\fa}(w,\tau)$.

\item
When $g=
\left( \begin{array}{cc}
0 & 1 \\
-1 & 0
\end{array} \right) $ the inverse element. In this case $g(\tau)=- 1/{\tau}$,
so we need a Poisson summation formula:

\begin{lem}
The Fourier tranform of $\tau^{1/2}\bar{\tau}exp[2\pi i(-\frac{w_+(x)^2}{2}\cdot
\frac{1}{\tau} - \frac{w_-(x)^2}{2} \cdot \frac{1}{\bar{\tau}})]$ is 
$exp[2 \pi i(\frac{w_+(x)^2}{2} \cdot \tau + \frac{w_-(x)^2}{2} \cdot \bar{\tau})]$,
and the Possion summation formula is
$$
\sum_{x \in \fa^3}f(x)=|D_F N(\fa)^2|^{3/2} \sum_{x \in ((\fa \fd_F)^{-1})^3} \hat{f}(x)
$$
\end{lem}

By this lemma, we have
$$
\begin{array}{lll}
\theta_{\fa}(w,g(\tau))& = &
\sum_{x \in \fa^3} exp(2\pi i (\frac{w_+(x)^2}{2} \cdot (\frac{-1}{\tau})
+ \frac{w_-(x)^2}{2} \cdot (\frac{-1}{\bar{\tau}})) \\
                   & = & 
|D_F N(\fa)^2|^{3/2}(\prod_i \tau^i)^{-1/2}\cdot (\prod_i\bar{\tau^i})^{-1}\cdot 
\sum_{x \in ((\fa\fd_F)^{-1})^3} 
exp(2\pi i (\frac{w_+(x)^2}{2} \cdot \tau
+ \frac{w_-(x)^2}{2} \cdot \bar{\tau } )) \\
                   & = &
(\prod_i \tau^i)^{-1/2}\cdot (\prod_i\bar{\tau^i})^{-1}\cdot 
\sum_{x \in (\fa)^3} 
exp(2\pi i (\frac{w_+(x)^2}{2} \cdot \tau
+ \frac{w_-(x)^2}{2} \cdot \bar{\tau } )) \\

\end{array}
$$
Hence $\theta_{\fa}(w,\tau)$ is a modular form of weight $(1,1/2)$ on the variable $\tau$.

\end{enumerate}
\end{proof}

Now for any variable $z=(z_1,z_2,\cdots,z_n)$ we use $z^{\bf a}$ to denote the product
$$
z^{\bf a}= \prod_{i=1}^n z_i
$$
Likewise we use $dz^{\bf a}$ to denote the product integration measure
$$
dz^{\bf a}= \prod_{i=1}^n dz_i
$$

Now let $G$ be a Hilbert modular form on $\cH^n$ of parallel weight 1/2, the
theta lifting then is defined as :

\begin{Def}
$$
T(G)=C\int_{\cH^n/{SL_2(\cO_F)}} \theta_F(w,\tau)\cdot G \cdot y^{\bf a} 
\frac{dx^{\bf a}dy^{\bf a}}{y^{2{\bf a}}}
$$
\end{Def}
where $C=\frac{1}{2^nD_FR_F}$ is the constant to make sure the volumn
under this integration measure of the fundamental domain $\cH^n/SL_2(\cO_F)$ is 1, and we write
$\tau = x + i y$.

Again this integral is divergent and we need to regulate it by consider
the analytic continuation of 
$$
\int_{\cH^n/{SL_2(\cO_F)}} \theta_{\fa}(w,\tau)\cdot G \cdot y^{\bf a} 
\frac{dx^{\bf a}dy^{\bf a}}{y^{(2+s){\bf a}}}
$$

from $\Re(s) >>0$, and we define the regulated integral as the constant
term of this at $s=0$.

\subsection{}
One of the advantage of the singular theta lifting is that we can calculate
the Fourier coefficients of the lifting by the Fourier coefficients of $G$. The idea is the
classical Rankin-Selberg method, basically we try to change the integration domain by rearrange
the series, and by applying the Poisson summation over a sublattice. In the end the integral is
changed to an integral over the region $\bR/{\bZ} \times \bR^+$, with the variable $x$ and $y$ seperated, which then
reduced to the classical Bessel's integral, which we can explicitly evaluate.

The only different feature in our case is that in the process of changing the integration domain 
we need to consider the action of unit group. So in order to reduce the integral to $\bR^n/{\cO_F} \times (\bR^+)^n$,
we need to make some restrictions on the form $G$. 
For this let's introduce some notation, for any Hilbert modular form $G$, let
$G(w)=\sum_{\lambda \in \cO_F}c(\lambda)exp(2 \pi i \lambda \cdot w)$ be its Fourier
expansion. We call $G$ uniform under the action $U_F$ if its coefficients
satisfying $c(\lambda)=c(u \lambda)$ for any $u \in U_F$.

\begin{thm}
Write $\tau=x+iy$, assume $G$ is a uniform  holomorphic form  with the Fourier expansion
$G(\tau)=\sum_{k \in \cO_F}\sum_{n \in \cO_F} c(n,k)exp(2 \pi i n \cdot \tau)y^{-k}$, then
the lifting $T(G)$ is well defined and is 
$$
\frac{1}{2^nD_FR_F}\sum_{(\alpha,\lambda)/{\approx}}c(\lambda^2/2)\frac{1}{N(\alpha)}
exp[2 \pi i \alpha \cdot (\lambda,w)]
$$
where the equivalence relation ``$\approx$'' is defined for the pair $(\alpha,\lambda) \in \cO_F \oplus \cO_F$
as
$(\alpha_1,\lambda_1) \approx (\alpha_2,\lambda_2)$ if and only if there is a unit $u \ U_F$ such that
$(\alpha_1,\lambda_1) = (\alpha_2 u^2,\lambda_2/u^2)$.

\end{thm}
\begin{proof}
The existence of the regulated integral followed from some general theory
(see for example \cite{Borcherds}).

To show the formula above, we need to explicit evaluate the regulated
integral, we will follow the standard method of Rankin-Selberg,
as in the case $F=\bQ$, which is given in Borcherds' paper(\cite{Borcherds}), until to the
point that we need to consider the units in $\cO_F$.

{\it Step 1}\\

Consider a null vector $\lambda_0=(x_1,x_2,x_3) \in \fa^3$,
such that $\lambda_0^2=0$, such $\lambda_0$
always exists. Then we make it into a null sublattice by
$V_0=\{\alpha \cdot \lambda_0, \alpha \in \cO_F \}$, i.e., by multiplying
$\cO_F$. Likewsie we choose a $\lambda_1$ such that $(\lambda_0, \lambda_1)=1$,
and define $V_1$ in the same manner. Finally we put $V_2 = \fa^3 /(V_0 \oplus V_1)$,
and we fix an uncanonical decomposition $L_F^3=V_0 \oplus V_1 \oplus V_2$.

{\it Step 2}\\

The kernel $\theta_{\fa}(w,\tau)$ is a sum over the lattice $\fa^3$, we 
are going to rewrite it
as $\theta_{\fa}=\sum_{x \in V_1 \oplus V_2} \sum_{\alpha \in V_0}g(x, \alpha)$,
and we take a partial Fourier transform in the direction of $V_0$. 

By definition, $g(x, \alpha)$ is then:
$$
g(x,\alpha)=exp[2\pi i(\frac{\lambda_0^2 \alpha^2}{2}\cdot (\tau -\bar{\tau})
+\alpha(x,w_+(\lambda_0))  \cdot \tau + \alpha(x,w_-(\lambda_0)) \cdot \bar{\tau}
+\frac{(w_+(x)^2 \cdot \tau + w_-(x)^2 \cdot \bar{\tau})}{2})]
$$

The Fourier tranform of $g(x,\alpha)$ in the direction of $\alpha$ is 
$$
\begin{array}{lll}
\hat{g}(x,\alpha)&=&(2^ny^{\bf a}N(w_+(\lambda_0))^2)^{-1/2}exp[2 \pi i(
\frac{w_+(x)^2}{2} \cdot \tau + \frac{w_-(x)^2}{2}
 \cdot \bar{\tau} \\
                 &  &
-\alpha(x,\frac{w_+(\lambda_0)
-w_-(\lambda_0)}{2w_+(\lambda_0)^2}) - \frac{1}{y} \cdot \frac{((x,\lambda_0)\tau + \alpha)^2}
{4i w_+(\lambda_0)^2})]\\
\end{array}
$$

The Poisson summation gives:
$$
\sum_{\alpha \in \cO_F}g(x,\alpha)=D_F^{1/2}\sum_{\alpha \in \fd^{-1}}\hat{g}(x,\alpha)
$$

Now write $d=\alpha$, and any element of $\fa^3/{V_0}$ can be written as
$c \lambda_1 + x$ with $x \in V_2$, then by the Poisson summation formula
we have
$$
\theta_{\fa}(z,\tau)=(2^ny^{\bf a}N(z_+(\lambda_0))^2)^{-1/2}\sum_{c,d \in \cO_F}
exp[2\pi i(\frac{1}{y} \cdot \frac{((x,c\tau + d)^2}
{4i z_+(\lambda_0)^2} + \lambda_1^2 \cdot \frac{cd}{2})]
\sum_{x \in V_2}exp[2\pi i (x^2 \cdot \frac{\bar{\tau}}{2})]
$$

{\it Step 3}\\

Introduce an equivalence relation $\sim$ for the pair $(c,d) \in \cO_F$ as
$(c,d) \sim (c_1,d_1)$ if there exists $\alpha$ such that $\alpha(c,d)=(c_1,d_1)$.
The equivalence classes under this relation can be identified with the
quotient $SL_2(\cO_F)/{\Gamma_2}$, where we recall $\Gamma_2$ is the matrix with
the form
$\left( \begin{array}{cc}
a & b \\
0 & d
\end{array} \right)$

Since $\cH^n/SL_2(\cO_F) \simeq (\cH^n /\Gamma_2)/(SL_2(\cO_F)/\Gamma_2)$, we can
use the tranformation property of the theta function $\theta_{\fa}$ to change
the integral domain to $\cH^n /\Gamma_2$ , this is the key idea of Rankin-Selberg.

So the theta lifting now is:
$$
\begin{array}{lll}
T(G,s)&=&
C\int_{\cH^n/SL_2(\cO_F)}(2^ny^{\bf a}N(w_+(\lambda_0))^2)^{-1/2}
\sum_{c,d \in \cO_F}
exp[2\pi i(\frac{1}{y} \cdot \frac{((x,c\tau + d)^2}
{4i w_+(\lambda_0)^2} + \lambda_1^2 \cdot cd/2)]\\
      &  & \times
\sum_{x \in V_2}exp[2\pi i (x^2 \cdot \bar{\tau}/2)]G(\tau) 
y^{\bf a}\frac{dx^{\bf a}dy^{\bf a}}{y^{2{\bf a}+s}}\\
      &=&
C\int_{\cH^n/SL_2(\cO_F)}(2^ny^{\bf a}N(w_+(\lambda_0))^2)^{-1/2}\sum_{\alpha \in \cO_F}
\sum_{(c,d)/{\sim} }
exp[2\pi i(\frac{1}{\Im(\frac{a\tau + b}{c\tau+d})} \cdot \frac{\pi \alpha^2}
{4i w_+(\lambda_0)^2} )]\\
      &  &  \times 
\sum_{x \in V_2}exp[2\pi i (x^2 \cdot \bar{\tau}/2)]G(\tau) 
y^{\bf a}\frac{dx^{\bf a}dy^{\bf a}}{y^{2{\bf a}+s}}\\
      &=& 
C\int_{\cH^n/\Gamma_2}(2^ny^{\bf a}N(w_+(\lambda_0))^2)^{-1/2}\sum_{\alpha \in \cO_F}
exp[2\pi i(\frac{1}{y} \cdot \frac{\pi \alpha^2}
{4i w_+(\lambda_0)^2} )]\\
      &  &  \times
\sum_{x \in V_2}exp[2\pi i (x^2 \cdot \bar{\tau}/2)]G(\tau) 
y^{\bf a}\frac{dx^{\bf a}dy^{\bf a}}{y^{2{\bf a}+s}}
\end{array}
$$

{\it step 4}\\

Now write $\cH^n \simeq \bR^n + i (\bR^+)^n$, the action of $\Gamma_2$ on $\cH^n$
becomes:
$$
\cH^n/\Gamma_1 \simeq \bR^n/\cO_F \oplus (\bR^+)^n
$$ 
with $\Gamma_2/\Gamma_1=\{
\left(\begin{array}{ll}
u^n & 0\\
0 & u^{-n}
\end{array}\right) \}$

In order to change the integral into $\cH^n/\Gamma_1$ we follow a similar
idea of Rankin-Selberg, this time about ``turning on'' the units group actions.

First we note that $V_2$ can be identified with $\cO_F$. We introduce another
equivalence relation $\approx$ for the pair $(\alpha,x) \in \cO_F$ as
$(\alpha,x) \approx (\alpha_1,x_1)$ iff $(\alpha/u^2, xu^2)=(\alpha_1,x_1)$
and we change the above sum
$$\sum_{(\alpha,x) \in \cO_F}$$
into
$$\sum_{(\alpha,x)/{\approx}}\sum_{u^2 \in U_F^2}$$

We note the expression doesn't change, hence we have
$$
\int_{\cH^n/\Gamma_1}(2^ny^{\bf a}N(w_+(\lambda_0))^2)^{-1/2}\sum_{(\alpha,x)/{\approx}}
exp[2\pi i(\frac{1}{y} \cdot \frac{\pi \alpha^2}
{4i w_+(\lambda_0)^2} )]
exp[2\pi i (x^2 \cdot \bar{\tau}/2)]G(\tau) 
y^{\bf a}\frac{dx^{\bf a}dy^{\bf a}}{y^{2{\bf a}+s}}
$$

{\it Step 5}\\

so we can seperate the integral into the x and y part. We integrate over x first.
Since
$G(\tau)=\sum_{k \in \cO_K} \sum_{n \in \cO_F} c(n,k)exp(2 \pi i \tau)y^{-k}$,
we integrate over x term by term, and it is easy to see that the x-integral
is zero unless $n=x^2/2$, so we  have
$$
\sum_k \sum_{(\alpha,\lambda)/{\approx}}\int_{(\bR^+)^n}
exp(-\frac{1}{y} \cdot \frac{\pi \alpha^2}{2w_+(x_0)^2} - 2\pi y \cdot \lambda^2)
c(\lambda^2/2,k)y^{-3/2-k-s}dy^{\bf a}
$$
Note here we need the assumption that our $G$ is uniform, otherwise the
rearrangment of the coefficients will not make sense.

The y-integral then is reduced to a product of the classical Bessel's integral
$$
K_s(a,b)=\int_0^{\infty} exp(-\frac{b^2}{t} - a^2t)t^s dt=(\frac{b}{a})^s K_s(ab)
$$
with
$$
K_s(c)=\int_0^{\infty} exp[-c(t + 1/t)] t^s \frac{dt}{t}
$$

So we finally get
$$
\frac{1}{2^n D_FR_F}\sum_{(\alpha,\lambda)/{\approx}}c(\lambda^2/2)\frac{1}{N(\alpha)}
exp[2 \pi i \alpha \cdot (\lambda,w)]
$$
\end{proof}

\subsection{}
The advantage of the integral representaion of like the theta lifting we
just defined is that we can read off the singularities of the lifting
more or less directly from the Fourier coefficents of $G$. This is first
observed by Harvey-Moore(\cite{Harvey}).

In our case however we need to make the following modifications. Since we 
want to know the local behavier of a function we need to have a good
local coordinates, and on the Hilbert modular varieties near the infinity
divisor the good local coordinates are provided by the Toroidal coordinates.

\begin{lem}
For real $r$ the function
$$
f(r)=\int_1^{\infty}e^{-r^2y}y^{s-1}dy=|r|^{-2s}\Gamma(s,r^2)
$$
has a singularity at $r=0$ of type $|r|^{-2s}\Gamma(s)$ unless $s$ is a non-positive
integer, in which case $f$ has a singularity of type $(-1)^{s+1}r^{-2s}log(r^2)/(-s)!$.
\end{lem}
For the proof see the paper of Borcherds(\cite{Borcherds}).

\begin{thm}
Let $G=\sum_{n,k}exp[2 \pi in\tau]y^{-k}$ be a Hilbert modular form of parallel
weight 1/2, and
$T(G)$ be the theta lifting, then $T(G)$ has the singularity of type
$$
\sum_{n,k}-c(\lambda^2/2,k)(-2\pi w_+(x_0)^2)^{1/2+k}log(w_+(x_0)^2)/(1/2+k)!
$$
\end{thm}
\begin{proof}
$T(G)$ is defined by the integral
$$
T(G)=C\int_{\cH^n/{SL_2(\cO_F)}} \theta_{\fa}(w,\tau)\cdot G \cdot y^{\bf a} 
\frac{dx^{\bf a}dy^{\bf a}}{y^{2{\bf a}}}
$$

Over any compact region the integral is convergent and defines an analytic
function of $w$, so the singularity of $T(G)$ is coming from the integral over
the neighborhood of the infinity divisors. Such a neighborhood can be constructed
as the following 
$W_d=\{(\tau_1,\tau_2,\cdots,\tau_n) | \prod_i \Im(\tau_i) \geq d \}
\subset \cH^n$

so we are reduced to the consideration of the integral over $W_1/SL_2(\cO_F)$,
but we know that $W_1/SL_2(\cO_F) \simeq Y_1/U_F^2 \oplus \bR^n/{\cO_F}$,
where 
$$Y_1=\{(y_1,y_2,\cdots,y_n) | y_1y_2 \cdots y_n \geq 1 \}\subset (\bR^+)^n$$

So the integral is 
$$
\int_{Y_1/U_F^2}\int_{\bR^n/{\cO_F}} \theta_{\fa}(z,\tau)\cdot G \cdot y^{\bf a} 
\frac{dx^{\bf a}dy^{\bf a}}{y^{2{\bf a}}}
$$

Substitute the expression of $G$ into the integral and carry over the x-integral
we have
$$
T(G)=\sum_{\lambda,k}c(\lambda^2/2,k)\int_{Y_1/U_F^2}exp(-2\pi y \cdot w_+(\lambda)^2)
y^{-2-k+1/2}dy^{\bf a}
$$

Now let's make a change of variable by putting $\prod_1^n y_i=t$, let
$$G=\{(y_1,y_2,\cdots,y_n) | y_1y_2 \cdots y_n = 1 \}\subset (\bR^+)^n$$ 
The unit group $U_F^2$ acts on $G$, and keeps $t$.
Let $G_0$ be the fundamental domain of $U_F^2$ action on $G$, and let
$d^*c$ be the invariant measure on $G$. Under this coordinates we have
$Y_1/U_F^2 \simeq [1,\infty] \times G_0$. 

So the y-integral becomes
$$
\int_1^{\infty}exp(-2 \pi t^{1/n} c \cdot w_+(\lambda)^2) t^{-3/2-k}dtd^*c
$$

Note however that we can not seperate the integral over the product $Y_1/U_F^2 \simeq [1,\infty] \times G_0$
into a product of
two integrals, because the variables are not seperated. But
the integral over $G_0$ is a finite integral, hence has no contribution to the
singularities. So as long as the singularities are concerned, the integral is reduced to

$$
\int_1^{\infty}exp(-2 \pi N(w_+(\lambda)^2))t^{-3/2-k}dt
$$
Now apply the above lemma we get the conclusion.
\end{proof}

\begin{rmk}
The non-seperateness of above integral over $[1,\infty] \times G_0$ actually reflects the basic
difficulty about this type of theta lifting, it means we only know ``roughly'' how the singularities
of the lifting are. The liftings are analytic functions coming from a Hilbert modular form, and as
the Hilbert modular varieties have a explicit compactification(\cite{AMRT}), it would be far better
if we can have a precise knowledge about the singularities of the lifting at the natural boundary.
\end{rmk}

\subsection{}

For $z=(z_1,\cdots,z_n) \in \cH^n$, let $q=(q_1,\cdots,q_n)$ with
$q_i=exp(2 \pi i z_i)$, and for $x \in \bR^n$ let
$(x,x)=(x_1^2,\cdots,x_n^2) \in \bR^n$, and let
$q^x= \prod_{i=1}^n q_i^{x_i}$.

Now let's define 
$$\theta_{\fa}(q)=\sum_{x \in \fa}q^{(x,x)}$$
this can be regarded as the generalization of the Jacobi theta function.

\begin{prop}
$\theta_{\fa}(q)$ is a Hilbert modular form of parallel weight 1/2.
\end{prop}
\begin{proof}
The method is exactly the same as we check the modular properties of Siegel theta functions.
\end{proof}

\begin{thm}
The theta lifting of $\theta_{\fa}(q)$,
$$
T(\theta_{\fa})=C\int_{\cH^n/{SL_2(\cO_F)}} \theta_{\fa}(w,\tau)\cdot \theta_{\fa} \cdot y^{\bf a} 
\frac{dx^{\bf a}dy^{\bf a}}{y^{2{\bf a}}}
$$
is the sum:
$$\sum_{\chi_F} \chi_F(\fa)h_{\chi_F}(w,\fa)$$
\end{thm}
\begin{proof}
We note the form $\theta_{\fa}(q)$ is a uniform form in our terminology,
hence it can be explicitly evaluated as
$$
\frac{1}{2^n D_FR_F}\sum_{(\alpha,\lambda)/{\approx}} \frac{1}{N(\alpha)}
exp[2 \pi i \alpha \cdot (\lambda,w)]
$$
write $\lambda \alpha=\beta \in \fa$, the sum becomes
$$
\frac{1}{2^n D_FR_F}\sum_{\beta \in \fa} \frac{1}{N(\beta)}
\sum_{\{\lambda | \beta\}/{\sim}}{N(\lambda)}
exp[2 \pi i \beta \cdot w]
$$

On the other hand we have the Fourier expansion of $h_{\chi_F}$:
$$
\begin{array}{lll}
h_{\chi_F}(w;\fa)&=&\frac{D_F N(\fa)}{2^{n-2}\pi^n h_F R_F}[\chi_F({\fd})L_F(2,\chi_F^{-1})y^{\bf a}\\
&+&\pi^nD_F^{-3/2}\sum_{0\neq b \in \fd^{-1}\fa} \sigma_{1,\chi_F}(bda)|N(b)|^{-1}
exp(2\pi i (\sum_{j=1}^n b_jx_j + i|b_jy_j|))]\\
\end{array}
$$
hence for $\sum_{\chi_F} \chi_F(\fa)h_{\chi_F}(w,\fa)$ we have for the first term

$$\sum_{\chi}\chi(\fa \fd)L_F(2,\chi^{-1})y^{\bf a}=\zeta_F(2,\fa \fd)y^{\bf a}$$

the second term is 
$$\pi^nD_F^{-3/2}\sum_{0\neq b \in \fd^{-1}\fa} \sum_{\chi}\chi(\fa)\sigma_{1,\chi}(bda)
|N(b)|^{-1}
exp(2\pi i b \cdot w)
$$

From the definition of $\sigma_{1,\chi_F}$ we know
$$ \sum_{\chi_F}\chi(\fa)\sigma_{1,\chi_F}(bda)= \sum_{\{\lambda | \beta\}/{\sim}}{N(\lambda)}$$

Compare the two expressions we get the theorem.
\end{proof}

\begin{thm}
There is a Hilbert modular form $\eta_K(w;\fa)$ such that
$\sum_{\chi_F} \chi_F(\fa)h_{\chi_F}(w,\fa) = log |\eta_K(w;\fa)|^2$.
\end{thm}
\begin{proof}
By above theorem, $\sum_{\chi_F} \chi_F(\fa)h_{\chi_F}(w,\fa)$ is the theta lifting of
$\theta_{\fa}(q)$, and we know that the singularity of $T(\theta_{\fa})$ is of the
type
$$log(N(w)^2)$$
with rational coefficients. 

On the other hand, on the variety $\cH^n /SL_2(\cO_F)$, or on its toroidal
compactification, let $\Omega$ be its cotangent sheaf. Then there exists
a section of $\wedge^n \Omega$ such that near the boundary divisors it is of the
singularity type $N(w)^2$. We can see such a section actually exists on
$\cH^n$ as $\prod_1^n w_idw_i$. We check this section is invariant under
the action of $SL_2(\cO_F)$, hence decents onto the vareity $\cH^n /SL_2(\cO_F)$.
Let this section be $\eta_K(w;\fa)$. 

Compare the type of singularity of both  $T(\theta_{\fa})$ and $log|\eta_K(w;\fa)|^2$,
and note they are both analytic functions, so we must have
$\sum_{\chi} \chi(\fa)h_{\chi}(w,\fa) = log |\eta_K(w;\fa)|^2$.

\end{proof}

\begin{rmk}
This characterization of $\eta_K(w;\fa)$ again has the problem of indirect. We still
don't know a canonical description of this object.
\end{rmk}

\subsection{}
It is interesting to try to look at the singularities of the $h$ function directly at the 
boundary. Since by (\cite{AMRT}) we have an explicit description of the infinite neighborhood
by the toroidal method, we may expect the toroidal boundary divisors have a simple defining
function, hence in the neighborhood of the infinity, we might be able to understand more
about $h$ function by studying its Fourier coefficients. The situation however is quite subtle,
exactly because of the unit group action.

In the standard toric geometry setting, we have the algebraic torus $(\bC^{\times})^n$ with the
character group $N \simeq \bZ^n$, let $V \subset \bR^n \simeq N \otimes \bR$ be a convex rational
cone, these data defines an affine toric variety $T_V$. Each 1-dimentional boundary $\ell \in V$
defines a boundary toric divisor $D_{\ell}$. If $\alpha \in \bZ^n \subset \bR^n$, let $X^{\alpha} \in N$
be the corresponding character. Then if we regard $X^{\alpha}$ as a function on $T_V$, it is
a rational function with the vanishing order on $D_{\ell}$ equals $<\ell, \alpha>$, where $<,>$ is
the standard inner product on $\bR^n$. In fact $\{X^{\alpha} | \alpha \in V \}$ generates the
affine coordinate ring of $T_V$.

In our case we have 
$$
\cH^n/{\Gamma_2}\simeq \cH^n /\cO_F \simeq \bR^n /\cO_F \oplus i(\bR^+)^n \simeq (\Delta^*)^n \subset
(\bC^{\times})^n
$$
with $U_F^+$ acts equivariantly, where $U_F^+$ is the group of positive units. Note $U_F^+$ acts
on the character group $N$ as well, hence acts on $(\bR^+)^n \subset \bR^n$ as well.

This means to compactify the infinite neighborhood $(\Delta^*)^n/U_F^+$ we need to have an equivariant
cone decomposition $(\bR^+)^n \simeq \coprod_{u \in U_F^+}u(V)$ of the rational convex cones. Such
decomposition of course has been shown by Shintani(\cite{Shintani1}). So we can get an equivariant
toric embedding, and after quotient off $U_F^+$, we get a toroidal boundary for the infinite neighborhood.

Now let's consider the problem of how to construct a rational function for this toroidal boundary. As before
for $\alpha \in \bZ^n \simeq \cO_F$ we denote $X^{\alpha}$ as the character function. This function however
can not be a right function for the toroidal boundary because it is not equivariant. To get a equivariant
function it seems the natural choice is
$$
Y^{\alpha}=\sum_{u \in U_F^+} X^{u\alpha}
$$

$Y^{\alpha}$ is invariant under $U_F^+$, so decents to a function on the toroidal neighborhood. Moreover
for any 1-dimensional boundary $\ell \in V$ we have a well defined vanishing order for each $X^{u \alpha}$,
it is $<\ell,u\alpha>$. So for $Y^{\alpha}$ the vanishing order for $D_{\ell}$ is the minimum
$min\{<\ell, u\alpha>|u \in U_F^+ \}$, we see that this is well defined. In fact for the boundary divisors defined by
the 1-dimensional ray of the cone $V$, the minimum is achieved by the lattice points contained in $\overline{V}$.

However we note $Y^{\alpha}$ is
already an infinte series, so unlikely to be a rational function.

The ``essential'' part of the Fourier expansion of the $h$ function can actually be written as a sum of such $Y^{\alpha}$:
$$
\frac{1}{D_FR_F}\sum_{\alpha \in V}\frac{1}{N(\alpha)}\sigma_{1,\chi}(\alpha)Y^{\alpha}
$$

Now disregard the coefficients $\frac{1}{D_FR_F}$, only looking at the sum 
$\sum_{\alpha \in V}\frac{1}{N(\alpha)}\sigma_{1,\chi}(\alpha)Y^{\alpha}$, if we ``specialize'' it to the normal cone
of the infinite divisors,i.e., we only take the lowest order term from the infinite series 
$Y^{\alpha}=\sum_{u \in U_F^+} X^{u\alpha}$, then we have:
$$
\sum_{\alpha \in V}\frac{1}{N(\alpha)}\sigma_{1,\chi}(\alpha)X^{\alpha} \sim \sum_{\alpha \in P(V)} 
log\prod_{n=1}^{\infty}(1-X^{n\alpha})
$$

where $P(V)$ is the set of primitive lattice points in $V$. 

In this way we can ``roughly'' see why the function $h(w;\fa)$ has a log-type singularities near the boundary divisors.

\begin{rmk}
We note the factor $D_FR_F$ can be interpretated as the volume of the ``link manifold''
of the toroidal boundary under some metric. This is natural as $D_F$ is the volume of $\bR^n/\cO_F$ while
$R_F$ is the volume of $\bR^{n-1}/U_F$ if we exponential the units to make it additive. This suggests that
if we can write the $h$ function as an integral over $((\bR^n/\cO_F) \oplus (\bR^+)^n)/U_F^+$, we might be able to
seperate the integral into the ``link'' part and ``polar'' part, which then might help us to understand the
nature of the $h(\w;\fa)$ functions better. 
\end{rmk}

\section{Further Discussions}

\subsection{}
As long as the Stark's conjecture concerns, the major problem now we are facing is the
integrality problem. Recall that we have shown that the quotient $\frac{\eta_K(R)}{\eta_K(\cO_K)}$
is an algebraic number, but in so far the analytic theory we developed, since we don't
have the explicit form about  $\eta_K(w;\fa)$, we don't have
the information about the integral properties about this number. Yet the Stark's conjecture
predits it is actually an unit in the class fields.

In some sense, as this quantity is such a canonical quantity, there is a chance that we can develop an algebraic
theory for them, such that the quotient $\frac{\eta_K(w_1;\fa_1)}{\eta_K(w_2;\fa_2)}$ is truely an
algebraic coordinate on the moduli space. By this we mean $\frac{\eta_K(R)}{\eta_K(\cO_K)}$ can ``tell''
the difference between $R$ and $\cO_K$, for example if $\wp$ is a prime ideal in $\cO_K$, which defines
the abelian variety $X_{\wp}$. The difference between $X_{\wp}$ and $X_0$ then is an ideal multiplication
by $\wp$, in other words, if we take reductions by primes other than $\wp$, $X_{\wp}$ and $X_0$ would be
isomorphic.
If $\frac{\eta_K(R)}{\eta_K(\cO_K)}$ is truely an algebraic coordinate, then 
$\frac{\eta_K(\wp)}{\eta_K(\cO_K)}$ should ``tell'' this difference, under any reduction by
any primes, i.e., it's absolute value's
prime factorization should only contain the primes from $\wp$. One is reminded Deuring's theory
of reduction of elliptic curves(\cite{Lang}). If we can show this, then the
integrality statement would follow immediately.

Of course so far these are just speculations.

If we can prove the Stark's conjecture in general, i.e., if we can show the quantities
$\frac{\eta_K(R)}{\eta_K(\cO_K)}$ are indeed units,
then the classical Iwasawa theory of cyclotomic fields and imaginery quadratic fields can be
generalized to arbitrary CM fields.

\subsection{}

In the classical theory of elliptic curves, Dedekind's
$\eta$ can be interpretated as the ``transcendental multipler'' of the period. Indeed for the
plane cubic curve
$$
y^2=4x^3 -g_2 x -g_3
$$
If we integrate the algebraic 1-form $\frac{dx}{y}$ over the 1-cycles, the classical elliptic integral
$$
\int_{-\infty}^{\infty} \frac{dx}{\sqrt{4x^3-g_2x-g_3}}
$$
gives the period as
$$
2 \pi\eta(\tau)^2 j(\tau)^{\frac{1}{12}}(12g_3)^{-\frac{1}{4}} (\bZ + \bZ\tau)
$$

This raises the question whether
our function $\eta_K(w,\fa)$ also has a similar interpretation. Since we regard $\eta_K(w;\fa)$ as
the coordinate function on the moduli space of primitive Hodge strucutres of middle degree, so the
natural choice is that we should expect this $\eta_K(w;\fa)$ as the transcendental multiplier of the
period of this primitive Hodge structure.

But then we have the following problem. Since we are dealing with the period integral over $X$, it's
natural to quotient off the automorphism group $X/(Aut(X))$. However we know in our case 
$Aut(X) \simeq U_F \simeq \bZ^{n-1} \oplus Torsion$, this is an infinite group when $ n>1$.
The quotient $X/(Aut(X))$ is not a variety, and in fact it is not even an algebraic space, it is just
a ``thing'', some kind of topological space. 
I am not sure how much algebraic sense can be made on them, in particular, I don't know
how to define an ``algebraic n-form'' on this quotient. 

On the other hand, as we explained before, the primitive classes of the middle degree can be regarded
as the equivariant classes under the action $Aut(X)$, so over there the Hodge structure is well-defined.
So although we don't know how to define the algebraic period, we know the quotient of the periods.
We just don't know the ``trancsendental multiplier'' for the Hodge structures.

Also in the elliptic curve theory, $\eta$ can be interpretated as the ``analytic torsion'' of the
flat metric. Can we do the same for the higher dimensional case?  This I believe is  possible,
for we can consider the space of invariant differential forms on $X$ under $U_F$, and since the flat metric is
invariant, it defines the Laplacian on the space of invariant forms. Over there we can define the
analytic torsion. Since the spectral zeta funtions of  invariant space can be identified (I believe) with zeta functions of
$F$, so our $\eta_K(w;\fa)$ can be identified with the analytic torsions in this equivariant sense.

\subsection{}
In his beautiful paper on the zeta functions of a totally real field, Shintani(\cite{Shintani1}) 
introduced the idea of using the multiple Gamma functions to represent the special zeta values.
The multiple Gamma function $\Gamma_n(x,w)$ is defined in a rather complicated manner(see a detailed
discussion in \cite{Yoshita}). Here we only note that if $w=(w_1,\cdots,w_n) \in (\bR^+)^n$, we
define the Riemann-Hurwitz zeta function as
$$
\zeta_n(s,w,x)=\sum(x+m_1w_1+m_2w_2+\cdots+m_nw_n)^{-s}
$$

where the sum is over all the non-negatives integers $(m_1,\cdots,m_n)$. Then $\Gamma_n(x,w)$ satisfying:
$$
log\frac{\Gamma_n(x,w)}{\rho_n(w)}=\frac{\partial}{\partial s}\zeta_n(s,w,x)|_{s=0}
$$
and
$$
-log\rho_n(w)=\lim_{x\ra +0}\{\frac{\partial}{\partial s}\zeta_n(s,w,x)|_{s=0} + logx \}
$$

In this sense Shintani's formula can be regarded as a generalization of the Lerch formula(\cite{Weil}):
$$
\frac{\partial}{\partial s}\zeta_1(s,w,x)|_{s=0}=log\frac{\Gamma(x)}{2\pi}
$$

Classically the product of Gamma values is related to the Dedekind $\eta$ by the Chowla-Selberg
formula (see \cite{Weil}). The same mechanism works equally well in the higher dimensional case,
multiple Gamma's related to our $\eta_K(w;\fa)$ in the same fashion(\cite{Yoshita}).

But Gamma function is also related to the Euler beta function
$$
B(a,b)=\int_0^1 t^{a-1} (1-t)^{b-1}dt = \frac{\Gamma(a)\Gamma(b)}{\Gamma(a+b)}
$$
which can be interpretated as some kind of ``period''. This is in fact the idea of Gross(\cite{Gross})
 that we can prove
the Chowla-Selberg formula geometrically, by explicitly calculate the period of Fermat curve
$x^n + y^n=1$.

Can we have the similar theory in the higher dimension? In view of the relation between multiple Gamma function
and our $\eta_K(w;\fa)$, this is the same problem as ``can we have a period interpretation of $\eta_K(w;\fa)$?''
I don't know any way of expressing multiple Gamma functions as the above simple integrals. In any case even
such interpretation exists, it seems not likely to be an ordinary period.

Also in his work on the real quadratic fields, Shintani(\cite{Shintani2}) discovered a formula that expressing
some combination of the multiple Gamma function as an infinite product. Precisely let $F$ be the real
quadratic field, $\epsilon$ be the fundamental unit, 
then the formula is
$$
\frac{\Gamma_2(z)}{\Gamma_2(1 + \epsilon -z)}=i^{1/2}exp(\frac{\pi}{12}i(\epsilon + \epsilon^{-1}))
\frac{\prod_{n=1}^{\infty}(1-q^nexp2\pi i z)}{\prod_{n=1}^{\infty}(1-(q^{\prime})^n exp\frac{2\pi i}{\epsilon}z)}
exp \frac{\pi i}{2}(\frac{z^2}{\epsilon} -(1+\epsilon^{-1})z)
$$
where $q=exp2 \pi i \epsilon$ and $q^{\prime}=exp(-2 \pi i \epsilon^{-1})$.

This beautiful formula is the main reason of my old belief that the function $\eta_K(w;\fa)$ should have
a similar infinite product formula as well. Now I no longer hold that belief. In some sense this formula
should be compared with the classical
$$
\Gamma(x) \Gamma(1-x)=\frac{\pi}{sin \pi x}=x^{-1}\prod_{n=1}^{\infty}(1-\frac{x^2}{n^2})^{-1}
$$ 

Let's note some peculiar feature of Shintani's infinite product, for the variable 
$q=exp2 \pi i \epsilon$ and $q^{\prime}=exp(-2 \pi i \epsilon^{-1})$, they are not inside the unit disc,
but somehow on the boundary, so strictly speaking it is a function of the upper half plane but
evaluated at the real axis. This raises the possibility that our function $\eta_K(w;\fa)$, when taking
the limit to the boundary, may also have an infinite product expression. Is this true or not? 
After all, it's unclear to me what's the meaning of this infinite product.

What is clear though is  much remains to be discovered.

\bigskip                    


\end{document}